%% file: lecture_ode.tex
\documentclass[12pt]{article}

\input{headlines}

\title{\large\sffamily\bfseries Lecture Notes \\ 
\LARGE{Numerical Solution of ODEs}}
\author{\Large\sffamily Davoud Mirzaei\\ Uppsala University}
\date{\sffamily February 16, 2023}

\begin{document}

\setlength{\abovedisplayskip}{3pt}
\setlength{\belowdisplayskip}{3pt}
\setlength{\abovedisplayshortskip}{0pt}
\setlength{\belowdisplayshortskip}{0pt}
\maketitle

\definecolor{contcol1}{HTML}{72E094}
\definecolor{contcol2}{HTML}{24E2D6}
\definecolor{convcol1}{HTML}{C0392B}
\definecolor{convcol2}{HTML}{8E44AD}

\begin{tcolorbox}[title=Contents, fonttitle=\huge\sffamily\bfseries\selectfont,interior style={left color=contcol1!10!white,right color=contcol2!10!white},frame style={left color=contcol1!30!white,right color=contcol2!30!white},coltitle=black,top=2mm,bottom=2mm,left=2mm,right=2mm,drop fuzzy shadow,enhanced,breakable]
\makeatletter
\@starttoc{toc}
\makeatother
\end{tcolorbox}

\vspace*{10mm}

\thispagestyle{empty}
\newpage
\pagenumbering{arabic}

\input{lec3_part1}
\input{lec3_part2}

\input{lec3_part3}

\input{lec3_part4}

\input{appendix}

\end{document}

%% file: headlines.tex
\usepackage[english]{babel}
\usepackage{graphicx}
\usepackage{framed}
\usepackage[normalem]{ulem}
\usepackage{indentfirst}
\usepackage{amsmath,amsthm,amssymb,amsfonts}
\usepackage[italicdiff]{physics}
\usepackage[T1]{fontenc}
\usepackage{mathtools,extarrows}

\usepackage{setspace}
\usepackage{wrapfig}
\usepackage{lmodern,mathrsfs}
\usepackage[inline,shortlabels]{enumitem}
\setlist{topsep=2pt,itemsep=2pt,parsep=3pt,partopsep=2pt}

\usepackage[dvipsnames]{xcolor}
\usepackage[utf8]{inputenc}
\usepackage[a4paper, top=0.8in,bottom=0.8in, left=0.8in, right=0.8in, footskip=0.3in, includefoot]{geometry}
\usepackage[most]{tcolorbox}
\usepackage{tikz,tikz-3dplot,tikz-cd,tkz-tab,tkz-euclide,pgf,pgfplots}
\pgfplotsset{compat=newest}
\usepackage{multicol}
\usepackage[bottom,multiple]{footmisc} 

\usepackage[backend=bibtex,style=numeric]{biblatex}
\usepackage{hyperref}

\usepackage{graphicx}
 \graphicspath{{figures/}} 

\numberwithin{equation}{section}

\usepackage{framed,color}
\definecolor{shadecolor}{rgb}{1,0.9,0.6}
\usepackage{caption}
\usepackage{tcolorbox}

\usepackage{multirow}
\usepackage{epstopdf}
\usepackage{listings}
\usepackage[]{mcode}
\usepackage{algorithm}
\usepackage{latexsym}
\usepackage{showidx}
\usepackage{latexsym}
\usepackage{amssymb}

\newcommand{\ignore}[1]{}
\usepackage{todonotes}

\usepackage{tikz}

\newcommand{\hide}[1]{} 

\newcommand{\myrightleftarrows}[1]{\mathrel{\substack{\xrightarrow{#1} \\[-.9ex] \xleftarrow{#1}}}}

\renewcommand{\qed}{\hfill\blacksquare}

\newcommand{\ep}{\varepsilon}

\newcommand{\C}{\mathbb{C}}

\newcommand{\R}{\mathbb{R}}

\newcommand{\wh}{\widehat}
\newcommand{\wt}{\widetilde}

\newcommand{\diagg}{\mathrm{diag}}

\def\bfm#1{\protect{\makebox{\boldmath $#1$}}}

\def\y {\bfm{y}}
\def\u {\bfm{u}}
\def\f {\bfm{f}}
\def\g {\bfm{g}}
\def\v {\bfm{v}}
\def\c {\bfm{c}}

\def\ee{\mathrm{e}}

\let\Re\relax

\DeclareMathOperator{\Re}{\text{Re}}

\newtheoremstyle{mystyle}{}{}{}{}{\sffamily\bfseries}{.}{ }{}
\newtheoremstyle{cstyle}{}{}{}{}{\sffamily\bfseries}{.}{ }{\thmnote{#3}}
\makeatletter
\renewenvironment{proof}[1][\proofname] {\par\pushQED{\qed}{\normalfont\sffamily\bfseries\topsep6\p@\@plus6\p@\relax #1\@addpunct{.} }}{\popQED\endtrivlist\@endpefalse}
\makeatother
\theoremstyle{mystyle}{\newtheorem{definition}{Definition}[section]}
\theoremstyle{mystyle}{\newtheorem{theorem}[definition]{Theorem}}
\theoremstyle{mystyle}{}
\theoremstyle{mystyle}{}
\theoremstyle{mystyle}{}
\theoremstyle{mystyle}{}
\theoremstyle{mystyle}{}
\theoremstyle{mystyle}{\newtheorem{workout}[definition]{Workout}}
\theoremstyle{mystyle}{\newtheorem{labexercise}[definition]{Lab Exercise}}
\theoremstyle{mystyle}{}

\usepackage{amsmath}
\newtheorem{remark}{Remark}[section]
\newtheorem{example}{Example}[section]

\theoremstyle{cstyle}{}

\newtheoremstyle{warn}{}{}{}{}{\normalfont}{}{ }{}
\theoremstyle{warn}

\newcommand{\warningsign}[1]{\tikz[scale=#1,every node/.style={transform shape}]{\draw[-,line width={#1*0.8mm},red,fill=yellow,rounded corners={#1*2.5mm}] (0,0)--(1,{-sqrt(3)})--(-1,{-sqrt(3)})--cycle;
\node at (0,-1) {\fontsize{48}{60}\selectfont\bfseries!};}}

\tcolorboxenvironment{definition}{boxrule=0pt,boxsep=0pt,colback={red!10},left=8pt,right=8pt,enhanced jigsaw, borderline west={2pt}{0pt}{red},sharp corners,before skip=10pt,after skip=10pt,breakable}
\tcolorboxenvironment{proposition}{boxrule=0pt,boxsep=0pt,colback={Orange!10},left=8pt,right=8pt,enhanced jigsaw, borderline west={2pt}{0pt}{Orange},sharp corners,before skip=10pt,after skip=10pt,breakable}
\tcolorboxenvironment{theorem}{boxrule=0pt,boxsep=0pt,colback={blue!10},left=8pt,right=8pt,enhanced jigsaw, borderline west={2pt}{0pt}{blue},sharp corners,before skip=10pt,after skip=10pt,breakable}
\tcolorboxenvironment{lemma}{boxrule=0pt,boxsep=0pt,colback={blue!10},left=8pt,right=8pt,enhanced jigsaw, borderline west={2pt}{0pt}{blue},sharp corners,before skip=10pt,after skip=10pt,breakable}
\tcolorboxenvironment{corollary}{boxrule=0pt,boxsep=0pt,colback={violet!10},left=8pt,right=8pt,enhanced jigsaw, borderline west={2pt}{0pt}{violet},sharp corners,before skip=10pt,after skip=10pt,breakable}
\tcolorboxenvironment{proof}{boxrule=0pt,boxsep=0pt,blanker,borderline west={2pt}{0pt}{CadetBlue!80!white},left=8pt,right=8pt,sharp corners,before skip=10pt,after skip=10pt,breakable}
\tcolorboxenvironment{remark}{boxrule=0pt,boxsep=0pt,blanker,borderline west={2pt}{0pt}{Green},left=8pt,right=8pt,before skip=10pt,after skip=10pt,breakable}
\tcolorboxenvironment{remarks}{boxrule=0pt,boxsep=0pt,blanker,borderline west={2pt}{0pt}{Green},left=8pt,right=8pt,before skip=10pt,after skip=10pt,breakable}
\tcolorboxenvironment{example}{boxrule=0pt,boxsep=0pt,blanker,borderline west={2pt}{0pt}{Black},left=8pt,right=8pt,sharp corners,before skip=10pt,after skip=10pt,breakable}
\tcolorboxenvironment{examples}{boxrule=0pt,boxsep=0pt,blanker,borderline west={2pt}{0pt}{Black},left=8pt,right=8pt,sharp corners,before skip=10pt,after skip=10pt,breakable}
\tcolorboxenvironment{cthm}{boxrule=0pt,boxsep=0pt,colback={gray!10},left=8pt,right=8pt,enhanced jigsaw, borderline west={2pt}{0pt}{gray},sharp corners,before skip=10pt,after skip=10pt,breakable}
\tcolorboxenvironment{workout}{boxrule=0pt,boxsep=0pt,colback={Cyan!10},left=8pt,right=8pt,enhanced jigsaw, borderline west={2pt}{0pt}{Cyan},sharp corners,before skip=10pt,after skip=10pt,breakable}
\tcolorboxenvironment{miniproject}{boxrule=0pt,boxsep=0pt,colback={gray!10},left=8pt,right=8pt,enhanced jigsaw, borderline west={2pt}{0pt}{gray},sharp corners,before skip=10pt,after skip=10pt,breakable}
\tcolorboxenvironment{labexercise}{boxrule=0pt,boxsep=0pt,colback={green!10},left=8pt,right=8pt,enhanced jigsaw, borderline west={2pt}{0pt}{green},sharp corners,before skip=10pt,after skip=10pt,breakable}

\newenvironment{talign*}{\let\displaystyle\textstyle\csname align*\endcsname}{\endalign}

\usepackage[explicit]{titlesec}
\titleformat{\section}{\fontsize{18}{30}\sffamily\bfseries}{\thesection}{18pt}{#1}
\titleformat{\subsection}{\fontsize{15}{18}\sffamily\bfseries}{\thesubsection}{15pt}{#1}
\titleformat{\subsubsection}{\fontsize{10}{12}\sffamily\large\bfseries}{\thesubsubsection}{8pt}{#1}

\titlespacing*{\section}{0pt}{5pt}{5pt}
\titlespacing*{\subsection}{0pt}{5pt}{5pt}
\titlespacing*{\subsubsection}{0pt}{5pt}{5pt}

\DeclareMathAlphabet\mathbfcal{OMS}{cmsy}{b}{n}
\newcommand{\vsp}{\vspace{0.3cm}}

\usepackage{matlab-prettifier}
\usepackage{textcomp,xcolor}
\definecolor{MatlabCellColour}{RGB}{252,251,220}
\lstdefinestyle{matlab}{
frame=single,
numbers=left,
basicstyle=\lstconsolas\small,
escapechar=`,
style=Matlab-editor,
backgroundcolor=\color{MatlabCellColour}
}
\lstdefinestyle{matlab1}{
frame=single,
numbers=right,
basicstyle=\lstconsolas\small,
escapechar=`,
style=Matlab-editor,
backgroundcolor=\color{MatlabCellColour}
}

%% file: lec3_part1.tex
Welcome to a {\em beautiful} subject in scientific computing: numerical solution of ordinary differential equations (ODEs) with initial conditions. More precisely, we consider
$$
y'(t) = f(t,y(t)), \quad y(t_0)=y_0, \quad t\in[t_0,b]
$$
where $t$ is an independent variable (usually plays the role of time), $y$ is the unknown solution of the problem (possibly a vector) which is sought, $f$ is a known function and $y_0$ is the initial condition at $t=t_0$.
We first give a brief introduction to the theory of ODEs and existence and uniqueness of solutions. Then some standard numerical techniques are derived and their advantages and disadvantages for solving different differential equations are outlined. In some parts of this lecture we follow
\cite{Atkinson-et-al:2009} and \cite{LeVeque:2007}.

\section{An introduction to ODEs}
\begin{figure}[!th]
\centering
\includegraphics[scale=0.75]{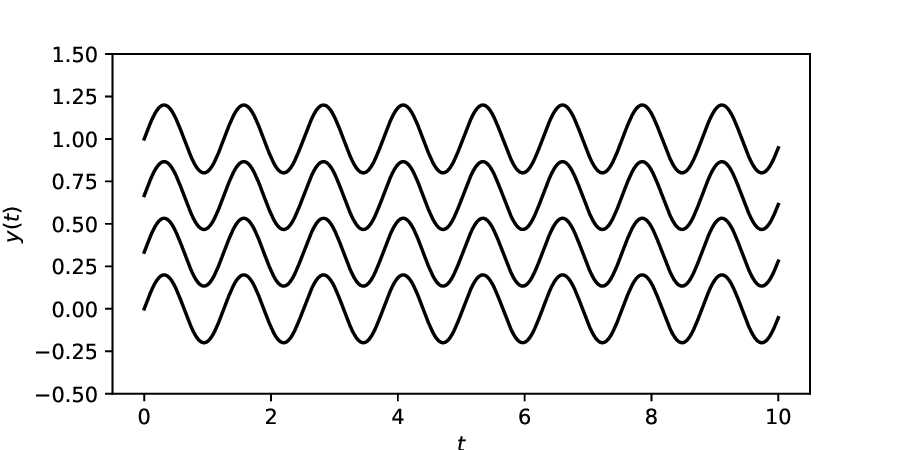}
\caption{Integration curves}\label{fig:ode0}
\end{figure}

The simplest ordinary differential equation (ODE) has the form
\begin{equation}\label{ode:simplest}
y'(t) = g(t)
\end{equation}
for a given function $g$. The {\em general solution} of this equations is
$$
y(t) = \int g(\tau) d\tau + c
$$
with $c$ an arbitrary integration constant. Here $\int f(\tau)d\tau$ denotes any fixed antiderivative of $g$.
In Figure \ref{fig:ode0}, with the case of $g(t)=\cos 5t$, the plots of $y(t)$ for four different values of $c$ are shown.
Such plots are sometimes called {\em integration curves}. 
For the simple ODE \eqref{ode:simplest}, the integration curves are just copies (shifts) of each other along the $y$-axis.

The constant $c$ can be obtained by specifying the value of $y$ at some given point
$$
y(t_0) = y_0,
$$
where $y_0$ is known.
The {\em particular solution} of the differential equation then can be written as
$$
y(t) = y_0 + \int_{t_0}^{t}g(\tau)d\tau,
$$
provided that $g$ be (for example) a continuous function. Another simple differential equation is
$$
y'(t) = \lambda y(t),\quad y(t_0) = y_0, \quad \lambda\in \R
$$
which possesses the exponential solution
$$
y(t) = y_0 e^{\lambda (t-t_0)}.
$$
A more general form which contains both the above differential equations reads as
\begin{equation}\label{ivp:linear1}
y'(t) = \lambda y(t)+g(t), \quad y(t_0) = y_0.
\end{equation}
The general solution of this equation can be obtained by the so-called {\em method of integration factors}. Multiplying the above linear equation by integration factor $e^{-\lambda t}$, the equation is reformulated as
$$
\frac{d}{dt}\big(e^{-\lambda t}y(t)\big)=e^{-\lambda t}g(t).
\vsp
$$
Integrating both sides from $t_0$ to $t$, we obtain
$$
y(t)= e^{\lambda t}\left[ c + \int_{t_0}^{t}e^{-\lambda \tau}g(\tau) d\tau  \right].
$$
Imposing the condition $y(t_0)=y_0$ gives $c = \exp(-\lambda t_0)y_0$, and therefore 
\begin{equation}\label{ivp:linsol}
y(t) = y_0e^{\lambda (t-t_0)} + \int_{t_0}^{t}e^{\lambda (t-\tau)}g(\tau)d\tau.
\end{equation}
In many applications of differential equations the independent variable $t$ plays the role of time, and $t_0$ can be interpreted as the initial time, and $y(t_0)=y_0$ is referred to the {\em initial condition}. A differential equation of the above type is called an
{\bf initial value problem (IVP)}. In a more general form an IVP has the form
\begin{equation}\label{ivp:form}
  \begin{split}
     \y'(t) =& \,\f(t,\y(t)) \\
      \y(t_0) =&\, \y_0
  \end{split}
\end{equation}
which is a system of IVPs for vector $\y = [y_1,\ldots,y_n]^T$. Here $\f$ may represent a nonlinear relation between the independent variable $t$ and the dependent variable $\y$.

For the linear IVP \eqref{ivp:linear1} the solution is given by \eqref{ivp:linsol} and it exists and is unique on any open interval where the data function $f$ is continuous. But for the nonlinear IVP \eqref{ivp:form} even if the right hand side function $f(t,y)$
has derivatives of any order the solution of IVP may exist on only a smaller interval. In some cases the solution is not unique.
\vsp 
\begin{example}
Consider the nonlinear equation
$$
y'(t) = -y(t)^2,\quad t\geqslant 0.
$$
This problem has the trivial solution $y(t)\equiv 0$ and a general solution
$$
y(t) = \frac{1}{t+c}
$$
with arbitrary constant $c$. Let the equation be accompanied by initial condition $y(0)=y_0$. If $y_0=0$ then $y(t)\equiv 0$ is the solution of the IVP for any $t\geqslant 0$. If $y_0\neq 0$ then the solution of IVP is
$$
y(t)=\frac{1}{t+y_0^{-1}}.
\vsp
$$
For $y_0>0$ the solution exists for any $t\geqslant0$ while for $y_0<0$ the solution exists only on interval $[0,-y_0^{-1}]$.

\end{example}
\vsp 

\begin{example}
Consider the IVP
$$
y'(t) = 2\sqrt{y(t)}, \quad t\geqslant 0, \quad y(0)=0.
$$
It is clear that both $y(t)\equiv 0$ and $y(t) = t^2$ are solutions of this IVP. In additions any $C^2$ function $y(\cdot;\alpha)$ of the form
$$
y(t;\alpha) = (t-\alpha)_{+}^2 = \begin{cases}
0, & 0\leqslant t\leqslant \alpha\\
(t-\alpha)^2, & t>\alpha
\end{cases}
$$
for any $\alpha \geqslant0$ is a solution for this IVP. In Figure \ref{fig:ode1} the solution with $\alpha=0,1,2$ are plotted. This example reveals the non-uniqueness of the nonlinear IVP
\eqref{ivp:form} for some right hand side functions $f$.
\begin{center}
  \includegraphics[scale=0.7]{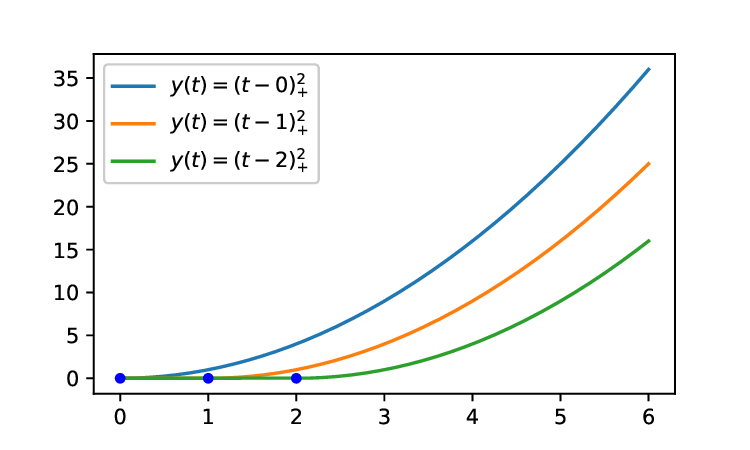}
  \captionof{figure}{Three sample solutions for IVP $y'=2\sqrt y$ with $y(0)=0$ (non-uniqueness).}\label{fig:ode1}
\end{center}
\end{example}
\vsp 

To guarantee that there is a unique solution it is
necessary to impose a certain amount of smoothness on function $f(t,y)$.
The following well-known theorem establishes the existence and uniqueness of the IVP \eqref{ivp:form}.
\begin{theorem}\label{thm:ivpexuq}
Let $D$ be an open and connected set in $\R^2$ and $f(t,y)$ be a continuous function in both $t$ and $y$ in $D$, and let $(t_0,y_0)$
be an interior point in $D$. Assume that $f$ satisfies the Lipschitz continuity in its second argument, i.e., there exists a constant
$L\geqslant 0$ such that
\begin{equation*}
  |f(t,y)-f(t,\tilde y)|\leqslant L|y-\tilde y|,\quad \forall (t,y),\,(t,\tilde y)\in D.
  \vsp
\end{equation*}
Then there exists a unique function $y(t)$ defined on an interval $[t_0-\beta,t_0+\beta]$ for some $\beta>0$ satisfying
$$
y'(t) = f(t,y(t)), \quad t_0-\beta\leqslant t\leqslant t_0+\beta, \quad y(t_0)=y_0.
$$
\end{theorem}
The Lipschitz continuity is slightly stronger than mere continuity, which only
requires that $|f(t,y)-f(t,\tilde y)|\to 0$ as $\tilde y\to y$. Lipschitz continuity requires that
$$
|f(t,y)-f(t,\tilde y)| = \mathcal{O}(|y-\tilde y|), \; {as}\; \tilde y\to y.
$$
If $f$ is differentiable with respect to $y$ in $D$ and this derivative $\frac{\partial f}{\partial y}(t,y)$ is
bounded then we can take
$$
L = \max_{(t,y)\in \overline D}\left|\frac{\partial f}{\partial y}(t,y)\right|
\vsp
$$
because the Taylor series representation gives
$$
f(y,t)=f(t,\tilde y) + (y-\tilde y)\frac{\partial f}{\partial y}(t,\eta),\; \mbox{for some}\, (t,\eta)\in D.
$$
The number $\beta$ in the statement of the theorem depends on the IVP \eqref{ivp:form}. For some equations, solutions exist for any  $t$, thus we can take $\beta$ to be infinity. However for many nonlinear equations solutions can exist only in a bounded interval.
\vsp 

\begin{example}
For IVP \eqref{ivp:linear1} we have $f(t,y)=\lambda y + g(t)$; hence $L=|\lambda|$.
This problem of course has a unique solution for any initial $y_0$ and for any $t$.
In particular, if $\lambda=0$ then $f(t,y)=g(t)$. In this case $f$ is independent of $y$.
The solution is then obtained by simply integrating the function $g(t)$.
\end{example}
\vsp 

\begin{example}
Consider the IVP
$$
y'(t)=2ty(t)^2,\quad y(0)=1.
$$
For this equation $f(t,y)=2ty^2$ and $\frac{\partial f}{\partial y}(t,y)=4ty$. Both functions are continuous for all $(t,y)$. On any bounded domain $D=(a,b)\times(c,d)$ containing $(t_0,y_0)=(0,1)$ we can take $L=4bd$. According to Theorem \ref{thm:ivpexuq} there is a unique solution to this IVP for $t$ in some neighborhood of $t_0=0$. This solution is
$$
y(t)=\frac{1}{1-t^2}, \quad -1<t<1.
$$
As a side note, this example shows that the continuity of $f(t,y)$ and $\frac{\partial f}{\partial y}(t,y)$ for all $(t,y)$ does not imply the existence of a solution $y(t)$ for all $t$.
\end{example}
\vsp

\subsection{Modelling with ODEs: a funny example}
Imagine that you are jogging along a given path. Suddenly a dog in a nearby garden sees you and begins chasing you at full speed with constant velocity $w$. What is the
trajectory of the dog if we assume it is always running directly toward you?\footnote{This example is taken form: 
W. Gander, M. J. Gander, F. Kwok, Scientific Computing, An Introduction using Maple and MATLAB, Springer (2014).}

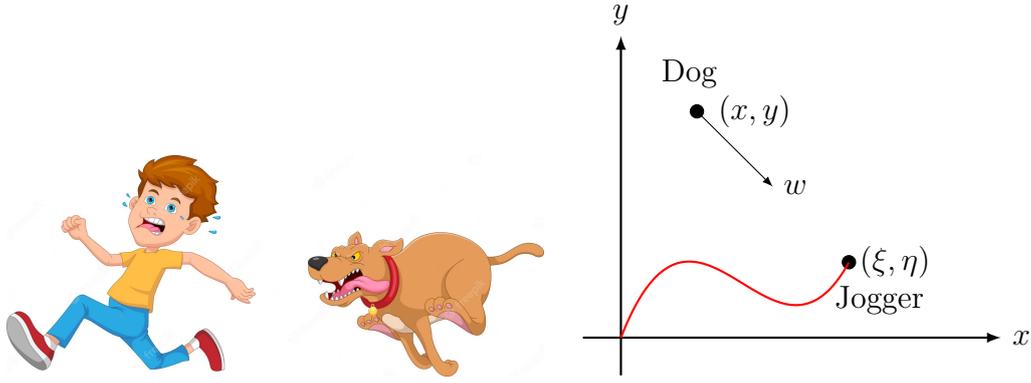
\begin{figure}[ht!]
\centering
\includegraphics[scale=.1]{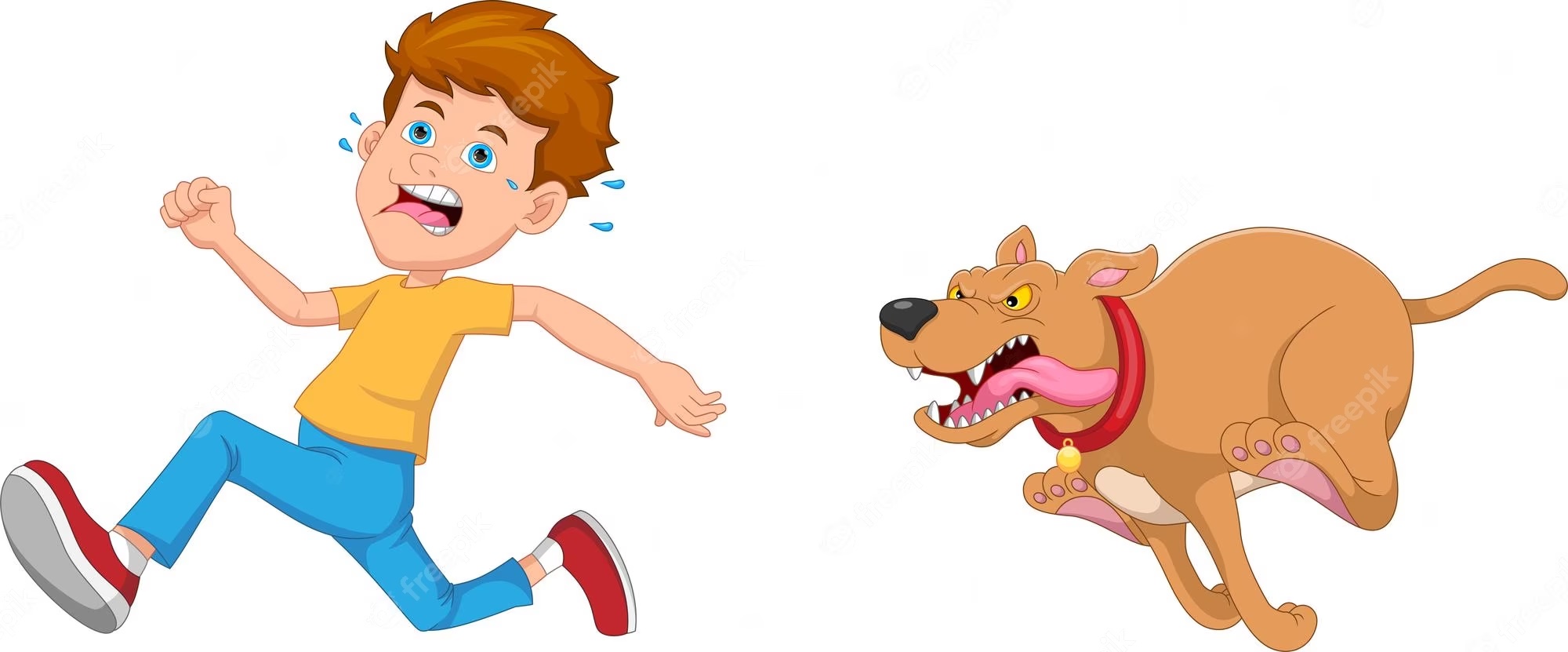} \;\; 
\begin{tikzpicture}[>=latex,scale=1]
\draw[->,thick] (-.5,0) -- (5,0) node[right] {$x$};
\draw[->,thick] (0,-.5) -- (0,4) node[above] {$y$};
\filldraw  (3,1) circle (2.5pt) node[right] {$(\xi,\eta)$};
\filldraw  (1,3) circle (2.5pt) node[right] {$~(x,y)$};
\draw[->] (1,3) -- (2,2) node[right] {$w$};
\node[color=black] at (.9,3.5) {Dog};
\node[color=black] at (3.4,.5) {Jogger};
\draw[samples=200, domain=0:3, color=red,thick] plot (\x,{\x-3/4*\x*(\x-1)+5/12*\x*(\x-1)*(\x-2)});
\end{tikzpicture}
\caption{Dog chasing a jogger. (left image from freepik.com)}\label{fig:dog}
\end{figure}

This situation is depicted in Figure \ref{fig:dog}. We assume that the trajectory of you is represented by $(\xi(t),\eta(t))$ and the trajectory of dog by  $(x(t),y(t))$. 
Since the dog is running with constant speed $w$, we have
\begin{equation}\label{dog-speed}
[x'(t)]^2 + [y'(t)]^2 = w^2, \quad \forall t\geqslant 0.
\end{equation}
Since the dog is always running toward you,
the velocity vector of the dog is proportional to the difference vector between the
position of you and the dog, i.e., 
\begin{equation}\label{dog-velocity}
\begin{bmatrix}
x'(t)\\ y'(t)
\end{bmatrix} = 
\lambda(t) 
\begin{bmatrix}
\xi(t)- x(t)\\ \eta(t)-y(t)
\end{bmatrix} ,
 \quad \forall t\geqslant 0, 
\end{equation}
where $\lambda(t)>0$ is a constant of proportionality. To find $\lambda(t)$ we can substitute \eqref{dog-velocity} in to \eqref{dog-speed} to obtain 
$$
\lambda^2 =  \frac{w^2}{(\xi-x)^2 + (\eta-y)^2}. \vsp
$$
The trajectory of the dog therefore satisfies the following system of ODEs
\begin{equation}\label{dog-ode}
\begin{bmatrix}
x'(t)\\ y'(t)
\end{bmatrix} = 
\frac{w}{\sqrt{(\xi(t)-x(t))^2 + (\eta(t)-y(t))^2}}
\begin{bmatrix}
\xi(t)- x(t)\\ \eta(t)-y(t).
\end{bmatrix}  
\end{equation}
The initial condition for this ODE is the initial position of the dog, i.e.,
$(x(0),y(0))=(x_0,y_0)$.

To find the trajectory of the dog, it remains to solve this IVP. But there is no hope for finding a closed-form solution for a general jogging path $(\xi(t),\eta(t))$. Soon we will solve this equation using some numerical methods, but let see the dog trajectories for different choices of jogger's path in Figure \ref{trajdog_fig}. This results are obtained by the Euler's method. See section \ref{sect-numerical-methods} below. 

\begin{figure}[ht!]
\centering
\includegraphics[scale=.5]{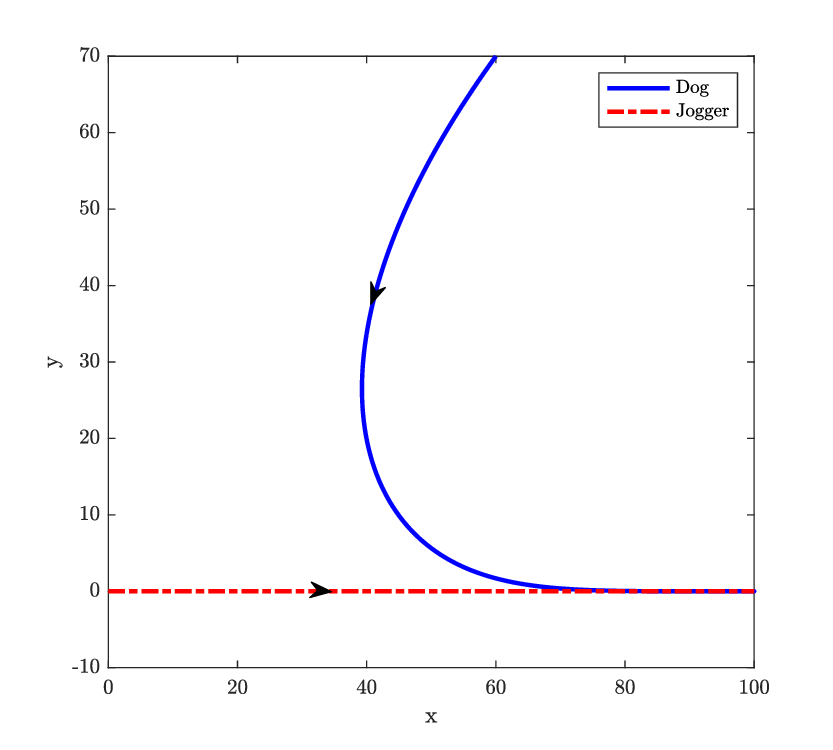} \includegraphics[scale=.5]{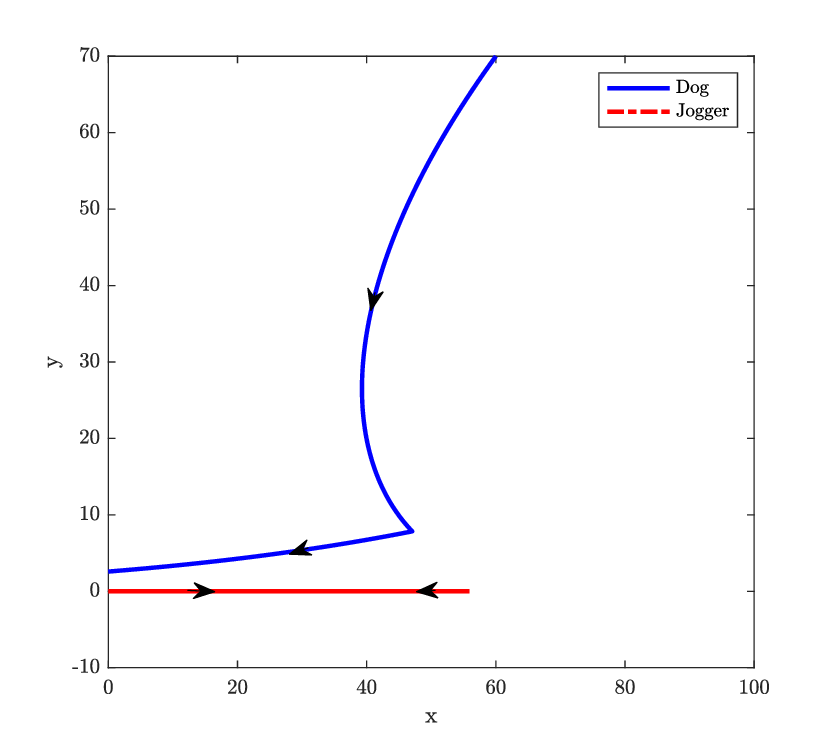}\\
\includegraphics[scale=.5]{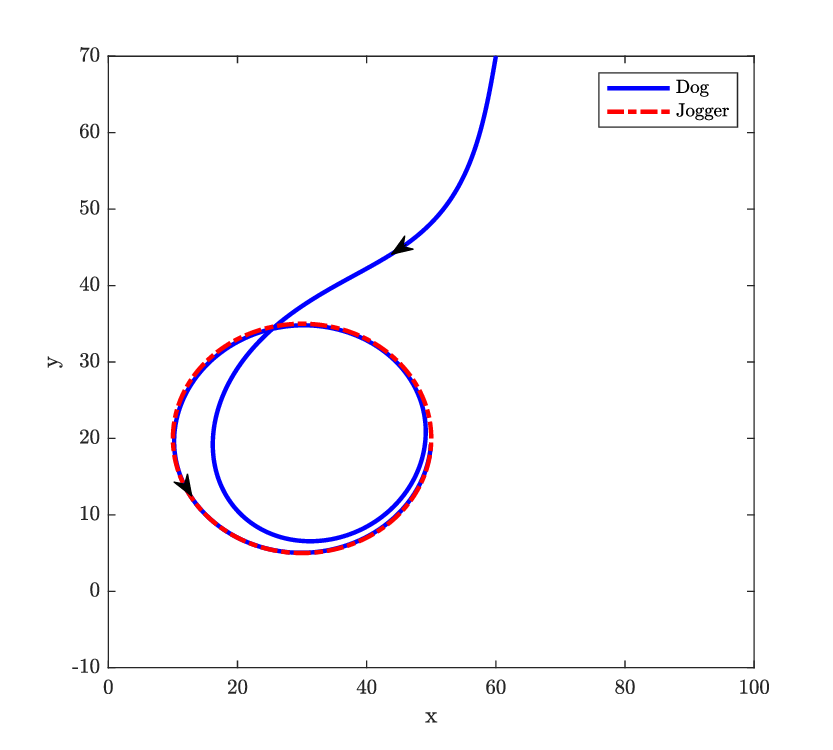} \includegraphics[scale=.5]{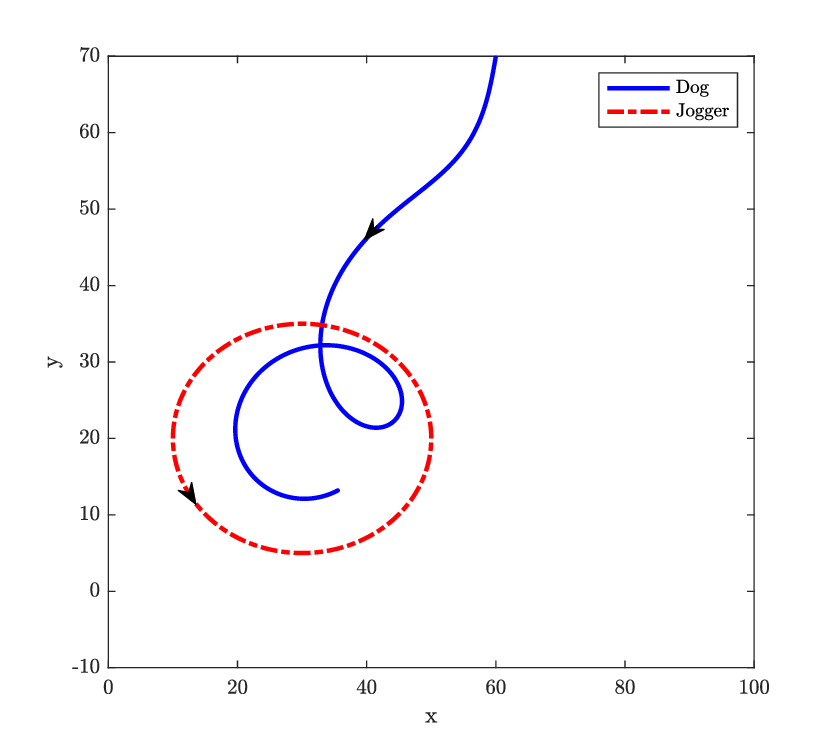}
\caption{Jogging path and the numerically computed
trajectory of the dog: The jogger running in a straight path (top-left), the jogger notices the dog and tries to run back (top-right), the jogger running
on a circular track (bottom-left), the jogger running
on a circular track but the dog is slow (bottom-right).}\label{trajdog_fig}
\end{figure}
\vsp

\subsection{Higher order ODEs}
The highest-order derivative appearing in an ODE determines the order of the
ODE. A $k$-th order ODE in the most general form can be written as
\begin{equation*}
  f(t,y,y',\ldots,y^{(k)}) = 0,
\end{equation*}
where $f:\R^{k+2}\to \R$ is a known function and  $y(t)$ is to be determined.
A $k$-th order ODE is said to be explicit if it can be written in the form
\begin{equation}\label{ode:kthorder}
  y^{(k)}(t) = f(t,y,y',\ldots,y^{(k-1)}).
\end{equation}
As a necessary condition for a unique solution, this ODE should accompany with $k$ initial conditions
\begin{equation}\label{ode:kthorder_initial}
y(t_0)=y_0,\quad y'(t_0)=y'_0, \quad \ldots\quad  y^{(k-1)}{(t_0)}=y^{(k-1)}_{0}.
\end{equation}
Many ODEs arise naturally of the form \eqref{ode:kthorder}, and many others can be transformed into it.
So far we considered the case $k=1$.
This is not a real restriction because
a higher-order ODE can always be reduced to a system of first-order equations
as follows.
For the explicit $k$-th order ODE \eqref{ode:kthorder} define $k$ new
variables
$$
y_1(t) = y(t),\quad y_2(t) = y'(t),\quad \ldots \quad y_k(t) = y^{(k-1)}(t),
$$
so that the original
$k$-th order equation becomes a system of $k$ first-order equations
$$
\y'(t) = \begin{bmatrix}
          y'_1(t) \\
          y'_2(t) \\
          \vdots \\
          y'_{k-1}(t) \\
          y'_{k}(t)
        \end{bmatrix}
 = \begin{bmatrix}
          y_2(t) \\
          y_3(t) \\
          \vdots \\
          y_{k}(t) \\
          f(t,y_1,y_2,\ldots,y_k)
        \end{bmatrix}=: \g(t,\y),
$$
with initial condition
$$
y_1(t_0)=y_0,\quad y_2(t_0)=y'_0, \quad \ldots\quad  y_k{(t_0)}=y^{(k-1)}_{0}.
$$
or simply
$$
\y(t_0) = \y_0,
\vsp 
$$
for $\y_0=[y_0,y'_0,\ldots,y_0^{(k-1)}]$. In general, we will focus on explicit first-order system of ODEs with initial conditions of the form
\begin{equation}\label{ode:sysform}
\begin{split}
   \y'(t) =&\,  \f(t,\y) \\
   \y(t_0) = &\, \y_0
\end{split}
\end{equation}
where $\f:\R^{n+1}\to \R^{n}$ and $\y_0\in \R^n$. If $\f$ is not explicitly depend on $t$, i.e., $\f(t,\y)=\f(\y)$, the system is called
{\em autonomous} and can be written in the form
$$\y' = \f(y).$$
 A nonautonomous ODE $\y' = \f(t, \y)$ can
always be converted to autonomous form by introducing an additional dependent
variable $y_{n+1}(t) = t$, so that $y'_{n+1}(t)=1$ and $y_{n+1}(t_0)=t_0$ yielding the autonomous ODE
$$
\begin{bmatrix}
          \y'(t) \\
          y'_{n+1}(t) \\
        \end{bmatrix}
 = \begin{bmatrix}
          \f(y_{n+1},\y) \\
          1 \\
        \end{bmatrix}.
$$
\ \\
It
is often convenient to assume $\f$ is of this form since it simplifies notation.
\vsp
\begin{example}
Consider the second order ODE
$$
\theta''(t) = -\frac{g}{\ell}\sin(\theta(t)),
$$
which models the motion of a pendulum with mass $m$ at the end of a rigid
bar of length $\ell$ by ignoring the mass of the bar and forces of friction and air resistance.
Here $\theta(t)$ is the angle of the pendulum from vertical at time $t$, and $g$ is the gravitational constant.
The motion is independent of the mass of pendulum.
Let $v = \theta'$ be the velocity and define
$$
\y = \begin{bmatrix}
  \theta \\
  v
\end{bmatrix}
$$
to obtain the following first-order linear system of equations
$$
\begin{bmatrix}
  \theta' \\
  v'
\end{bmatrix} =
\begin{bmatrix}
  v\\
  -(g/\ell)\sin\theta
\end{bmatrix}=:
\begin{bmatrix}
  f_1(\theta,v)\\
  f_2(\theta,v)
\end{bmatrix}.
$$

\end{example}
\vsp 
\begin{workout}
Convert the following system of third order equations to a system of first order equations:
\begin{align*}
  & u'''(t)+4u''(t)+5u'(t)+2u(t)=2t^2+10t+8 \\
  & u(0)=1, \; u'(0)=-1,\; u''(0)=3. 
\end{align*}
\end{workout}
\vsp 

\begin{workout}
The following system of second order equations arises from studying the gravitational attraction of one mass by another. Convert it to a system of first order equations.
\begin{align*}
x''(t)=-\frac{cx(t)}{r(t)^3}, \quad y''(t)=-\frac{cy(t)}{r(t)^3}, \quad z''(t)=-\frac{cz(t)}{r(t)^3}, \quad
\end{align*}
Here $c$ is a positive constant and $r(t)=\sqrt{x(t)^2+y(t)^2+z(t)^2}$ with $t$ denoting time. 
\end{workout}
\vsp 

\subsection{First order system of ODEs}
The study of first order system of ODEs is essentially important not only in solving a high order scaler equation using the techniques for first order ODEs but also in variety of other applications in which the system is obtained directly from the problem model.
Such systems also appear in the procedure of solving parabolic and hyperbolic partial differential equations using the method of lines (MOL).

The system of ODEs \eqref{ode:sysform} is linear if
$$
\f(t,\y)=A(t)\y + \g(t)
$$
where $A(t)\in \R^{n\times n}$ and $\g\in \R^n$. An important special case is the constant coefficient
linear system
\begin{equation*}
  \y'(t) = A\y(t) + \g(t)
  \vsp 
\end{equation*}
where $A \in\R^{n\times n}$ is a constant matrix. If $\g(t)=0$, then the equation is homogeneous. The
solution to the homogeneous system $\y'(t)=A\y(t)$ with initial data $\y(t_0)=\y_0$ is
\begin{equation*}
  \y(t) = e^{A(t-t_0)}\y_0.
\end{equation*}
where $e^{A(t-t_0)}$ is matrix exponential. 
\vsp

\begin{example}[Chemical Reaction Kinetics]
Let $X$ and $Y$ represent chemical compounds and consider a reaction of the form
$$
X \xlongrightarrow{k_1}Y
$$
which represents a reaction in which $X$ is transformed into $Y$ with rate $k_1 > 0$. If we let $y_1$
represent the concentration of $X$ and $y_2$ represent the concentration of $Y$ (often denoted by
$y_1=[X]$ and $y_2=[Y]$, then the ODEs for $y_1$ and $y_2$ are
\begin{align*}
  y'_1 =& \,-k_1y_1 \\
  y'_2 =& \,+k_1 y_1
\end{align*}
If there is also a reverse reaction at rate $k_2$, we write
$$
X \myrightleftarrows{\rule{0.8cm}{0cm}}_{{}_{{}_{\hspace{-.7cm}k_2}}}^{{}^{\hspace{-.7cm}k_1}} Y
$$
and we have the system of ODEs
\begin{align*}
  y'_1 =& \,-k_1y_1+k_2y_2 \\
  y'_2 =& \,+k_1 y_1-k_2y_2
\end{align*}
which with given initial concentrations $y_1(0)$ and $y_2(0)$ forms a linear system of initial value problems with constant coefficient matrix:
\begin{equation}\label{chem:eq1}
\begin{bmatrix} y'_1\\y'_2  \end{bmatrix}=
\begin{bmatrix}
  -k_1 & k_2  \\
  k_1 & -k_2
\end{bmatrix}
\begin{bmatrix} y_1\\y_2  \end{bmatrix}, \quad \begin{bmatrix} y_1(0)\\y_2(0)  \end{bmatrix}=
\begin{bmatrix} y_{1,0}\\y_{2,0}  \end{bmatrix}.
\end{equation}
\ \\
Another simple system arises from the decay process
$$
X \xlongrightarrow{k_1}Y\xlongrightarrow{k_2}Z.
$$
If $y_1=[X]$, $y_2=[Y]$ and $y_3=[Z]$ then we have the following equations
\begin{equation}\label{chem:eq2}
\begin{array}{rl}
  y'_1 =& \,-k_1y_1 \\
  y'_2 =& \,k_1 y_1 - k_2y_2\\
  y'_3 = &\, k_2y_2
\end{array}, \quad or \quad
\begin{bmatrix} y'_1\\y'_2\\y'_3  \end{bmatrix}=
\begin{bmatrix}
  -k_1 & 0 & 0 \\
  k_1 & -k_2 & 0 \\
  0 & k_2 & 0
\end{bmatrix}
\begin{bmatrix} y_1\\y_2\\y_3  \end{bmatrix}.
\end{equation}
\end{example}
\vsp 

Consider a linear system $\y'=A\y$ with initial condition $\y(t_0)=\y_0$, where $A$ is a constant $n\times n$ matrix, and suppose for
simplicity that $A$ is diagonalizable, i.e., $A$ has a complete set of $n$ linearly independent eigenvectors $\v_k$ corresponding to eigenvalues $\lambda_k$ for $k=1,\ldots,n$ such that
$$
A\v_k = \lambda_k \v_k, \quad k=1,\ldots,n
$$
or equivalently
$$
A = VDV^{-1},
$$
where $V=[\v_1\; \v_2\ldots \v_n]$ and $D = \diagg(\lambda_1,\ldots,\lambda_n)$. Applying the change of variables $\u(t) = V^{-1}\y(t)$, the linear system $\y'=A\y$ will transfer to
$$
\u' = D\u, \quad \u(t_0) = V^{-1}\y_0=:\u_0.
$$
This is a decoupled system of ODEs because $D$ is diagonal. We may write
\begin{equation}\label{sys:solveAy}
u'_k = \lambda_k u_k, \quad u_k(t_0)=u_{k,0}, \quad k=1,2,\ldots,n,
\end{equation}
with solutions $u_k(t)=u_{k,0}e^{\lambda_k (t-t_0)}$ or in a matrix form
$$
\u(t) = e^{D(t-t_0)}\u_0
$$
or equivalently 
$$ V^{-1}\y(t) = e^{D(t-t_0)}V^{-1}\y_0 $$
which finally gives
$$
 \y(t)=Ve^{D(t-t_0)}V^{-1}\y_0 = e^{A(t-t_0)}\y_0.
$$
Keep in mind that $e^{A(t-t_0)}$ is a matrix and $\y_0$ is a vector, thus $e^{A(t-t_0)}\y_0$ is a matrix-vector multiplication. 
\ \\
\begin{example}\label{ex:chem_solve}
Consider the linear system of ODEs \eqref{chem:eq1} with $k_1=2$ and $k_2=1$ and 
initial conditions $y_1(0)=5$ and $y_2(0)=2$: 
\begin{equation*}
\begin{bmatrix} y'_1\\y'_2  \end{bmatrix}=
\begin{bmatrix}
  -2 & 1  \\
  2 & -1
\end{bmatrix}
\begin{bmatrix} y_1\\y_2  \end{bmatrix}, \quad \begin{bmatrix} y_1(0)\\y_2(0)  \end{bmatrix}=
\begin{bmatrix} 5\\2  \end{bmatrix}.
\end{equation*}
It can be simply shown that the eigenvalues and eigenvectors of $A$ are
\begin{align*}
  \lambda_1 &= -3, \quad \v_1=[1,-1]^T\\
  \lambda_2 &=0, \qquad \v_2=[1,2]^T,
\end{align*}
which gives 
$$
V = \begin{bmatrix}
  1 & 1  \\
  -1 & 2
\end{bmatrix}, \quad 
D = \begin{bmatrix}
  -3 & 0  \\
  0 & 0
\end{bmatrix}.\vsp
$$
Therefore, the solution can be written as
$$
\y(t) = Ve^{Dt}V^{-1}\y_0 =\frac{1}{3} \begin{bmatrix}
                            1 & 1 \\
                            -1 & 2
                          \end{bmatrix}
\begin{bmatrix}
  e^{-3t} & 0 \\
  0 & e^{0t}
\end{bmatrix}
\begin{bmatrix}
  2 & -1 \\
  1 & 1
\end{bmatrix}
\begin{bmatrix} 5 \\ 2  \end{bmatrix}
$$
which shows that
\begin{align*}
  y_1(t) & \,= \frac{5}{3}\left[2e^{-3t}+1\right]+ \frac{2}{3}\left[-e^{-3t}+1\right]\\
  y_2(t) & \,= \frac{5}{3}\left[-2e^{-3t}+2\right]+ \frac{2}{3}\left[e^{-3t}+2\right].
\end{align*}
The plots of solutions are depicted in Figure
\ref{fig:ode2}.
\begin{center}
  \includegraphics[scale=0.75]{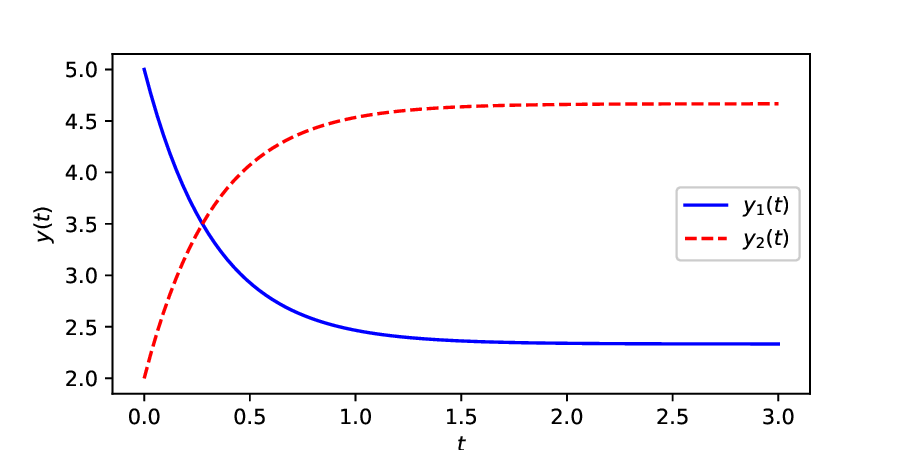}
  \captionof{figure}{Solutions of a chemical reaction kinetics problem.}\label{fig:ode2}
\end{center}
Since, in this example, the reaction rate $k_1$ is larger than $k_2$, the amount $y_1$ decreases while $y_2$ increases. Both solutions reach the steady states $y_1(t)\to \frac{1}{3}(y_1(0)+y_2(0))=\frac{7}{3}$ and $y_2(t)\to \frac{2}{3}(y_1(0)+y_2(0))=\frac{14}{3}$ when $t$ increases.
\end{example}
\vsp 

\begin{workout}
In the chemical reaction model \eqref{chem:eq2} let $k_1=2$ and $k_2=1$ and
$[y_1(0), y_2(0), y_3(0)] = [1, 3, 2]$. Obtain the solution. 

Hint: Calculate the eigenvalues and eigenvectors of the coefficient matrix and put them into the formula.  

\end{workout}
\vsp 

The situation will be more complicated for a nonlinear system of equations. 
In the sequel we study several numerical algorithms for solving different types of ODEs; both linear and nonlinear equations. 

\subsection{Stability of solutions}\label{sect-stabsol}
Depending on given data $f$ and $y_0$, solutions of IVP may behave differently as $t\to\infty$.
The Lipschitz constant measures how much $f(t,y)$ is changed if $y$ is perturbed (at some fixed
time $t$). Since $f(t,y)=y'(t)$ is the slope of the line tangent to the solution curve through
the value $y$, this indicates how the slope of the solution curve is varied if $y$ is perturbed.
\vsp

\begin{example}
The solutions of $y'(t)=g(t)$ with Lipschitz constant $L=0$ are parallel curves each for a prescribed initial condition $y_0$. See Figure \ref{fig:ode0}.
The ODE $y'(t)=\lambda y(t)$ with $L=|\lambda|$ possesses solutions $y(t)=y_0e^{\lambda t}$. Depending on the sign of $\lambda$
solutions decay exponentially to zero (if $\lambda <0$), grow exponentially to infinity (if $\lambda >0$), or stay
parallel lines (if $\lambda=0$) for different values of $y_0$. See Figure \ref{fig:ode345}.
\begin{center}
  \includegraphics[scale=0.67]{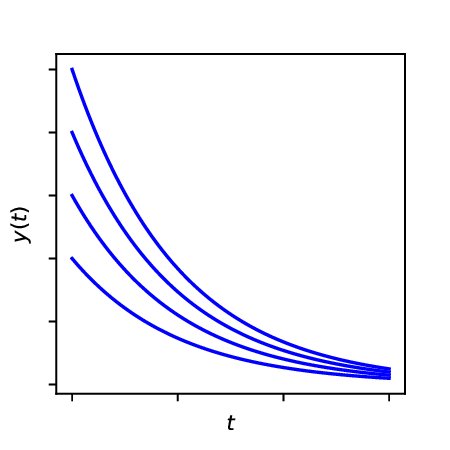}\includegraphics[scale=0.67]{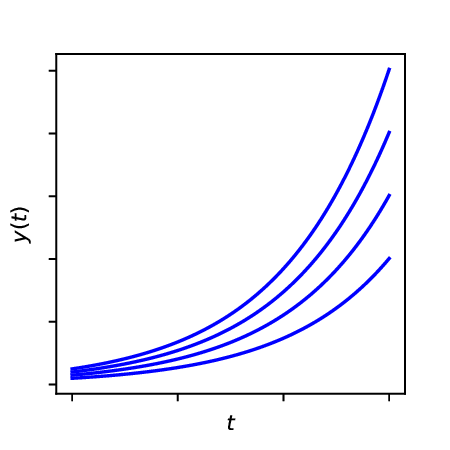}\includegraphics[scale=0.67]{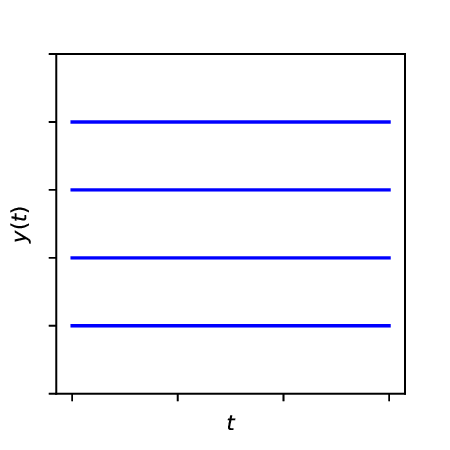}
  \captionof{figure}{Solutions of IVP $y'=\lambda y$ with different initial values $y_0$ for $\lambda<0$ (left), $\lambda>0$ (middle), and $\lambda=0$ (right).}\label{fig:ode345}
\end{center}

For the case of $\lambda>0$
any two solutions diverge away from each other, i.e.,
a small perturbation in the initial condition results in a substantial difference between the final solutions at time $t>t_0$. For the case of $\lambda<0$ the situation is different; any two solutions converge toward each other and even large perturbation in initial data will finally diminish in the solutions. For the case of $\lambda=0$
 perturbations in data and solution are of the same order. If $\lambda$ is a complex number, $\lambda=a+ib$ say, then
 $$
 y(t)=y_0e^{at}(\cos bt + i\sin bt).
 $$
The behaviour now depends on the sign of $\mathrm{Re}(\lambda)=a$. We have exponential decay for $\mathrm{Re}(\lambda)<0$, exponential growth for $\mathrm{Re}(\lambda)>0$ and oscillatory solutions (parallel curves) for $\mathrm{Re}(\lambda)=0$.
\end{example}
\vsp 

A solution of the ODE $\y'(t) = \f(t, \y)$ is said to be {\bf stable} if for every $\epsilon>0$
there is a $\delta>0$ such that if $\wh \y(t)$ satisfies the ODE and $\|\wh \y(t_0)-\y(t_0)\|\leqslant \delta$, then
\begin{equation}\label{ivp-stab}
\|\wh \y(t)-\y(t)\|\leqslant \epsilon, \quad  \mbox{for all} \quad t \geqslant t_0.
\end{equation}
Hence, for a stable solution, if the initial value is perturbed, the perturbed solution remains close to the original solution.
This roles out the exponential divergence solutions allowed by the ODE.
A stable solution is said to be {\bf asymptotically stable}
if
$$\|\wh \y(t)-\y(t)\| \to 0, \; as \; t\to \infty.$$
This stronger form of stability means
that the original and perturbed solutions not only remain close to each other, they
converge toward each other over time. As we will soon see in detail, the significance
of these concepts for the numerical solution of ODEs is that any errors introduced
during the computation can be either amplified or diminished over time, depending on the stability of the solution.
\vsp 

\begin{workout}
Consider the IVP 
$$
y'(t)=-[y(t)]^2, \quad y(0)=1. \vsp 
$$
Show that the solution is $y(t)=1/(1+t)$. Then solve the perturbed
problem
$$
\wh y'(t)=-[\wh y(t)]^2, \quad \wh y(0)=1+\delta
$$
and show this IVP is stable, and even asymptotically stable. 

Hint: This is a separable differential equation, so the analytical solution can be easily obtained. 
\end{workout}
\vsp 

\begin{example}\label{ex:eigenvalAnal}
For system of ODEs
$$\y'=A\y$$
for $n\times n$ constant diagonalizable matrix $A$, the stability of solutions depends on the sign of the real part of eigenvalues of $A$. See \eqref{sys:solveAy}. In this case eigenvalues with negative real parts yield exponentially decaying solution components,
eigenvalues with positive real parts yield exponentially growing solution components, and eigenvalues with zero real parts give oscillatory solutions.
This means that the solutions of this
ODE are stable if $\mathrm{Re}(\lambda_k)\leqslant 0$ for all $k=1,2,\ldots,n$, and asymptotically stable
if $\mathrm{Re}(\lambda_k)< 0$ for all $k=1,2,\ldots,n$, but unstable if there is any eigenvalue such that
$\mathrm{Re}(\lambda_k)>0$.
\end{example}
\vsp 

For the general case $\y'(t)=A(t)\y(t)$ where $A(t)$ is a time dependent matrix or for the nonlinear equation $\y'(t)=\f(t,\y)$ the stability analysis is more complicated.

\section{Basic numerical methods}\label{sect-numerical-methods}
Although the exact solution of some ODEs can be obtained by few analytic methods, such as method of integration factors,
the solution of most ODEs arising in applications are so complicated that should be computed only by numerical methods.
Even when a solution formula is available, it may involve integrals that can be calculated only by numerical integration formulas.
An analytical
solution of an ODE is a continuous (and sometimes a close form) function in an infinite-dimensional space, while a numerical solution is a table of approximate values of the solution function
at a discrete set of points which can be considered as a vector in a finite dimensional space.

In this section some basic numerical methods are given for solving \eqref{ivp:linear1}. In all methods, starting from the initial condition
$y_0$ at $t=t_0$, the approximate solutions at times $t_1,t_2, \ldots$ are obtained successively by solving an algebraic (system of) difference equations obtained by discretizing the differential equation.

\subsection{Euler's method}
The simplest technique is the {\bf Euler's method}, also called {\em explicit} or {\em forward Euler's method}.
For a given time step $h>0$ assume that $t_k = t_0 + kh$, $k=1,2,\ldots,N$, is a partitioning of time domain $[t_0,b]$ with $b=t_N$.
Let the derivative $y'(t)$ in the ODE
$y'(t) = f(t, y)$ at $t=t_k$ be approximated by the first-order forward difference approximation
$$
y'(t_k) = \frac{y(t_{k+1})-y(t_k)}{h} + \frac{h}{2}y''(\xi_k), \quad t_k\leqslant\xi_k\leqslant t_{k+1}.
$$
By dropping the error term and using the approximate values $y_k$ instead of $y(t_k)$,
we obtain
\begin{shaded}
\vspace*{-0.3cm}
\begin{equation}\label{euler:method}
\begin{split}
 & y_{k+1} = y_{k} + hf(t_k,y_k), \quad k=0,1,\ldots, N-1,\\
 & y_0 = y(t_0)
\end{split}
\end{equation}
\vspace*{-0.3cm}
\end{shaded} 
We use $y_0 = y(t_0)$ or some close approximation of it. Formula \eqref{euler:method} gives a rule for computing $y_1,y_2,\ldots, y_N$ in succession. We bring an example from \cite{Heath:2018}. 
\vsp 

\begin{example}
Consider the simple IVP $y'=y$ with initial value $y_0$ at initial time $t_0=0$. The approximate value $y_1$ is obtained as
$y_1 = y_0 + hy_0 = (1+h)y_0$. It is obvious that $y_1\neq y(t_1)$, thus, the value $y_1$
lies on a different solution curve of the ODE from the one on which we started, as shown in the left side of Figure \ref{fig:euler1}.
The approximate value $y_1$ is the slope of the tangent line for the new (perturbed) solution curve.
Now we continue to obtain the approximate value $y_2$ at $t=t_2$ by starting from $y_1$ with Euler rule $y_2 = y_1+hy_1=(1+h)y_1$.
The approximate value $y_2$ carries both previous approximation error in $y_1$ and a new error introduced in the current step of the Euler's method. In fact, in this step we have moved to still another solution of the ODE, as is again shown in Figure \ref{fig:euler1}.
We can advance to future times $t_3,t_4,\ldots$ until reaching the final time $b=t_N$. In each step a new {\em local truncation error} is introduced and the approximate solution falls down (or rises up) to another solution curve.

Since the solutions of this ODE are unstable, the errors are amplified with time. For an equation with stable solutions, on the other
hand, the errors in the numerical solution do not grow, and for an equation with
asymptotically stable solutions, such as $y'=-y $, the errors diminish with time, as
is shown in the right hand side of Figure \ref{fig:euler1}.

\begin{center}
\includegraphics[scale=0.35]{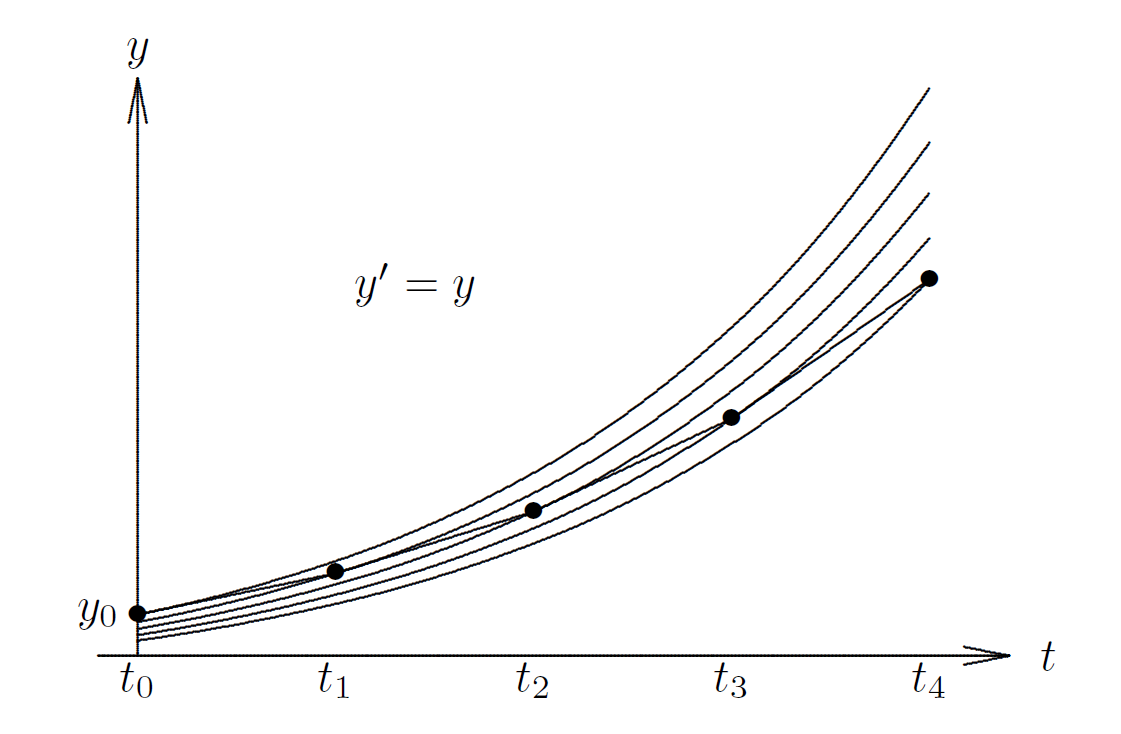}\includegraphics[scale=0.25]{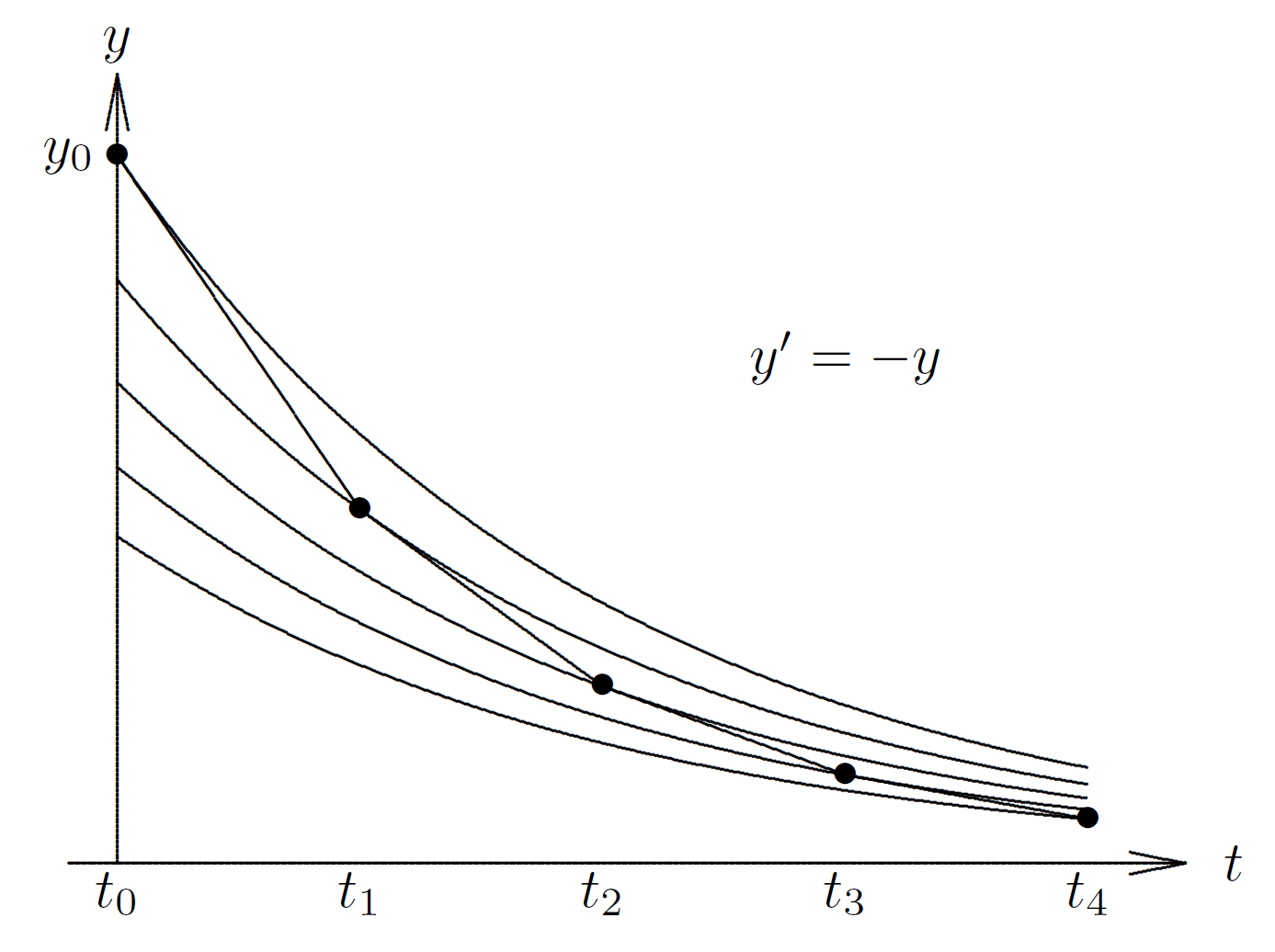}
\captionof{figure}{Steps of Euler's method for $y'=y$ (left) and $y'=-y$ (right).}\label{fig:euler1}
\end{center}
\end{example}
\vsp  
\begin{workout}
Write down Euler's formula for the following IVPs. Compute 1 iteration for the third ODE with $h=0.1$.
\begin{itemize}
\item (a) $y'(t)=te^{-t}-y(t), \; y(0)=1$,
\item (b) $y'(t)=[\cos (y(t))]^2,\; y(0)=0$,
\item (c) $y'(t)=t^3/y(t),\; y(0)=1$.
\end{itemize} 
\end{workout}
If the ODE is a system then $y$ and $f$ in formula \eqref{euler:method} are simply replaced by vectors $\y$ and $\f$. The ODE algorithms are usually simple to program. 
Here is the function of explicit Euler's method. 
\begin{lstlisting}[style = matlab]
function [T,Y] = ExEuler(f,y0,tspan,h)
% Euler's method for IVP system y' = f(t,y) with y(t0) = y0
% Inputs: 
%  f: right hand side function f(t,y)
%  y0: initial condition of size (1 x m)
%  tspan: [t0, tfinal] 
%  h: stepsize
% Output:
%  T: vector of time step
%  Y: solution of size (N x m) 
Y = y0; T = tspan(1);
for t = tspan(1):h:tspan(2)-h
   y = y0 + h*f(t,y0);
   Y = [Y y]; y0 = y; T = [T t+h];
end
\end{lstlisting}

The input \verb+f+ is a function of $t$ and $y$ variables, and can be defined as a separated function in the main script. 
For example we use the following commands for solving the chemical reaction kinetics ODE \eqref{chem:eq1} on interval $[0,3]$ with $k_1=2$ and $k_2 = 1$ and initial conditions 
$y_1(0)=5$ and $y_2(0)=2$. We also set $h = 0.01$. 
\begin{lstlisting}[style = matlab]
%% Setup of the problem
y0 = [5;2]; tspan = [0 3]; h = 0.01;
% ODE call
[t,y] = ExEuler(@func,y0,tspan,h); 
% plot solutions y1 and y2
plot(t,y(1,:),'-b',t,y(2,:),'--r')
set(gca,'TickLabelInterpreter','latex')
xlabel('Time $t$', Interpreter='latex');
ylabel('Solution $y$',Interpreter='latex');
leg = legend('$y_1$','$y_2$'); set(leg,Interpreter='latex'); 
%% definition of right-hand side function f
function yprime = func(t,y)
yprime = [-2*y(1) + y(2); 2*y(1)-y(2)];
end
\end{lstlisting}
 The exact solution was given in Example \ref{ex:chem_solve}
 and Figure \ref{fig:ode2}. If you run this script on your computer you will receive the same plot. 
 \vsp
 \begin{labexercise}
 Solve the dog-jogger problem using the Euler's method and reproduce the plots of Figure \ref{trajdog_fig}. Use $(x(0),y(0))= (60,70)$ and
 \begin{itemize}
 \item $(\xi(t),\eta(t))=(8t,0)$ for $t\in[0,12]$ and $w=10$,
 \item $(\xi(t),\eta(t))=\begin{cases} (8t,0), & t\in[0,7)\\ (8(7-t),0), & t\in[7,12] \end{cases}$ and $w=10$,
 \item $(\xi(t),\eta(t)) = (30+20\cos t, 20 + 15\sin t)$ and $t\in[0,4\pi]$ and $w=10$,
 \item  $(\xi(t),\eta(t)) = (30+20\cos t, 20 + 15\sin t)$ and $t\in[0,4\pi]$ and $w=18$. 
 \end{itemize}
Solve the problem with other jogger's paths and dog speeds.
 \end{labexercise}
 \vsp
\subsubsection{Error analysis of Euler's method}\label{sect:error_anal_euler}

The analysis of Euler's method is useful to understand how it works, to predict the error when using it and perhaps to accelerate its convergence. Moreover, it gives an insight to how analyze other more efficient numerical methods.

\begin{figure}[!th]
\centering
\includegraphics[scale=0.32]{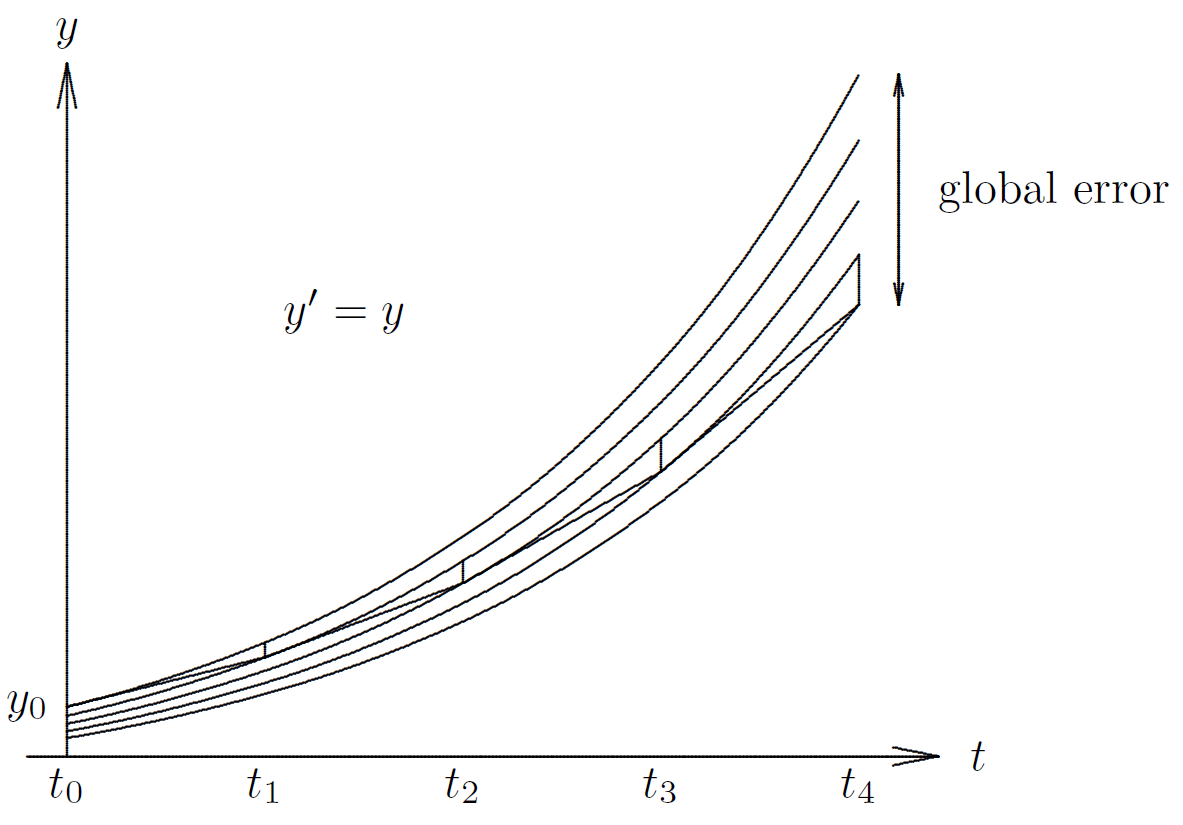}\includegraphics[scale=0.35]{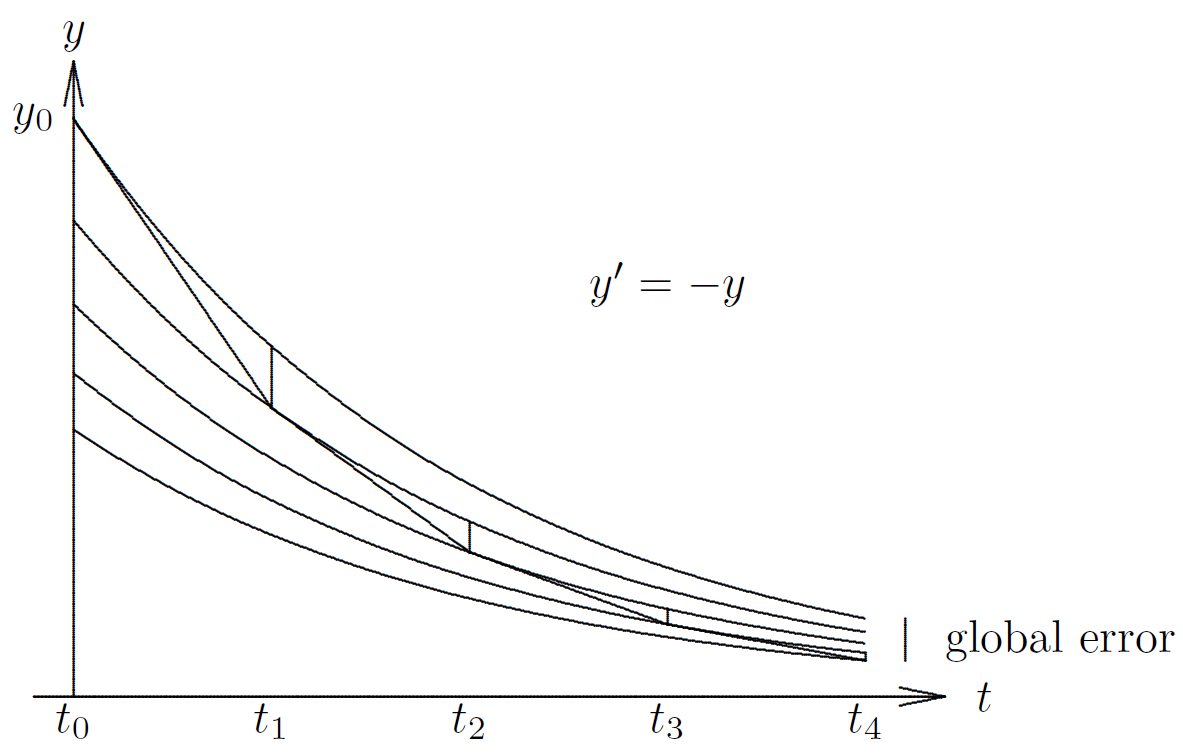}
\caption{Local and global errors of Euler's method for $y'=y$ (left) and $y'=-y$ (right).}\label{fig:euler2}
\end{figure}

We analyze the scaler IVP $y' = f(t,y)$ with $y(t_0)=y_0$ by assuming that
it has a unique solution $y(t)$ on $t_0\leqslant t\leqslant b$ and this solution has a bounded second derivative $y''(t)$ over this interval.
Using the Taylor series formula we have
$$
y(t_{k+1})=y(t_k) + h y'(t_k)+\frac{1}{2}h^2y''(\xi_k)
$$
for some $t_k\leqslant \xi_k\leqslant t_{k+1}$. Using the fact that $y'(t)=f(t,y(t))$ we can write
\begin{equation}\label{euler:tylorser}
y(t_{k+1})=y(t_k) + h f(t_k,y(t_k))+\frac{1}{2}h^2y''(\xi_k).
\end{equation}
The term
$$
\tau_{k+1} = \frac{1}{2}h^2y''(\xi_k)
$$\ \\
is a {\em local truncation error} (or one-step error) for the Euler's method introduced in step $k+1$. Subtracting \eqref{euler:tylorser} from the Euler's rule
$$
y_{k+1} = y_k + hf(t_k,y_k),
$$
we will obtain
\begin{equation}\label{euler:errlocalglobal}
y(t_{k+1})-y_{k+1}=(y(t_k)-y_k) + h [f(t_k,y(t_k))-f(t_k,y_k)]+\tau_{k+1},
\end{equation}
which shows that the error in $y_{k+1}$ consists of two parts:
\begin{itemize}
  \item (1) the newly introduced local truncation error $\tau_{k+1}$,
  \item (2) the {\em propagated error} $(y(t_k)-y_k) + h [f(t_k,y(t_k))-f(t_k,y_k)]$.
\end{itemize}
If we assume that
$$
e_k = y(t_k)-y_{k} \vsp 
$$
then $e_N$ would be the {\em global error} at the final time $t=t_N$. The global
error reflects not only the
local error at the final step, but also the compounded effects of the local errors at all
previous steps. Unless in some specific situations where $f$ is independent of $y$, {\bf the global error is not simply the sum of the local errors.} For example, for ODE $y'=y$ since the solutions are diverging, the local errors at each step are
magnified over time, so that the global error is greater than the sum of the local
errors, as shown in Figure \ref{fig:euler2}, where the local errors are indicated by small vertical
bars between solutions and the global error is indicated by a bar at the end. On the other hand for ODE $y'=-y$ since the solutions are converging, the global error is less than the sum of the local errors.
It is obvious that for the only case $f(t,y)=g(t)$, where the solutions are parallel curves, the global error is the direct sum of local errors. Because in this case we have $[f(t_k,y(t_k))-f(t_k,y_k)]=g(t_k)-g(t_k)=0$ and \eqref{euler:errlocalglobal}
reduced to $e_{k+1}=e_k + \tau_{k+1}$ for $k=0,1,\ldots,N-1$. We then simply have $e_N=\tau_1+\tau_2+\cdots+\tau_N$. However, for a general $f(t,y)$ we have to analyze the effect of $[f(t_k,y(t_k))-f(t_k,y_k)]$ in each step. Considering $f(t,y)$ as a function of $y$ and using the {\em mean value theorem}, we can write
$$
f(t_k,y(t_k))-f(t_k,y_k) = \frac{\partial f}{\partial y}(t_k,\eta_k)(y(t_k)-y_k)
$$
for some $\eta_k$ between $y(t_{k})$ and $y_k$. Then, \eqref{euler:errlocalglobal} yields
\begin{equation}\label{euler:errlocalglobal2}
e_{k+1}=\left(1 + h \frac{\partial f}{\partial y}(t_k,\eta_k) \right)e_k +\tau_{k+1},
\end{equation}
which shows that the amplification (or diminishing) factor for propagation error is $\left(1 + h \frac{\partial f}{\partial y}(t_k,\eta_k) \right)$ which is related to the stability of solutions via the Lipschitz constant of $f$.

We assume that the function $f(t,y)$ satisfies the following stronger Lipschitz condition: there exists a constant $L>0$ such that
\begin{equation*}
  |f(t,y)-f(t,\tilde y)|\leqslant L|y-\tilde y|,\quad \forall (t,y),\,(t,\tilde y)\in [t_0,b]\times \R.
\end{equation*}
In \eqref{euler:errlocalglobal2} if we go through absolute value of $e_k$, the term $\frac{\partial f}{\partial y}(t_k,\eta_k)$ can be replaced by the Lipschitz constant $L$ to obtain
\begin{equation}\label{euler:errlocalglobal3}
|e_{k+1}|\leqslant(1 + h L)|e_k| +|\tau_{k+1}|, \quad k=0,1,\ldots,N-1
\end{equation}
from \eqref{euler:errlocalglobal2}. We assume that
$$
\tau(h) = \frac{1}{2}h\|y''\|_\infty = \frac{1}{2}h\max_{t_0\leqslant t\leqslant b}|y''(t)|.
$$
Then we have $|\tau_k|\leqslant h\tau(h)$ for $k=1,\ldots,N$.
Now, apply \eqref{euler:errlocalglobal3} recursively, we obtain
$$
|e_k|\leqslant (1+hL)^k|e_0| + [1+(1+hL)+(1+hL)^2+\cdots+(1+hL)^{k-1}] h\tau(h).
$$
Using the formula for sum of geometric series,
$$
1+r+r^2+\cdots+r^{n-1} = \frac{r^n-1}{r-1}, \quad r\neq 1,
$$
we obtain
$$
|e_k|\leqslant (1+hL)^k|e_0| +\left[ \frac{(1+hL)^k-1}{L} \right] \tau(h).
$$
Now we use the standard formula
$$
1+x \leqslant e^x, \quad x\geqslant 0
$$
with $x=hL$ to obtain
$$
|e_k|\leqslant e^{khL}|e_0| +\left[ \frac{e^{khL}-1}{L} \right] \tau(h).
$$
On the other hand $kh = t_k-t_0$, so we can write
$$
|e_N|\leqslant e^{L(t_N-t_0)}|e_0| +\left[ \frac{e^{L(t_N-t_0)}-1}{L} \right] \tau(h).
$$
If $y_0=y(t_0)$ or $|e_0|=|y(t_0)-y_0|\leqslant c_1 h$ then we have a global error of order $h$ for the Euler's method, i.e.,
\begin{equation}\label{euler:errorbound}
|e_N| \leqslant Ch, \quad C = c_1e^{L(b-t_0)} + \frac{1}{2}\left[ \frac{e^{L(b-t_0)}-1}{L} \right]\|y''\|_\infty.
\end{equation}\ \\
Therefore, the Euler's method is said to converge with order $1$. This order of convergence is obtained by assuming $y$ to have a continuous second derivative $y''$ over interval $[t_0,b]$. When such assumption fails the error bound  \eqref{euler:errorbound} no longer holds. See Workout \ref{wo_eulersmoothness}.

Since the error bound \eqref{euler:errorbound} uses the Lipschitz constant $L$ instead of $\frac{\partial f}{\partial y}$ and ignores the sign of $\frac{\partial f}{\partial y}$, it sometimes produces a very pessimistic numerical bound for the error. If
\begin{equation}\label{euler:negativefpar}
\frac{\partial f}{\partial y}(t,y)\leqslant 0
\end{equation}
then we may have a smaller than $1$ amplification factor $\left(1 + h \frac{\partial f}{\partial y}(t_k,\eta_k) \right)$ instead of
$(1+hL)$ which is always bigger than $1$. In this case (negative partial derivative of $f$) if we assume that
$$
L = \sup_{t\in[t_0,b],y\in\R}\left| \frac{\partial f}{\partial y}(t,y) \right|
$$
and $h$ is chosen so small that $1-hL\geqslant -1$ then we have
$$
1\geqslant 1+h\frac{\partial f}{\partial y}(t_k,\eta_k)\geqslant 1-hL\geqslant -1,
$$
and
from \eqref{euler:errlocalglobal2} we can write
$$
|e_{k+1}|\leqslant |e_k|+|\tau_k|,\quad k=0,1,\ldots,N-1.
$$
Applying this bound recursively, we obtain
\begin{equation}\label{euler:bound2}
|e_N|\leqslant |e_0| + (b-t_0)\tau(h) = Ch, \quad C = c_1 + (b-t_0)\|y''\|_\infty
\end{equation}
where $|e_0|\leqslant c_1h$ is assumed. The constant $C$ behind $h$ in bound \eqref{euler:bound2} is much smaller than that in bound
\eqref{euler:errorbound} which contains the exponential terms. But, the error bound \eqref{euler:bound2} is valid with restrictive assumption \eqref{euler:negativefpar}.
\vsp 

\begin{example}
For simple IVP
$$
y'(t)=-y(t), \quad y(0)=1, \quad 0\leqslant t\leqslant b,
$$
we have $\partial f(t,y)/\partial y=-1$ and $L=1$. The true solution is $y(t)=e^{-t}$, hence $\|y''\|_\infty=1$. With $y_0=y(0)=1$ ($|e_0|=0$). From the bound \eqref{euler:errorbound} we have
$$
|e_N|\leqslant \frac{1}{2}h(e^b-1)
$$
which shows the convergence with order $h$. However the constant in the bound grows exponentially in $b$. For example with $b=5$ the bound becomes approximately $73.7 h$ which is far larger than the actual error in Table \ref{tb:euler1}. The error bound \eqref{euler:bound2}, on the other hand gives
$$
|e_N|\leqslant \frac{1}{2}bh
$$
which is very close to the actual error with different $h$ in Table \ref{tb:euler1}. The fifth column of the table also confirms the theoretical order $1$ for the method.

\begin{center}
\captionof{table}{Euler's method: numerical solutions, errors, orders, and error bounds for IVP $y'=-y$ with $y_0=1$ at time $t=b=5$.}\label{tb:euler1}
\begin{tabular}{l|cccccc}
  \hline
  $h$ & $y_N$ & $|e_N|$ & $|e_N|/|y(b)|$ & order & $\frac{1}{2}h(e^b-1)$ & $\frac{1}{2}bh$ \\
  \hline
  $0.2$     & $3.778\ee-3$ & $2.960\ee-3$ & $4.393\ee-1$& $-$    & $14.74$ & $0.5 $ \\
  $0.1$     & $5.154\ee-3$ & $1.584\ee-3$ & $2.351\ee-1$& $0.90$ & $7.37 $ & $0.25$ \\
  $0.05$    & $5.921\ee-3$ & $8.174\ee-4$ & $1.213\ee-1$& $0.95$ & $3.69 $ & $0.12$ \\
  $0.025$   & $6.323\ee-3$ & $4.149\ee-4$ & $6.158\ee-2$& $0.98$ & $1.84 $ & $0.06$ \\
  $0.0125$  & $6.529\ee-3$ & $2.090\ee-4$ & $3.102\ee-2$& $0.99$ & $0.92 $ & $0.03$ \\
  $0.00625$ & $6.633\ee-3$ & $1.049\ee-4$ & $1.557\ee-2$& $0.99$ & $0.46 $ & $0.02$ \\
  \hline
\end{tabular}
\end{center}
The results of this table are obtained by executing the following code:
\begin{lstlisting}[style = matlab1]
ExactSol = @(t) exp(-t);
b = 5; 
h = 0.2;
for n = 1:6
    [t,y] = ExEuler(@(t,y) -y, 1, [0 b], h); 
    AppSol(n) = y(end);
    ABSerr(n) = abs(AppSol(n)-ExactSol(b));
    RELerr(n) = ABSerr(n)/ExactSol(b);
    h = h/2;
end
Order = log2(ABSerr(1:5)./ABSerr(2:6));
fprintf('y_N = \n'); fprintf('  %1.3e\n',AppSol); 
fprintf('abs_err = \n'); fprintf('  %1.3e\n',ABSerr);
fprintf('rel_err = \n'); fprintf('  %1.3e\n',RELerr);
fprintf('order = \n'); fprintf('  %1.2f\n',Order);
\end{lstlisting}
Try understanding why the numerical orders are computed using that logarithmic formula in line 11 of the script above!

\begin{center}
\captionof{table}{Euler's method: numerical solutions, errors, orders, and error bounds for IVP $y'=+y$ with $y_0=1$ at time $t=b=5$}\label{tb:euler2}
\begin{tabular}{l|cccccc}
  \hline
  $h$ & $y_N$ & $|e_N|$ & $|e_N|/|y(b)|$ & order & $\frac{1}{2}h(e^b-1)\|y''\|_\infty$ & \\
  \hline
  $0.2$     & $9.540\ee+1$ & $5.302\ee+1$ & $3.572\ee-1$& $-$    & $2.187\ee+3$ & \\
  $0.1$     & $1.174\ee+2$ & $3.102\ee+1$ & $2.090\ee-1$& $0.77$ & $1.093\ee+3 $ & \\
  $0.05$    & $1.315\ee+2$ & $1.691\ee+1$ & $1.140\ee-1$& $0.88$ & $5.470\ee+2$ & \\
  $0.025$   & $1.396\ee+2$ & $8.849\ee+0$ & $5.963\ee-2$& $0.93$ & $2.735\ee+2 $ & \\
  $0.0125$  & $1.439\ee+2$ & $4.529\ee+0$ & $3.052\ee-2$& $0.97$ & $1.367\ee+2$ & \\
  $0.00625$ & $1.461\ee+2$ & $2.291\ee+0$ & $1.544\ee-2$& $0.98$ & $6.837\ee+1$ & \\
  \hline
\end{tabular}
\end{center}

Now consider the IVP
$$
y'(t)=y(t), \quad y(0)=1, \quad 0\leqslant t\leqslant b. 
\vsp 
$$
For this problem the
the error bound \eqref{euler:bound2} is not applicable because $\partial f(t,y)/\partial y=1>0$. However, the error bound \eqref{euler:errorbound}
is nearly sharp for this IVP. The exact solution is $y(t)=e^t$ and $\|y''\|_\infty=e^b$. See the results in Table \ref{tb:euler2}.

\end{example}
\vsp

\begin{labexercise}\label{py_ex_euler}
Solve the following problems using Euler's method with stepsizes $h=0.2,0.1,0.05,0.025,0.0125,0.00625$. Compute the relative errors using the true solutions $y(t)$. In each case plot the error function in terms of $h$ in the log-log scale, and compute the computational order of convergence.
\begin{itemize}
\item (a) $y'(t)=te^{-t}-y(t), \; 0\leqslant t\leqslant 10,\; y(0)=1$, with exact solution $y(t)=(1+0.5t^2)e^{-t}$.
\item (b) $y'(t)=[\cos (y(t))]^2, \; 0\leqslant t\leqslant 10,\; y(0)=0$, with exact solution $y(t)=\tan^{-1}(t)$.
\item (c) $y'(t)=t^3/y(t), \; 0\leqslant t\leqslant 10,\; y(0)=1$, with exact solution $y(t)=\sqrt{0.5t^4+1}$.
\end{itemize} 
\end{labexercise}
\vsp

\subsection{General explicit one-step methods}\label{sect:general-explicit}
The Euler's method is a one-step method (counterpoise to multistep methods) as in each time level the approximate solution is obtained from merely the previous time level.
A general explicit one-step method has the form
\begin{equation}\label{onestep:general}
  y_{k+1} = y_k + h\psi(t_k,y_k,h)
\end{equation}
for a more general function $\psi$ instead of $f$. We will assume that $\psi(t,y,h)$ is
continuous in $t$ and $h$ and Lipschitz continuous in $y$, with Lipschitz constant $\tilde L$ that is
generally related to the Lipschitz constant $L$ of $f$.
\vsp 

\begin{example}
The choice $\psi(t,y,h)=f(t,y+\frac{h}{2}f(y))$ results in the two-stage {\em Runge-Kutta method} that will be addressed later. For this scheme
we can simply show that $\psi$ has Lipschitz constant $\tilde{L}=L+\frac{h}{2}L^2$ where $L$ is the Lipschitz constant of $f$.
\end{example}
\vsp 

A one-step method is said to be {\em consistent} if
\begin{equation}\label{one-step:consistency}
\psi(t,y,0) = f(t,y),
\end{equation}
for all $t,y$, and $\psi$ is continuous in $h$. The consistency, indeed, implies that the local truncation error of method \eqref{onestep:general} is at least of order $h^2$, because
\begin{align*}
  \tau_{k+1}&=y(t_{k+1})-y(t_k)-h\psi(t_k,y_k,h)\\
            &=hy'(t_k)-h\psi(t_k,y(t_k),h) + \mathcal{O}(h^2)\\
            &=h\psi(t_k,y(t_k),0)-h[\psi(t_k,y(t_k),0)+\mathcal O(h)] + \mathcal{O}(h^2)\\
            &=\mathcal{O}(h^2)
\end{align*}
The error analysis of general one-step methods can be obtained in a similar way as was done for the Euler's method. First, like as \eqref{euler:tylorser}, the truncation error is obtained as
\begin{equation*}
  \tau_{k+1} = y(t_{k+1})-y(t_k)-h\psi(t_k,y(t_k),h),
\end{equation*}
and then
\eqref{euler:errlocalglobal} is modified to
\begin{equation*}
e_{k+1}=e_k + h [\psi(t_k,y(t_k),h)-\psi(t_k,y_k,h)]+\tau_{k+1}.
\end{equation*}
The remaining parts of analysis follow a same direction only $L$ should be replaced by $\tilde L$ in new error bounds.
We then can conclude the following theorem.
\begin{theorem}\label{thm:converg-onesteps}
If $\psi(t,y,h)$ is continuous in all its arguments and is Lipschitz continuous in its second argument, and the consistency condition \eqref{one-step:consistency} holds, then
the explicit one-step method \eqref{onestep:general} is convergent with a global error of at least order $h^1$. If the local truncation errors
$\tau_{k+1}$ behave as $h^{p+1}$, then the global order is of order $h^p$.
\end{theorem}
\vsp 

\subsection{Zero-stability}

In the convergence proof of the one-step methods we observed the effect of amplification factor in propagating the local errors which finally was summed up to factor $C$ in the error bound \eqref{euler:errorbound}. Although this factor grows in $b$, it is bounded independent
of $h$ as $h\to 0$. Consequently the method is stable. This form of stability for a numerical method is often called {\bf zero-stability},
since it is concerned with the stability of the method in the limit as $h$ tends to zero.
To see this observation in a more relevant presentation to stability of the original IVP \eqref{ivp:form}, assume that the initial condition $y_0$ is perturbed by $\ep$ and define numerical solutions of the perturbed problem by
$$
z_{k+1}=z_k + hf(t_k,z_k), \quad z_0 = y_0+\ep.
$$
For comparing two numerical solutions $y_k$ and $z_k$, let $e_k = z_k-y_k$. Then $e_0=\ep$ and subtracting from $y_{k+1}=y_k+hf(t_k,y_k)$ we obtain
$$
e_{k+1}=e_k+h[f(t_k,z_n)-f(t_k,y_k)].
$$
This is exactly the same form as \eqref{euler:errlocalglobal} with $\tau_{k+1}$ set to be zero. Using the same procedure we obtain
$$
|e_{k+1}|\leqslant (1+hL)^k|e_0|\leqslant e^{L(t_k-t_0)}|\ep|.
$$
Consequently, we can write
\begin{equation*}
  \max_{0\leqslant k\leqslant N}|z_k-y_k| \leqslant e^{L(b-t_0)}|\ep| \vsp
\end{equation*}
which is the analog to the stability result \eqref{ivp-stab} for the original IVP \eqref{ivp:form}.
The {\em zero-stability} is different from
other forms of stability which are of equal importance in practice. The
fact that a method is zero-stable (and converges as $h\to 0$) is no guarantee that it will give
reasonable results on the particular grid with $h>0$ that we want to use in practice. Other
stability issues of a different nature will be taken up in the next sections.
\vsp 

\subsection{Absolute stability}

The zero-stability is needed to guarantee
convergence of a numerical method as $h\to 0$. In practice, however, we need to perform a single calculation using a given positive stepsize $h>0$. Moreover, to minimize the computational cost a larger as possible $h$ (consistent with our desired accuracy) is usually preferred.
A stronger form of stability than the zero-stability is required in this case to force the method to work for this particular stepsize
$h$.
Let's illustrate the situation in three numerical examples borrowed from \cite{LeVeque:2007}.
\vsp 

\begin{example}\label{ex:absstab}
We apply the Euler's method on a simple IVP of the form
$$
y'(t) = -\sin(t), \quad 0\leqslant t\leqslant 2, \quad y(0)=1
$$
with exact solution $y(t)=\cos t$. Since $f(t,y)=\sin(t)$ is independent of $y$, ($L=0$) the global error is the sum of local errors
$$
|\tau_k| = \frac{h^2}{2}|y''(\xi_k)|\leqslant \frac{h^2}{2}.
$$
Indeed we have
$$
|e_N|\leqslant (b-t_0)\tau(h) = 2\tau(h) = h.
$$
Suppose we want to compute the solution at $t=2$ with a global error less than $0.001$. According to the error bound it suffices to
take $h = 0.001$ and obtain the approximate solution after $2000$ time steps. The computed solution $y_{2000}\doteq -0.4156921$
has error $|e_{2000}|=|y_{2000}-\cos(2)|\doteq 0.45\times 10^{-3}$.

Now, we change the IVP to
\begin{equation}\label{exeulerlam}
y'(t) =\lambda(y-\cos t)-\sin(t), \quad 0\leqslant t\leqslant 2, \quad y(0)=1,
\end{equation}
for some constant $\lambda$. The exact solution is $y(t)=\cos t$, as before. Let $\lambda=-10$. The error bound
\eqref{euler:bound2} suggests again the global error $|e_N|\leqslant h$. For this reason, we again choose $h=0.001$ for a global error less than $0.001$.
The computed solution now is $y_{2000}\doteq -0.4161629$ with error
$|e_{2000}|\doteq 0.16\times 10^{-4}$ which is even better than the previous one.

Let us examine some larger (in magnitude) $\lambda$. Let $\lambda = -2100$. Executing the Euler's method gives $y_{2000}\doteq 0.15\times 10^{77}$ which is far away from the exact solution and shows a blown up in computations. The method is zero-stable and we proved that when
$h\to 0$ it is convergent. Indeed, for sufficiently small stepsizes we achieve accurate results as reported in Table \ref{tb:euler3}.

\begin{center}
   \captionof{table}{Global errors for the Euler's method with different stepsizes.}\label{tb:euler3}
  \begin{tabular}{l|ccccc}
    \hline
    $h$   & $0.001$ & $0.00097$ & $0.00095$ & $0.0008$ & $0.0004$ \\
    \hline
    $|e_N|$ & $0.15\ee+77$ & $0.77\ee+26$ & $0.40\ee-07$ & $0.79\ee-07$ & $0.40\ee-08$ \\
    \hline
  \end{tabular}
\end{center}

Something dramatic happens for values of $h$ between $0.00095$ and $0.00097$. For smaller values of $h$ we get very good results, whereas for larger values of $h$ the solution blows up. To find the reason, we come back to
\eqref{euler:errlocalglobal2} where we have for the linear IVP with $f(t,y)=\lambda(y-\cos t)-\sin(t)$ the recursion
$$
e_{k+1} = (1+\lambda h) e_{k} + \tau_{k+1}.
$$
This means that in each time step the previous error is multiplied by factor $(1+\lambda h)$. With $\lambda=-2100$
and $h = 0.001$ we have $|1+\lambda h| = 1.1$. After $2000$ steps the truncation error introduced in
the first step has grown by a factor of roughly $(1.1)^{2000}\approx 10^{82}$, which is consistent
with the error actually seen. Note that with $\lambda=10$, we have $|1+\lambda h|=0.99$ causing a decay in
the effect of previous errors in each step. For the first case, i.e., $\lambda=0$, the amplification factor is $1$ the reason why we got a worse result in this case than the case of $\lambda=-10$.

Consequently, we can argue that for values of $h$ satisfying 
$$
|1+\lambda h|\leqslant 1
$$
the Euler method produces stable and accurate results for IVP \eqref{exeulerlam}. In the case of $\lambda=-2100$ the above criterion suggests the values of $h$ smaller than $2/2100\doteq 0.000952$.

\end{example}
\vsp 

\begin{remark}
Note that the exponential growth of errors for some positive values of $h$ in the previous example does not contradict zero-stability or convergence
of the method in any way. The method does converge as $h\to 0$.
\end{remark}
\vsp 

\begin{figure}[!th]
\begin{center}
\begin{tikzpicture}[scale=2]
    \node[color=black] at (-.5,1.2) {{\tiny Explicit Euler}};
    \draw[rectangle] (-3/2,-1) rectangle (0.5cm,1cm);
        \draw[fill=blue!10] (-0.5,0) circle (0.5cm);
    \draw[-] (-3/2,0) -- (1/2,0)  node[right] {} ;
    \draw[-] (0,-1) -- (0,1)  node[right] {} ;
    \node[color=black] at (-3/2,-1.1) {{\tiny $-3$}};
    \node[color=black] at (-1,-1.1) {{\tiny $-2$}};
    \node[color=black] at (-1/2,-1.1) {{\tiny $-1$}};
    \node[color=black] at (0,-1.1) {{\tiny $0$}};
    \node[color=black] at (1/2,-1.1) {{\tiny $1$}};
    \node[color=black] at (-1.7,-1) {{\tiny $-2$}};
    \node[color=black] at (-1.7,-1/2) {{\tiny $-1$}};
    \node[color=black] at (-1.6,0) {{\tiny $0$}};
    \node[color=black] at (-1.6,1/2) {{\tiny $1$}};
    \node[color=black] at (-1.6,1) {{\tiny $2$}};
\end{tikzpicture}
\begin{tikzpicture}[scale=2]
    \node[color=black] at (-.5,1.2) {{\tiny Implicit Euler}};
    \draw[fill=blue!10] (-3/2,-1) rectangle (0.5cm,1cm);
    \draw[fill=blue!0] (-0.5,0) circle (0.5cm);
    \draw[-] (-3/2,0) -- (1/2,0)  node[right] {} ;
    \draw[-] (-1,-1) -- (-1,1)  node[right] {} ;
    \node[color=black] at (-3/2,-1.1) {{\tiny $-1$}};
    \node[color=black] at (-1,-1.1) {{\tiny $0$}};
    \node[color=black] at (-1/2,-1.1) {{\tiny $1$}};
    \node[color=black] at (0,-1.1) {{\tiny $2$}};
    \node[color=black] at (1/2,-1.1) {{\tiny $3$}};
    \node[color=black] at (-1.7,-1) {{\tiny $-2$}};
    \node[color=black] at (-1.7,-1/2) {{\tiny $-1$}};
    \node[color=black] at (-1.6,0) {{\tiny $0$}};
    \node[color=black] at (-1.6,1/2) {{\tiny $1$}};
    \node[color=black] at (-1.6,1) {{\tiny $2$}};
\end{tikzpicture}
\end{center}
\caption{Absolute stability regions of explicit and implicit Euler's methods}\label{fig:eulerregabs}
\end{figure}
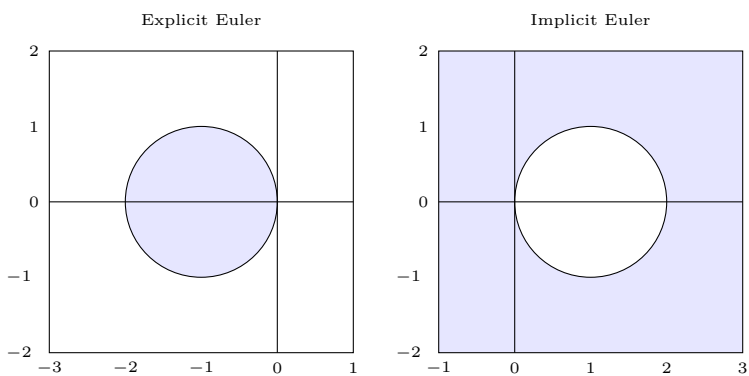

Example \ref{ex:absstab} shows that another notion of stability is needed to force a numerical method to produce stable and accurate  results with a given step length $h>0$. There exists a wide variety of ``stability'' notions but one that is most basic is the {\bf absolute stability}.
This kind of stability is based on the linear {\em test equation}
\begin{equation}\label{lineareq_test}
  y'(t)=\lambda y(t), \quad \lambda\in \C.
\end{equation}
The restriction on the step length $h>0$ on which the method will work for test ODE \eqref{lineareq_test} is called the
absolute stability conditions. For example, the Euler's method if is applied on \eqref{lineareq_test} yields
$$
y_{k+1}=(1+\lambda h)y_k = (1+\lambda h)^2y_{k-1} = \cdots = (1+\lambda h)^{k+1}y_0.
\vsp 
$$
In order to prevent a blown up in the solution when $k\to\infty$ we should impose the condition
\begin{equation*}
  |1+\lambda h|\leqslant 1.
\end{equation*}
If $\lambda \in \R$, the stability condition will be $-2\leqslant z\leqslant 0$ for $z\equiv \lambda h$. This implies that
$$
\lambda \leqslant 0 \;\; \mbox{and}\;\; 0\leqslant h\leqslant \frac{2}{-\lambda}.
$$
In the general case $\lambda\in\C$, the complex number $z$ should satisfy $|1+z|\leqslant 1$ which means that $z$ should lie inside and on a circle with center $(-1,0)$ and radius $1$ in the complex plane. This region is called the {\bf region of absolute stability} of the Euler's method. See shadow part on the left hand side of Figure \ref{fig:eulerregabs}.

\begin{definition}
By applying a one-step method on the test problem \eqref{lineareq_test} we get
$$
y_{k+1}= R(z)y_k, \quad z = \lambda h
$$
for some function $R(z)$,
and then the absolute stability region of the method is defined to be
$$
S = \{z\in\C: |R(z)|\leqslant 1\}.
$$
\end{definition}
The stability region of the explicit Euler's method is rather small compared to other (usually) implicit methods. This will impose a serious restriction on stepsize $h$ to guarantee the convergence. This restriction together with slow convergence rate of Euler's method convince us to study and develop other more efficient ODE solvers. For comparison, the absolute stability region of the implicit Euler's method (see the next section) is drown on the right hand side of Figure \ref{fig:eulerregabs}. It contains all the complex plane except a unit circle with center $(1,0)$.

Although the absolute stability region is determined by testing the method on simple linear ODE \eqref{lineareq_test}, it yields information that is typically useful in determining an appropriate
step length in nonlinear problems as well.

For a system of ODEs of the form
\begin{equation}\label{y1Ay}
\y'(t) = A\y(t), \quad A\in\R^{n\times n}
\end{equation}
where $A$ is diagonalizable with eigenvalues $\lambda_\ell$, $\ell=1,2,\ldots,n$,
a numerical method is absolutely stable if $z_\ell = \lambda_\ell h $ all lie in the absolute stability region of the method in the scaler case. The proof is simple, because as we observed in \eqref{sys:solveAy} the system can be decoupled to $n$ scaler ODE
$$
u'_\ell = \lambda_\ell u_\ell, \quad \ell=1,2,\ldots,n.
$$
Now, we investigate a numerical solution of a simple partial differential equation (PDE) using the method of lines (MOL) which results in
a linear system of ODEs.
\ \\
\begin{example}\label{ex:heat}
Consider the linear {\em diffusion} equation
$$
\frac{\partial u(x,t)}{\partial t} = \frac{\partial^2 u(x,t) }{\partial x^2},\quad 0\leqslant x\leqslant 1, \quad t\geqslant 0
$$
with homogeneous Dirichlet boundary conditions $u(0,t)=u(1,t)=0$ and initial condition $u(x,0)=u^0(x)$.
The method of lines (MOL) solution if is applied on this problem
with the central difference approximation
$$
\frac{\partial^2 u}{\partial x^2}(x_k,t)\approx \frac{u(x_{k+1},t)-2u(x_k,t)+u(x_{k-1},t)}{(\Delta x)^2}, \quad \Delta x = \frac{1}{m+1},
$$
with $x_k=k\Delta x$,
leads to a system of equations of the form \eqref{y1Ay} with
$$
\y(t) = \begin{bmatrix}
      u(x_1,t) \\
      u(x_2,t) \\
      \vdots \\
      u(x_{m-1},t) \\
      u(x_m,t)
    \end{bmatrix}, \quad
A = \frac{1}{(\Delta x)^2}\begin{bmatrix}
                        -2 & 1 &  &  &  \\
                         1 & -2& 1 &  &  \\
                         &\ddots &  \ddots & \ddots &  \\
                         & & 1 & -2& 1 \\
                         & &   & 1 & -2
                         \end{bmatrix}.
$$
The matrix $A$ is symmetric and tridiagonal. There exists a close formula for its eigenvalues:
\begin{align*}
    \lambda_\ell = \frac{2}{(\Delta x)^2}(\cos(\pi\ell \Delta x)-1)= \frac{-4}{(\Delta x)^2}\sin^2(\frac{\pi}{2}\ell \Delta x), \quad \ell=1,\ldots,m.
\end{align*}
The distribution of eigenvalues for two different matrix size $m=10,20$ (or $\Delta x=1/11,1/21$) are displayed in Figure \ref{fig:eulereig}.
\begin{center}
  \includegraphics[width=15cm]{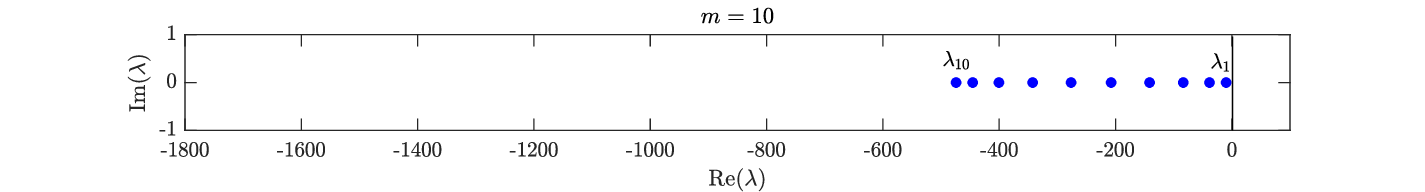}\\
  \includegraphics[width=15cm]{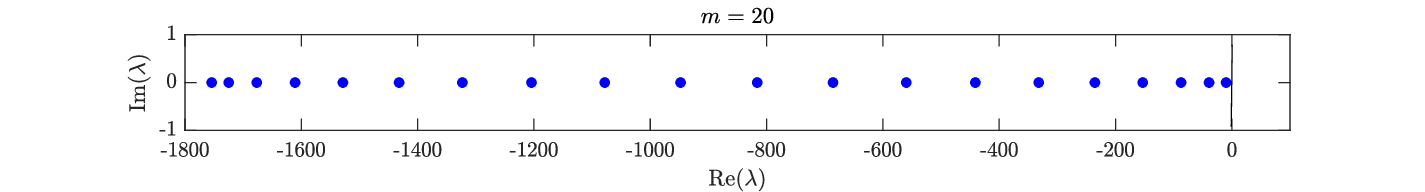}
  \captionof{figure}{Distribution of eigenvalues of matrix $A$.}\label{fig:eulereig}
\end{center}

All eigenvalues are real (because $A$ is symmetric) and fall on the left-half (complex) plane. If one insists to apply the explicit Euler's method for solving this system then the step length  $h$ should be chosen small enough such that
all $\lambda_\ell h$ lie in the absolute stability region of the method. Since the largest (in magnitude) eigenvalue is $\lambda_m$, the absolute stability is guaranteed if the step length  is chosen equal to or less than
$$
\frac{2}{-\lambda_m} = \frac{(\Delta x)^2}{2\sin^2(\frac{\pi}{2}m\Delta x)}
$$
Since $\sin^2(\frac{\pi}{2}m\Delta x)<1$, it is enough to take the step length  equal to or less than $\frac{1}{2}(\Delta x)^2$. This is a serious restriction for numerical solution of such PDE.
\end{example}
\vsp
\begin{remark}
A method is zero-stable if the origin belongs to its region of absolute stability.
\end{remark}

\subsection{Implicit methods}
Euler's method is an explicit method in that it uses the already known information at time $t_k$ to
advance the solution to time $t_{k+1}$. However, this method has a rather small stability region.
The {\bf implicit Euler's method} (backward Euler's method)
is obtained by approximating $y'(t_k)$ by the first-order backward difference approximation
$$
y'(t_k) = \frac{y(t_{k})-y(t_{k-1})}{h} - \frac{h}{2}y''(\xi_k), \quad t_{k-1}\leqslant\xi_k\leqslant t_{k}.
$$
By dropping the error term and using the approximate values $y_k$ instead of $y(t_k)$,
we obtain an algebraic equation $y_k = y_{k-1}+hf(t_k,y_k)$ for $k=1,2,\ldots$. Shifting the index by $1$, the implicit Euler's method
is obtained as
\begin{shaded}
\vspace*{-0.3cm}
\begin{equation}\label{imeuler:method}
\begin{split}
  &y_{k+1} = y_{k} + hf(t_{k+1},y_{k+1}), \quad k=0,1,\ldots, N-1, \\
  & y_0=y(t_0).
  \end{split}
\end{equation}
\vspace*{-0.3cm}
\end{shaded} 
This scheme is implicit because we must evaluate $f$ with the
argument $y_{k+1}$ before we know its value. If $f$ is a
nonlinear function in $y$ then a rootfinding method such
as fixed-point iteration or Newton's method can be used. A good starting guess
for the iteration is the solution at the previous time step or one step solution of the explicit Euler's method. If $f$ is Lipschitz continuous and $h$ is small enough it can be proved that
$y-y_k+hf(t_{k+1},y)=0$ has a unique root.

Usually, a simple iteration method is efficient for solving the nonlinear equation  in each step. In step $k+1$, given an initial guess
$y_{k+1}^{(0)}$, we define $y_{k+1}^{(1)},y_{k+1}^{(2)},\ldots$ by
\begin{equation}\label{imeuler:iter}
  y_{k+1}^{(j+1)}=y_k + hf(t_{k+1},y_{k+1}^{(j)}),\quad j=0,1,2,\ldots.
\end{equation}
Subtracting \eqref{imeuler:iter} from \eqref{imeuler:method}, we obtain
$$
 y_{k+1}-y_{k+1}^{(j+1)}=h\left[f(t_{k+1},y_{k+1})-f(t_{k+1},y_{k+1}^{(j)})\right].
$$
If we assume that $f$ is Lipschitz continuous with constant $L$ then we can write
$$
 |y_{k+1}-y_{k+1}^{(j+1)}|\leqslant hL |y_{k+1}-y_{k+1}^{(j)}|.
$$
This means that if $h$ is chosen small enough such that
\begin{equation}\label{imeuler:Lh}
  hL\leqslant 1
\end{equation}
then the error will converge to zero for a sufficiently good initial guess $y_{k+1}^{(0)}$. In practice, usually one step of the explicit Euler's method, i.e.,
$$
y_{k+1}^{(0)} = y_k+hf(t_k,y_k)
$$
is used as an initial guess in each step of the implicit Euler's method. This is called a {\em predictor formula} as predicts the root of the implicit method. Besides, $h$ is chosen so small such that \eqref{imeuler:Lh} is much smaller than $1$ to have a rapid fixed-point convergence. Often few iterates (sometimes only one iterate) need(s) to obtain a satisfactory result.

Another practical way is to assume $y_{k+1}^{(0)}=y_{k}$ and do the iteration \eqref{imeuler:iter} twice. This two-point iteration is equivalent to  following two-step scheme
\begin{align*}
  z & =  y_k + hf(t_{k+1},y_k) \\
  y_{k+1} & = y_k + hf(t_{k+1},z),
\end{align*}
or, by writing it in a one line,
\begin{equation}\label{imeuler:1iter}
  y_{k+1} = y_k + hf(t_{k+1},y_k+hf(t_{k+1},y_k)).
\end{equation}
This method is still of first-order accuracy but has some absolute stability limitations. However, the implicit Euler's method \eqref{imeuler:method} if is applied on test equation \eqref{lineareq_test} gives
$$
y_{k} = \frac{1}{(1-\lambda h)^k}y_0.
$$
The instability never happens if
$$
\frac{1}{|1-\lambda h|}\leqslant 1.
\vsp
$$
Therefore, the region of absolute stability of the method is $S = \{z\in \C: |1-z|\geqslant 1\}$ which is shown on the right hand side of Figure \ref{fig:eulerregabs}.

\begin{figure}[!th]
\begin{center}
\begin{tikzpicture}[scale=2]
    \node[color=black] at (-.5,1.2) {{\tiny Trapezoidal}};
    \draw[fill=blue!10] (-3/2,-1) rectangle (0.25cm,1cm);
    \draw[fill=blue!0] (-1/2,-1) rectangle (0.5cm,1cm);
    \draw[-] (-3/2,0) -- (1/2,0)  node[right] {} ;
    \node[color=black] at (-3/2,-1.1) {{\tiny $-2$}};
    \node[color=black] at (-1,-1.1) {{\tiny $-1$}};
    \node[color=black] at (-1/2,-1.1) {{\tiny $0$}};
    \node[color=black] at (0,-1.1) {{\tiny $1$}};
    \node[color=black] at (1/2,-1.1) {{\tiny $2$}};
    \node[color=black] at (-1.7,-1) {{\tiny $-2$}};
    \node[color=black] at (-1.7,-1/2) {{\tiny $-1$}};
    \node[color=black] at (-1.6,0) {{\tiny $0$}};
    \node[color=black] at (-1.6,1/2) {{\tiny $1$}};
    \node[color=black] at (-1.6,1) {{\tiny $2$}};
\end{tikzpicture}
\end{center}
\caption{Absolute stability region of the trapezoidal method}\label{fig:regionabstrap}
\end{figure}
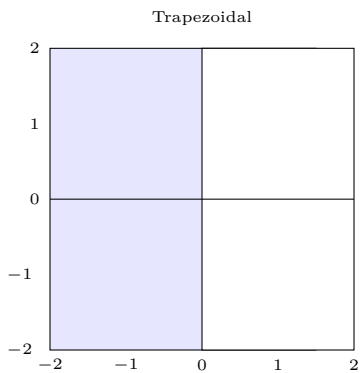

One drawback of both explicit and implicit Euler's methods is the low convergence order. Before presenting new higher order schemes let us discuss another approach for obtaining Euler's formulas. If we integrate the  equation
$y'(t)=f(t,y(t))$ from $t_k$ to $t_{k+1}$, we obtain
\begin{equation}\label{imeuler:integral}
  y(t_{k+1})=y(t_k)+\int_{t_k}^{t_{k+1}}f(\tau,y(\tau))d\tau.\vsp 
\end{equation}
The explicit Euler's method will be resulted if the integral in \eqref{imeuler:integral} is approximated by the {\em box rule}
$$
\int_{a}^{b}g(\tau) d\tau \approx (b-a)g(a)
\vsp
$$
while the implicit Euler's method follows from the box quadrature
$$
\int_{a}^{b}g(\tau) d\tau \approx (b-a)g(b).
$$
More accurate quadratures can be used to obtain more accurate methods. For instance, we can implement the {\em trapezoidal rule} (with the error term)
\begin{equation}\label{trap-rule}
\int_{a}^{b}g(\tau) d\tau = \frac{1}{2}(b-a)[g(a)+g(b)] -\frac{1}{12}(b-a)^3g''(\xi),
\end{equation}
for some $a\leqslant \xi\leqslant b$. Applying \eqref{trap-rule} to \eqref{imeuler:integral}, we obtain
\begin{equation}\label{trap1}
  y(t_{k+1})=y(t_k)+\frac{h}{2}[f(t_k,y(t_k))+f(t_{k+1},y(t_{k+1}))]-\frac{h^3}{12}y^{(3)}(\xi_k)
\end{equation}
for some $t_k\leqslant \xi_k\leqslant t_{k+1}$. The second derivative in the error term of the trapezoidal rule \eqref{trap-rule} is replaced by the third derivative of $y$ in \eqref{trap1} because $f(t,y)=y'(t)$. By dropping the error term in \eqref{trap1} and replacing
$y(t_k)$ by approximate values $y_k$, the {\bf trapezoidal method} is obtained as
\begin{shaded}
\vspace*{-0.2cm}
\begin{equation}\label{trap-method}
\begin{split}
  & y_{k+1}=y_k+\frac{h}{2}[f(t_k,y_k)+f(t_{k+1},y_{k+1})], \quad k=0,1,2,\ldots\\
  & y_0=y(t_0).
  \end{split}
\end{equation}
\vspace*{-0.2cm}
\end{shaded} 
 The local truncation error for this method is
\begin{equation}\label{trap-trunc}
  \tau_{k+1} = -\frac{h^3}{12}y^{(3)}(\xi).
\end{equation}
It can be proved that the trapezoidal method is of second-order accuracy and its global error satisfies
$$
|e_N|\leqslant Ch^2 \vsp
$$
for all sufficiently small $h$. The proof follows the same sketch as the proof of the explicit Euler's method.
In additions, we can simply show that the region of absolute stability of the trapezoidal method is the left half plane as shown in Figure
\ref{fig:regionabstrap}.

\begin{workout}
Show that the region of absolute stability of the trapezoidal method is the left half complex plane.
\end{workout}
\vsp 

The convergence order $2$ and the absolute stability of the trapezoidal method are two advantages that make this method an important tool for solving ordinary differential equations.

When $f(t,y)$ is nonlinear in $y$, the discussion for the solution of the implicit Euler's method applies to the solution of the trapezoidal method \eqref{trap-method} with a slight variation. The iteration formula \eqref{imeuler:iter} is replaced by
\begin{equation}\label{imeuler:itertrap}
  y_{k+1}^{(j+1)}=y_k + \frac{h}{2}[f(t_k,y_k)+f(t_{k+1},y_{k+1}^{(j)}),\quad j=0,1,2,\ldots.
\end{equation}
The convergence condition \eqref{imeuler:Lh} is replaced by
$$
\frac{hL}{2}\leqslant 1.
$$
The usual choice of the initial guess $y_{k+1}^{(0)}$ for \eqref{imeuler:itertrap} is the Euler's solution
$$
y_{k+1}^{(0)} = y_k+hf(t_k,y_k).
$$
With this choice the resulting global error will be still of order $h^2$, because a one step error of the Euler's method is of order $h^2$. If only one iteration of \eqref{imeuler:itertrap} is used the resulting new scheme is
\begin{equation}\label{heun-method}
  y_{k+1} = y_k + \frac{h}{2}[f(t_k,y_k)+ f(t_{k+1},y_k+hf(t_k,y_k))],
  \vsp
\end{equation}
which is also known as {\em Heun's method}. This method is still of second order accuracy, but with a more restricted region of absolute stability. The Heun's method is identical with one of the {\em Runge-Kutta methods} of order two. We will address various types of Runge-Kutta method in a forthcoming section.
\vsp 

\begin{workout}
Show that the absolute stability regions of schemes \eqref{imeuler:1iter} and \eqref{heun-method} are bounded. 
Especially both of them do not include the whole negative real line. 
\end{workout}
\vsp


\begin{workout}
Let $\theta\in[0,1]$ and consider the $\theta$-method
$$
y_{k+1}=y_k + h[(1-\theta)f(t_k,y_k)+\theta f(t_{k+1},y_{k+1})].
$$
(a) Which values of $\theta$ correspond to the explicit Euler, implicit Euler, and trapezoidal methods?
(b) Separate two cases $\theta\in[0,1/2)$
and $\theta\in [1/2,1]$. In which case the left-half plane lies in the absolute stability regions of this method?  

Optional: answer the above questions for the {\em generalized midpoint} method 
$$
y_{k+1}=y_k+hf(t_{k+\theta},(1-\theta)y_k+\theta y_{k+1})
$$
where $t_{k+\theta}=(1-\theta)t_k+\theta t_{k+1}$. 
\end{workout}
\vsp

\begin{workout}\label{wo:leapfrog}
Consider the scheme 
$$
y_{k+1}=y_{k-1}+2hf(t_{k},y_k),
$$
for solving $y'(t) = f(t,y)$. 
Show that the local truncation error of this scheme is of order $h^3$. 
The absolute stability region for this scheme is
$S = \{z=\alpha+i\beta\in\C: \alpha=0, -1\leqslant\beta\leqslant 1\}$, which is a marginal stability region (no need to prove!)
This method is known as {\em midpoint} or {\em leapfrog method}. Is the leapfrog method $A$-stable?

Hint: To obtain the truncation error, replace $y_k$ by exact values $y(t_k)$ and leave the formula by additional term $\tau_k$. Then use Taylor expansion to determine $\tau_k$. 
\end{workout}
\vsp
\begin{labexercise}
Write MATLAB functions for implicit and trapezoidal methods. In each case, use an iterative method and/or a predictor formula to handle the nonlinearity.
Then solve the examples in Lab Exercise \ref{py_ex_euler} using your codes and compare the errors and orders. Comment on your results.
\end{labexercise}
\vsp

\begin{labexercise}\label{wo_eulersmoothness}
Consider the IVP
$$
y'(t)=\frac{1}{t}y(t)+(\alpha-1) t^{\alpha-1}, \quad y(0)=0, \quad t\in[0,1],
$$
with $\alpha>0$. The solution is $y(t)=t^\alpha$. To have $y$ twice continuously differentiable, we need $\alpha\geqslant 2$.
Use your MATLAB codes for explicit and implicit Euler and trapezoidal methods  for $\alpha=2.5,1.5,1.1$ with stepsizes $h=0.2,0.1,0.05,0.025,0.0125$. Determine the computational convergence orders. Compare with theoretical orders and report a reason for your observation.
\end{labexercise}

%% file: lec3_part2.tex
\subsection{Taylor series methods}
Euler's methods can be formulated by using a Taylor series approximation when $y'$ is replaced by $f(t,y)$ and higher order derivative in the series are dropped. One can of course use higher order terms but then $y''$, $y'''$, $\ldots$ should be obtained by differentiating the differential equation
$$
y'(t)=f(t,y),
$$
successively. From the Taylor series expansion of order $p$ we have
\begin{equation}\label{tayl-ser}
  y(t_{k+1}) \approx y(t_k) + hy'(t_k)+\frac{h^2}{2!}y''(t_k)+\cdots + \frac{h^p}{p!}y^{(p)}(t_k)
\end{equation}
where the truncation error is
\begin{equation}\label{tayl-trunc}
  \tau_{k+1} = \frac{h^{p+1}}{(p+1)!}y^{(p+1)}(\xi_k), \quad t_k\leqslant \xi_k\leqslant t_{k+1}.
\end{equation}
The term $y'(t)$ in \eqref{tayl-trunc} can be replaced by $f(t,y)$ as we have done in Euler's methods. For higher order derivatives we can write
\begin{equation*}
  \begin{split}
      y''(t)& = f_t + f_yf\\
      y^{(3)}(t)& = f_{tt}+2f_{ty}f + f_{yy}f^2+f_y(f_t+f_yf)\\
      \vdots
  \end{split}
\end{equation*}
provided that partial derivatives of $f(t,y)$ with respect to $y$ exist.
Substituting these formulas into \eqref{tayl-ser}, we obtain
\begin{equation}\label{tayl-method}
  y_{k+1} =y_k + hy'_k + \frac{h^2}{2}y''_k + \cdots + \frac{h^p}{p!}y^{(p)}_k,
  \vsp 
\end{equation}
which is called the {\bf Taylor series method}. The derivatives formulas in \eqref{tayl-method} are
$$
y'_k = f(t_k,y_k), \quad y''_k = (f_t+f_yf)(t_k,y_k), \; \mbox{and so on}.
$$
The formulas for higher order derivatives rapidly become too
complicated, so Taylor series methods of higher order
have not often been used in practice. Recently, however, the availability of symbolic
manipulation and automatic differentiation systems have made these methods more
feasible.

If the solution $y$ and the derivative function $f(t,y)$ are sufficiently differentiable then it can be proved that the global error for the scheme \eqref{tayl-method} satisfies
\begin{equation*}
  |e_N|\leqslant Ch^p\|y^{(p+1)}\|_\infty,
\end{equation*}
which means that the method is of $p$-th order accuracy. The constant $C$ is something similar to that was obtained for the explicit Euler's method (the first order Taylor method).

\begin{labexercise}
Construct the Taylor series methods of orders $2$ and $3$ for IVP
$$
y'(t) = [\cos y(t)]^2, \quad 0\leqslant t\leqslant 10, \quad y(0) = 0.
$$
 Write a MATLAB code to compute the results for stepsizes $h=0.2$, $0.1$,$0.05$, $0.025$, $0.0125$, $0.00625$.
Plot the error functions in each case and calculate numerical orders.
Compare the results with Euler and trapezoidal methods. 
The exact solution of the above IVP is $y(t)=\tan^{-1}(t)$. Use this information to calculate errors and orders. 
\end{labexercise}
\vsp

\subsection{Runge-Kutta methods}
The calculation of higher order partial derivatives of $f(t,y)$ makes the Taylor series methods complicated and time-consuming. {\bf Runge-Kutta} methods (abbreviated by RK methods) replace higher derivatives by more evaluations of $f(t,y)$ to have finite difference
approximations for derivatives while retain the accuracy of Taylor series methods. The RK methods are {\em one-step} but {\em multi-stage} and are fairly easy to program not only for a scaler ODE but also for a system of ODEs.

To derive a second order RK method, consider the second order Taylor method
\begin{equation}\label{rk-tayl}
y_{k+1}=y_k + hy'_k + \frac{h^2}{2}y''_k
\end{equation}
where $y' = f(t,y)$ and $y''=f_t + f_yf$ both evaluated at $(t_k,y_k)$. We aim to approximate $f_t+f_yf$ by expanding $f$ in a bivariate Taylor series as
$$
f(t+h,y+hf) = \big(f+ hf_t+hff_y\big){(t,y)} + \mathcal O(h^2).
$$
This simply shows that
$$
y''(t)=(f_t + f_yf)(t,y)  =  \frac{1}{h}[f(t+h,y+hf(t,y))-f(t,y)] + \mathcal O(h).
$$
Dropping the $\mathcal O(h)$ term and substituting in \eqref{rk-tayl}, we obtain
\begin{equation}\label{rk-heun}
\begin{split}
  y_{k+1} & = y_k + hf(t_k,y_k) + \frac{h^2}{2}\frac{1}{h}[f(t_{k}+h,y_k+hf(t_k,y_k))-f(t_k,y_k)] \\
   & = y_k + \frac{h}{2}[f(t_k,y_k) + f(t_{k}+h,y_k+hf(t_k,y_k))].
\end{split}
\end{equation}
This method was previously derived as the Heun's method in \eqref{heun-method}. As an RK2 method it is usually written in the following two step pattern:
\begin{shaded}
\vspace*{-0.3cm}
\begin{equation}\label{rk2}
\begin{split}
  z_1 & = y_k\\
  z_2 & =y_k+hf(t_k,z_1) \\
  y_{k+1} & = y_k + \frac{h}{2}[f(t_k,z_1)+f(t_k+h,z_2)].
\end{split}
\end{equation}
\vspace*{-0.2cm}
\end{shaded} 
This is not the only order $2$ explicit RK method. As we discussed in section \ref{sect:general-explicit}, a general explicit method can be written as
\begin{equation*}
  y_{k+1} = y_k + h\psi(t_k,y_k,h), \quad y_0=y(t_0).
\end{equation*}
In the RK2 method \eqref{rk-heun} we derived $\psi(t,y,h)$ as
$$
\psi(t,y,h) = \frac{1}{2}f(t,y) + \frac{1}{2}f(t+h,y+hf(t,y)).
$$
This formula can be generalized to ansatz
\begin{equation}\label{rk2:ansatz}
  \psi(t,y,h) = b_1 f(t,y) + b_2 f(t+\alpha h,y+\beta hf(t,y)),
\end{equation}
with unknown coefficients $b_1,b_2, \alpha, \beta$ that can be determined such that the local truncation error
\begin{equation*}
  \tau_{k+1} = y(t_{k+1})-[y(t_k)+h\psi(t_k,y(t_k),h)] \vsp
\end{equation*}
will satisfy $\tau_{k+1} = \mathcal{O}(h^3)$ just as with the Taylor method of order $2$. After some manipulations with the bivariate Taylor expansion, we will obtain the relations  
\begin{equation*}
  b_2\neq 0, \quad b_1=1-b_2, \quad \alpha = \beta = \frac{1}{2b_2}.\vsp
\end{equation*}
between the coefficients in order to have $\tau_{k+1} = \mathcal{O}(h^3)$. 
Depending on the choice of $b_2$, there exists a family of RK methods of order $2$. The case $b_2=1/2$ results in RK method \eqref{rk-heun}. Another choice $b_2=1$, $b_1=0$ and $\alpha=\beta = \frac{1}{2}$, results in
\begin{equation}\label{rk2-2}
  y_{k+1} = y_k + hf(t_k+\tfrac{1}{2}h, y_k + \tfrac{1}{2}hf(t_k,y_k)).
\end{equation}
or in a multi-stage format
\begin{shaded}
\vspace*{-0.3cm}
\begin{equation}\label{rk2_3}
\begin{split}
  z_1 & = y_k\\
  z_2 & = y_k+\frac{h}{2}f(t_k,z_1) \\
  y_{k+1} & = y_k + hf(t_k+\tfrac{h}{2},z_2).
\end{split}
\end{equation}
\vspace*{-0.2cm}
\end{shaded} 
A general explicit RK method is defined as below. 
\begin{definition}
An explicit RK method with $s$ stages has the form
\begin{equation}\label{rk:sstageform}
  \begin{split}
     z_1 &= y_k ,\\
     z_2 &= y_k + ha_{2,1}f(t_k,z_1) ,\\
     z_3 &= y_k + h\big[a_{3,1}f(t_k,z_1)+a_{3,2}f(t_k+c_2h,z_2)\big] ,\\
       & \vdots \\
     z_s& =y_k + h\big[a_{s,1}f(t_k,z_1)+a_{s,2}f(t_k+c_2h,z_2)+\cdots+a_{s,s-1}f(t_k+c_{s-1}h,z_{s-1})\big], \\
     y_{k+1} &= y_k + h\big[b_{1}f(t_k,z_1)+b_{2}f(t_k+c_2h,z_2)+\cdots+b_{s}f(t_k+c_{s}h,z_{s})\big].
  \end{split}
\end{equation}
\end{definition}
Such RK method is fully determined by coefficients $\{c_\ell,a_{\ell,j},b_j\}$. These coefficients are usually displayed in a table called {\bf Butcher's Tableau}\footnote{After John Charles Butcher (1933-present) who is a New Zealand mathematician and an specialist in numerical methods for ODEs.}

\begin{equation*}
  \begin{array}{r|ccccc}
    0=c_1 &  &  &  &  &  \\
    c_2 & a_{2,1} &  &  &  &  \\
    c_3 & a_{3,1} & a_{3,2} &  &  &  \\
    \vdots &\vdots & \vdots & \ddots &  &  \\
    c_s & a_{s,1} & a_{s,2} & \cdots & a_{s,s-1} &  \\ \hline
     & b_1 & b_2 & \cdots & b_{s-1} & b_s
  \end{array}
  \vsp 
\end{equation*}
RK methods can be expressed in the general form \eqref{onestep:general} with
$$
\psi(t,y,h) = \sum_{j=1}^{s}b_{j} f(t+c_jh,z_j), \quad z_j = y+h\sum_{\ell=1}^{j-1}a_{j\ell}f(t+c_\ell h,z_\ell).
$$
The consistency condition $\psi(t,y,0)=f(t,y)$ holds if
\begin{align}\label{rk:consistency}
  \sum_{j=1}^{s} b_{\ell}  = 1.
\end{align}
According to Theorem \ref{thm:converg-onesteps}, if $f$ is continuous in $t$ and Lipschitz continuous in $y$, and
the condition \eqref{rk:consistency} holds, then the RK method \eqref{rk:sstageform} is convergent.

Moreover, in a RK method we always assume that
\begin{align}\label{rk:stagecond}
  \sum_{j=1}^{\ell-1} a_{\ell j} =c_\ell, \quad \ell=1,2,\ldots, s,
\end{align}
which ensure that intermediate values $z_\ell$ provide approximations of order at least $1$ to exact values $y(t_k+c_\ell h)$.
Conditions \eqref{rk:stagecond} are called the {\em stage conditions} for RK methods.
\vsp 
\begin{example}
The Butcher's tableau of the RK2 method \eqref{rk2} is
\begin{equation*}
  \begin{array}{r|cc}
    0 &  &    \\
    1 & 1 &   \\\hline
     & 1/2 & 1/2
  \end{array}
\end{equation*}
while the RK2 method \eqref{rk2_3} has a Butcher's tableau of the form
\begin{equation*}
  \begin{array}{r|cc}
    0 &  &    \\
    1/2 & 1/2 &   \\\hline
     & 0 & 1
  \end{array}
\end{equation*}
\end{example}
\vsp 
There also exist a family of third-order RK methods. The Butcher's tableau of one of these schemes is: 
\begin{equation*}
  \begin{array}{r|ccc}
    0 &  &  &  \\
    1/2 & 1/2 &&   \\
    1   &-1 & 2& \\ \hline
      &1/6&2/3&1/6
  \end{array}\vsp
\end{equation*}
\begin{workout}
(a) Convert the above tableau into a 3-stage RK formula. 
(b) Search the internet to find another RK3 method and write down its stages. 
\end{workout}
A very popular explicit RK method is the following fourth-order scheme (RK4):
\begin{shaded}
\vspace*{-0.3cm}
\begin{equation}\label{rk4}
  \begin{split}
     z_1 &= y_k ,\\
     z_2 &= y_k + \tfrac{1}{2}hf(t_k,z_1) ,\\
     z_3 &= y_k + \tfrac{1}{2}hf(t_k+\tfrac{1}{2}h,z_2) ,\\
     z_4 &= y_k + hf(t_k+\tfrac{1}{2}h,z_3) ,\\
 y_{k+1} &= y_k + \frac{h}{6}\big[f(t_k,z_1)+2f(t_k+\tfrac{1}{2}h,z_2)+2f(t_k+\tfrac{1}{2}h,z_3)+f(t_k+h,z_4)\big].
  \end{split}
\end{equation}
\vspace*{-0.2cm}
\end{shaded} 
The Butcher's tableau for this method is
\begin{equation*}
  \begin{array}{r|cccc}
    0 &  &  &&  \\
    1/2 & 1/2 &&&   \\
    1/2 & 0 & 1/2&& \\
    1 &0&0&1& \\ \hline
      &1/6&1/3&1/3&1/6
  \end{array}
\end{equation*}
Using a similar but more tedious calculation than that was done for method \eqref{rk2}, we can show that the local truncation error for this $4$-stage method is of order $h^5$. Then by applying the sketch given in section \ref{sect:general-explicit} for analysing general explicit methods we can show that the global error of method \eqref{rk4} satisfies
$$
|e_N|\leqslant Ch^4
$$
which shows that this method is of fourth-order accuracy, for this reason is usually called RK4. This is not the only fourth-order explicit RK method; there exists a family of such methods with different Butcher's tableaus.
\vsp 

\begin{workout}
Show that the RK4 method \eqref{rk4} when is applied on simple differential equation $y'(t)=f(t)$ (with no dependence of $f$ on $y$)
reduces to Simpson's rule for numerical integration.
\end{workout}
\vsp 

\begin{figure}[!th]
  \centering
  \includegraphics[scale=0.5]{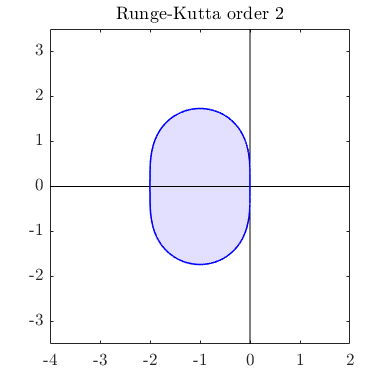}\includegraphics[scale=0.5]{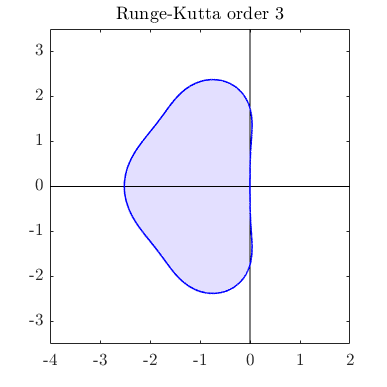}\includegraphics[scale=0.5]{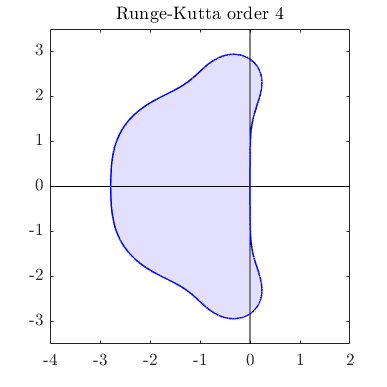}\\
  \caption{Absolute stability regions of Runge-Kutta methods of orders $2$, $3$ and $4$.}\label{fig:absstab-rk}
\end{figure}

Like as other explicit methods, RK methods have restricted regions of absolute stability. For example, by applying the RK2 formula \eqref{rk-heun} on test equation $y'(t)=\lambda y$, we obtain
$$
y_k = (1+\lambda h + \frac{1}{2}(\lambda h)^2) y_k,
$$
which shows that the region of absolute stability for this method is
$$
S = \{z\in \C: |1+z+\tfrac{1}{2}z^2|\leqslant 1\}.
$$
This region is shown in Figure \ref{fig:absstab-rk} (left panel). The stability regions for an RK3 method and the RK4 method \eqref{rk4},
obtained in a similar way, are shown in the middle and right sides of Figure \ref{fig:absstab-rk}.

\begin{labexercise}
Use a MATLAB code for solving the IVP
$$
y'(t) = \frac{1}{1+t^2}-2[y(t)]^2, \quad 0\leqslant t\leqslant 10, \quad y(0) = 0.
$$
using RK2 and RK4 formulas with stepsizes $h=0.2,0.1,0.05,0.025,0.0125,0.00625$.
Plot the error functions in each case and determine numerical orders of convergence.
The exact solution for this IVP is $y(t)=t/(1+t^2)$. Use this information to calculate errors and orders.
\end{labexercise}
\vsp

\begin{labexercise}
Use the RK2 method to solve
$$
y'(t) = -y(t) + t^{0.1}(1.1+t), \quad y(0)=0
$$
whose exact solution is $y(t)=t^{1.1}$. Solve the equation on interval $[0,5]$ and compute the solution and errors at times $t=1,2,3,4,5$. 
Use different stepsizes $h=0.1,0.05,0.025,0.0125,0.00625$. Compute the order of convergence and compare with theoretical order $2$ of the RK2 method. Explain your results. 
\end{labexercise}
\vsp 

\begin{workout}
What difficulty arises in attempting to use a Taylor series method of order $\geqslant 2$ to solve the equation 
$$
y'(t) = -y(t) + t^{0.1}(1.1+t), \quad y(0)=0. 
$$
What does it tell us about the solution? 
\end{workout}
\vsp 
\begin{workout}
(a) Write down the RK4 method for solving linear system of equations $\y'(t)=A\y(t)$ for $A\in\R^{n\times n}$ with initial condition $\y(0)=\y_0$. (b) How many matrix-vector multiplications should be performed in each step? (c) Optional: estimate the overall complexity of the method to approximate $\y(t_N)$. 
\end{workout}
\vsp 

\section{Stiff differential equations}\label{sect:stiff}
At the beginning of 1950's, a new difficulty was discovered in numerical solution of some practical ODEs, which has come to be known
as {\bf stiffness}, and led to some new concepts of {stability} by Germund Dahlquist\footnote{Germund Dahlquist (1925--2005) was a Swedish mathematician known primarily for his early contributions to the theory of numerical solution of ODEs.} and others.
Numerical methods with finite absolute stability regions (such as explicit methods) all fail to produce accurate and stable solutions for stiff problems unless the step size $h$ is chosen excessively small which is impractical and inefficient in many situations.
\subsection{What is stiffness?}
It is difficult to formulate a precise definition for stiffness. One may argue that a stiff equation includes some terms that can lead to rapid variation (fast transients) in the solution. However, there exist some ODEs with smooth solutions\footnote{In this section, by a `smooth solution' we mean a function without any rapid transition.} but are known as stiff problems. This means that the stiffness is independent of the solution but it is a property of the ODE itself. However, even if a stiff problem has a smooth solution,
a slight perturbation to the solution at any time results in another solution curve that has a rapid variation.
The following example from \cite{LeVeque:2007} will make this more clear.
\begin{example}\label{ex:stiff1}
Consider the ODE
\begin{equation}\label{eq:stiff1}
y'(t) =\lambda(y-\cos t)-\sin(t) \vsp
\end{equation}
from Example \ref{ex:absstab}.
For initial value $y(0)=1$,
this problem has smooth solution $y(t)=\cos t$ independent of the value of $\lambda$.
If we change the initial condition to $y(t_0)=y_0$ that does not lie on this curve, then the solution is
$$
y(t) = e^{\lambda(t-t_0)}(y_0-\cos t_0) + \cos t
$$
If $\Re(\lambda)<0$, this function approaches $\cos t$ exponentially quickly with decay rate $\Re(\lambda)$.
In Figure \ref{fig:stiff2} different solution curves for this equation with different
choices of $t_0$ and $y_0$ with $\lambda=-2$ and $\lambda=-20$ are plotted.
\begin{center}
\includegraphics[scale=0.68]{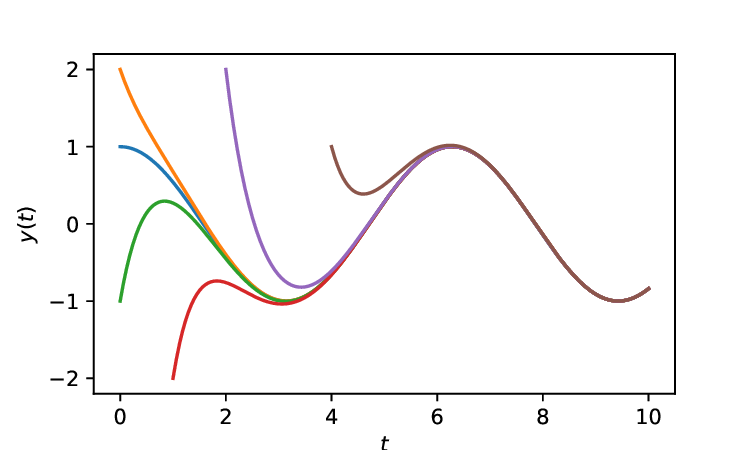}\includegraphics[scale=0.68]{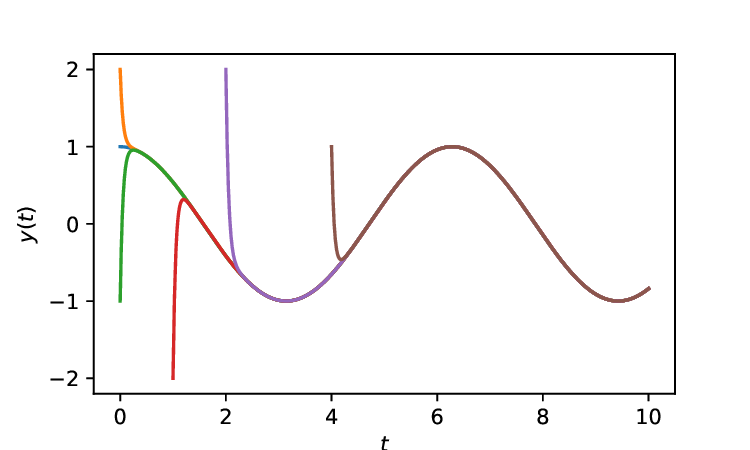}
\captionof{figure}{Solution curves of a stiff problem with different initial times and initial values for $\lambda=-2$ (left) and $\lambda=-20$ (right).}\label{fig:stiff2}
\end{center}
We observe rapid transients in solution curves for $\lambda=-20$.
The perturb solutions quickly approach
toward the particular solution $y(t)=\cos t$. This problem is known as a stiff problem for values of $\lambda$ with large real part magnitudes.
\end{example}
\vsp 

The phenomenon we observed in Example \ref{ex:stiff1} will cause a serious numerical difficulty even if the initial condition is chosen such that the exact solution does not exhibit any rapid transient (for example $y(t)=\cos t$ with $y(0)=1$ in ODE \eqref{eq:stiff1}).
{\em Because any numerical method is subjected to local truncation and roundoff errors which act as a perturbation to the solution and move us away from the smooth solution to a solution with a rapid transient.}
Numerical methods with finite absolute stability regions are unstable unless the time step is small
relative to the time scale of the rapid transient.
In the case of a smooth true solution it
seems that a reasonably large step length would work, but the numerical method must
always deal with the rapid transients introduced by truncation and roundoff errors in every time step. Consequently,
a very small step length  is needed to avoid the instability.
\vsp 

\begin{example}
Consider the ODE \eqref{eq:stiff1} on interval $[0,10]$ with initial condition $y(0)=1$. Let $\lambda=-10^4$.
The numerical solution using the explicit RK4 method \eqref{rk4} with step length  $h=0.00028$ blows up as is shown in the left hand side of  Figure \ref{fig:stiff3}. Smaller values of $h$ such as $h = 0.00025$ leads to a stable solution that is shown in the right hand side. However, this stable calculations requires lots of function evaluations in the procedure of RK method.
Note that, the complexity will quickly increase for a system of differential equations.
\begin{center}
\includegraphics[scale=0.68]{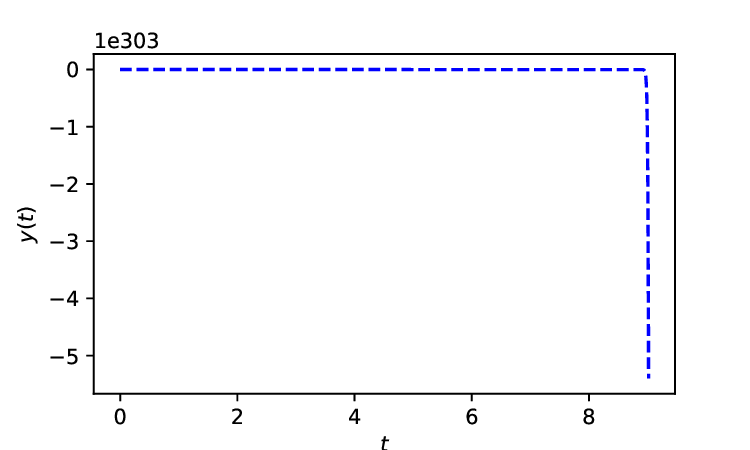}\includegraphics[scale=0.68]{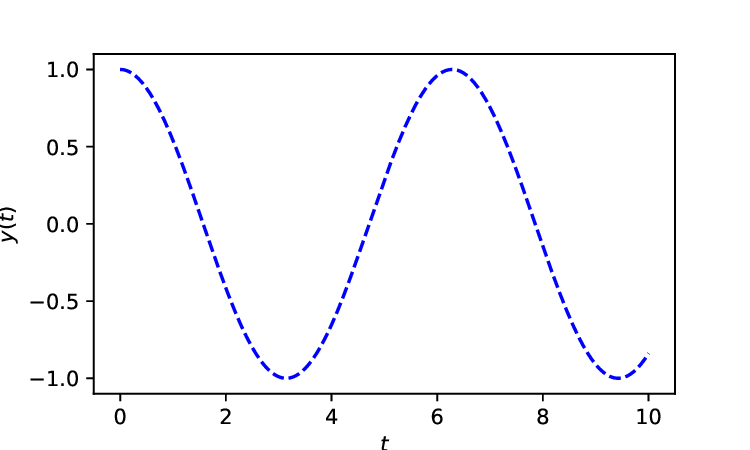}
\captionof{figure}{Numerical solution of ODE \eqref{eq:stiff1} with $\lambda=-10^4$ and initial condition $y(0)=1$ using the RK4 method: unstable solution with $h=0.00028$ (left) and a stable solution with $h=0.00025$ (right). In the left panel, the values on the $y$-axis are multiplied by huge number $10^{303}$.}\label{fig:stiff3}
\end{center}

It is better to look for other efficient numerical methods
that can solve such a stiff problem stably using a much larger step size. For instance, the implicit Euler's method will do this job perfectly (write and execute the code).
The trapezoidal method works much better than explicit methods but still introduces some limitations in the presence of rapid transients in the solution.
\end{example}
\vsp

Lambert, after examining some other statements, finally has suggested the following definition for stiff problems \cite{Lambert:1991}.
\begin{definition}\label{def:stiffness}
A stiff ODE is an equation for which certain numerical methods with finite absolute stability regions for solving the equation are numerically unstable, unless the step size is taken to be excessively small in relation to the smoothness of the exact solution.
\end{definition}
\vsp 

Explicit methods such as forward Euler's method and explicit RK methods (and in fact all explicit methods) are examples of numerical methods with finite absolute stability regions. Definition \ref{def:stiffness} reveals the fact that solving a stiff problem with an explicit method
(or an implicit method with a finite absolute stability region) is very costly.
We also note that the stiffness may vary over the total interval of integration.
\subsection{Stiff systems}
Now, let us look at a system of ODEs. Consider the linear system $\y'(t)=A\y(t)+\g(t)$ on interval $[0,b]$ with a constant matrix $A$ of size $n\times n$, and initial condition $\y(0)=\y_0$. The exact solution of this ODE is
$$
\y(t) = e^{At}{\y_0}  + \int_{0}^{t}e^{A(t-\tau)}{\g(\tau)}d\tau.
$$
If $A$ has distinct eigenvalues $\lambda_k\in \C$ and corresponding eigenvectors $\v_k\in \C^n$, then
$$
\y(t) = \sum_{k=1}^{n}c_{k}e^{\lambda_k t}\v_k + \wt \g(t)
$$
where $\c=V^{-1}\y_0$ for $V=[\v_1,\ldots,\v_n]$, and $\wt \g(t)$ is the integral term (particular solution) in solution $\y(t)$.
Let us assume that
$$
\Re(\lambda_k)< 0 ,\quad  k=1,2,\ldots,n,
$$
which imply that all terms $e^{\lambda_k t}\v_k$ go to $0$ as $t\to\infty$, so that the solution $\y$ approaches $\wt \g(t)$ asymptotically as $t\to\infty$. The term $e^{\lambda_k t}\v_k$ decreases monotonically  if $\lambda_k$ is real and sinusoidally if $\lambda_k$ is complex.
Thus, the term
$$
\sum_{k=1}^{n}c_{k}e^{\lambda_k t}\v_k \vsp
$$
can be viewed as the {\em transient solution} and the term $\wt \g(t)$ as the {\em steady state} solution.
If $\Re(\lambda_k)$ is large then the corresponding term $c_ke^{\lambda k t}\v_k$
will decay quickly as $t$ increases and is thus called a fast transient; if $\Re(\lambda_k)$ is small
the corresponding term $c_ke^{\lambda k t}\v_k$ decays slowly and is called a slow transient.
Let $\overline\lambda$ and $\underline \lambda$ be defined such that
$$
|\Re (\underline \lambda)|\leqslant |\Re (\lambda_k)|\leqslant |\Re(\overline\lambda)|, \quad k=1,2,\ldots n.
$$
If our aim is to reach
the steady-state solution, then we must keep integrating until the slowest transient is
negligible. The smaller $|\Re (\underline \lambda)|$ is, the longer we must keep integrating. If, however, the
method we are using has a finite region of absolute stability, we must ensure that the step size $h$ is sufficiently small for
$\lambda_kh\in S$, $k = 1,2,\ldots,n$ to hold. Clearly a large value of $|\Re (\overline \lambda)|$ implies a small step size $h$.
We therefore get into a difficult situation if $|\Re (\overline \lambda)|$ is very large and $|\Re (\underline \lambda)|$ is very small;
we are forced to integrate for a very long time with an excessively small step size. It seems natural to take the ratio
\begin{equation*}
r_S:=\frac{|\Re (\overline \lambda)|}{|\Re (\underline \lambda)|}
\end{equation*}
the {\em stiffness ratio}, as a measure of the stiffness of the system \cite{Lambert:1991}.
While $r_S$ is
often a useful quantity, one should not rely entirely on this measure to determine whether
a problem is stiff. For example, a scalar problem can be stiff while the $r_S$ is always $1$ since there
is only one eigenvalue. Or, in a system of equations if one eigenvalue is zero then
the contribution of this eigenvalue to the exact solution is a constant term.
If the moduli of the real parts of the remaining eigenvalues are not particularly large,
the system is not stiff, yet the stiffness ratio is now infinite.
\vsp 

\begin{example}\label{ex:stiff2}
Consider the following systems of ODEs \cite{Lambert:1991}
\begin{equation}\label{eq:stiff-sys1}
  \begin{bmatrix}
    y'_1 \\
    y'_2
  \end{bmatrix}
=
  \begin{bmatrix}
    -2 &1 \\
    1  &-2
  \end{bmatrix}
    \begin{bmatrix}
    y_1 \\
    y_2
  \end{bmatrix}
+
  \begin{bmatrix}
    2\sin t \\
    2(\cos t-\sin t)
  \end{bmatrix}
,\quad
    \begin{bmatrix}
    y_1(0) \\
    y_2(0)
  \end{bmatrix}
=
    \begin{bmatrix}
    2 \\
    3
  \end{bmatrix}
\end{equation}
and
\begin{equation}\label{eq:stiff-sys2}
  \begin{bmatrix}
    y'_1 \\
    y'_2
  \end{bmatrix}
=
  \begin{bmatrix}
    -2 &1 \\
    998  &-999
  \end{bmatrix}
    \begin{bmatrix}
    y_1 \\
    y_2
  \end{bmatrix}
+
  \begin{bmatrix}
    2\sin t \\
    999(\cos t-\sin t)
  \end{bmatrix}
,\quad
    \begin{bmatrix}
    y_1(0) \\
    y_2(0)
  \end{bmatrix}
=
    \begin{bmatrix}
    2 \\
    3
  \end{bmatrix}.\vsp
\end{equation}
Both problems have the same exact solution
$$
    \begin{bmatrix}
    y_1(t) \\
    y_2(t)
  \end{bmatrix}
=2e^{-x}
    \begin{bmatrix}
    1 \\
    1
  \end{bmatrix}
+
    \begin{bmatrix}
    \sin t \\
    \cos t
  \end{bmatrix}.
$$
The plots of these solutions on interval $[0,10]$ are given in Figure \ref{fig:stiff4}.
Both solutions are smooth (without any rapid transients).
\begin{center}
\includegraphics[scale=0.7]{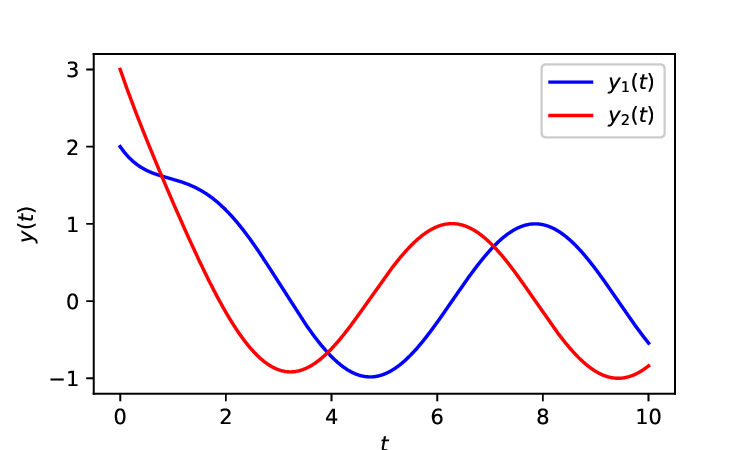}
\captionof{figure}{Exact solutions of systems \eqref{eq:stiff-sys1} and \eqref{eq:stiff-sys2}}\label{fig:stiff4}
\end{center}
We employ the explicit Euler's method for both systems. Everything is perfect for system \eqref{eq:stiff-sys1}, but unstable results are obtained for system \eqref{eq:stiff-sys2} unless the step size is chosen smaller than $0.002$. System \eqref{eq:stiff-sys2} is stiff but system \eqref{eq:stiff-sys1} is non-stiff. The eigenvalues of the matrix in \eqref{eq:stiff-sys1} are
$-1$ and $-3$, and if we consider  the general initial condition to $y(0)=y_0$ then the exact solution is
$$
    \begin{bmatrix}
    y_1(t) \\
    y_2(t)
  \end{bmatrix}
=c_1e^{-t}
    \begin{bmatrix}
    1 \\
    1
  \end{bmatrix}
+
c_2e^{-3t}
    \begin{bmatrix}
    1 \\
    -1
  \end{bmatrix}
+
    \begin{bmatrix}
    \sin t \\
    \cos t
  \end{bmatrix}.
$$
where $c_1$ and $c_2$ are determined by imposing the initial value $y_0$. The eigenvalues of the matrix in \eqref{eq:stiff-sys2} are $-1000$ and $-1$ and the exact solution for an arbitrary initial value $y_0$ is
$$
    \begin{bmatrix}
    y_1(t) \\
    y_2(t)
  \end{bmatrix}
=c_1e^{-t}
    \begin{bmatrix}
    1 \\
    1
  \end{bmatrix}
+
c_2e^{-1000t}
    \begin{bmatrix}
    1 \\
    -998
  \end{bmatrix}
+
    \begin{bmatrix}
    \sin t \\
    \cos t
  \end{bmatrix}.\vsp
$$
In the second solution the exponential term $e^{-1000t}$ produces a rapid transient in the solution.
Although the initial value $y_0=[2,3]^T$ gives $c_1=2$ and $c_2=0$ for both systems (therefore annihilates the rapid transient in the second system), slight perturbations (truncation and roundoff errors) in numerical solution at different time levels produce a component on the rapid transient term and introduce the mentioned numerical difficulties.
The same happens even if we choose initial conditions such that $c_1=c_2=0$; in that case
the explicit Euler's and explicit RK methods are unable to integrate even the very smooth solution $y(x) = [\sin t, \cos t]^T$ unless at very small stepsizes. Finally, we note that the stiffness ratio is $r_S=3$ for system \eqref{eq:stiff-sys1} while it is $r_S=1000$ for system \eqref{eq:stiff-sys2}.
\end{example}
\vsp

The situation would be more complicated for nonlinear system of equations.
For the scaler case the partial derivative $\frac{\partial f}{\partial y}$ or Lipschitz constant
of $f$ with respect to $y$, and for nonlinear system
the eigenvalues of the
Jacobian matrix $J_f(t,\y)$ may give an insight to discover the stiffness. As we pointed out at the end of subsection \ref{sect-stabsol}, the stability analysis through the Jacobian matrix has only a local validity. Nevertheless, the stiffness ratio
$$
r_S = \max_{t\in[t_0,b]} \frac{\max_{k}|\Re(\lambda_k(t))|}{\min_{k}|\Re(\lambda_k(t))|},
$$
where $\lambda_k(t)$ are eigenvalues of $J_f(t,\y)$,
may give an insight to stiffness of system $\y'(t)=f(t,\y)$.

\begin{remark}\label{remark_stiff0}
Sometimes we can indicate the stiffness of system of ODEs $\y'(t)=f(t,\y)$ by looking at the coefficients on the right-hand side.  
When these coefficients vary significantly in magnitude, it suggests that the system is likely stiff. For example the following system of equations for 
Robertson's auto-catalytic chemical reaction
\begin{align*}
& y_1' = -0.04y_1+10^4y_2y_3\\
& y_2' = 0.04y_1-10^4y_2y_3-3\times 10^7y_2^2\\
&y_3'= 3\times 10^7 y_2^2
\end{align*}
is a stiff system of ODEs, as it contains coefficients with different sizes ranging from $0.04$ to $3\times 10^7$. This is also the case for system \eqref{eq:stiff-sys2} in Example \ref{ex:stiff2}. 
\end{remark}
\vsp

A practical way to detect the stiffness of an ODE
is to attempt to solve it using a method with a finite
absolute stability region with a moderate step size and see whether the computed solution is blown up or is not.
\vsp
\begin{workout}
Consider solving Robertson's auto-catalytic chemical system of ODEs, as described in Remark \ref{remark_stiff0}, over the time interval $[0,500]$ with initial conditions $y_1(0)=1$, $y_2(0)=y_3(0)=0$, using the Explicit Euler and Explicit RK4 methods. Ensure that the chosen step length $h$ is sufficiently small to yield stable results. Upon plotting the solutions, provide your observations. What distinguishes the behavior of the solution $y_2$ from the others? 
\end{workout}
\vsp 
\subsection{A-stability}
It is clear from the considerations of this section that those methods with bounded stability regions are inappropriate for stiff problems. Such class of methods includes all explicit methods and even some implicit methods such as Adams-Moulton methods. On the other hand, an implicit method that its
region of absolute stability includes the whole of the left half-plane is efficient for stiff problems because there
will be no stability restriction on the step size $h$ provided that all the eigenvalues have
negative real parts, as is often the case in practice.

\begin{definition}
A numerical ODE solver is called {\bf A-stable} if its absolute stability region contains the whole left-half plane $\{z\in\C:\Re(z)\leqslant0\}$.
\end{definition}
Among the methods we have already studied, the implicit Euler's method and the trapezoidal method are A-stable.
As a negative result, the {\em Dahlquist's second barrier theorem} states that any A-stable linear multistep method is at most second order accurate, and in fact the trapezoidal method is the A-stable method with smallest truncation error.
However, higher order A-stable implicit RK methods (multi-stage methods) do exist.

For many stiff problems the eigenvalues are located near the negative real axis or at most
in a wedge $\pi-\alpha\leqslant \mathrm{arg}(z)\leqslant \pi+\alpha$ for an angle $\alpha\in [0,\pi/2]$.
For such problems, we just need the stability region to contain such a wedge rather than the whole
left-half plane. See Figure \ref{fig:Aalpha}.
\begin{definition}
A numerical ODE solver is called {\bf A($\alpha$)-stable}, for $\alpha\in[0,\pi/2]$, if its absolute stability region contains the wedge $\{z\in\C:
\pi-\alpha\leqslant \mathrm{arg}(z)\leqslant \pi+\alpha\}$.
\end{definition}

\definecolor{qqqqff}{rgb}{0.8,.8,1}
\begin{center}
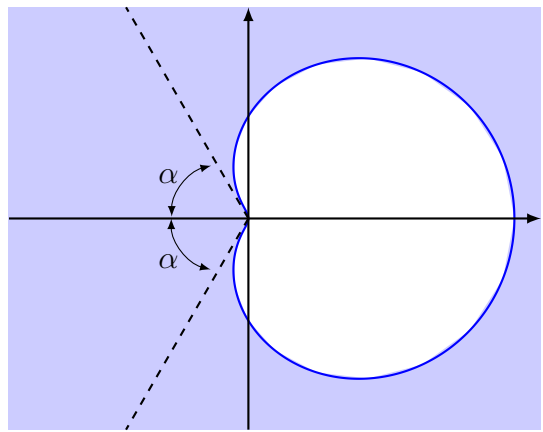

\begin{tikzpicture}[>=latex,scale=0.7]
\fill[fill=qqqqff] (-5, -4) rectangle (5,4);
\fill[fill=white] plot[domain=-pi:pi] (xy polar cs:angle=\x r,radius= {5/2+2*cos(\x r)});
\draw[thick,color=blue,domain=0:2*pi,samples=500,smooth] plot (xy polar cs:angle=\x r,radius= {5/2+4/2*cos(\x r)});
\draw[dashed,thick] (-.5,0) -- (-2.8,-4);
\draw[dashed,thick] (-.5,0) -- (-2.8,4);
\draw[<->] (-1.2,1) arc (105:180:1) ;
\draw[<->] (-1.95,0) arc (180:255:1) ;
\node at (-2,0.8) {$\alpha$};
\node at (-2,-0.8) {$\alpha$};
\draw[->,thick] (-5,0) -- (5,0);
\draw[->,thick] (-0.5,-4) -- (-0.5,4);
\end{tikzpicture}
\captionof{figure}{A region of A($\alpha$)-stability}\label{fig:Aalpha}
\end{center}

An A-stable
method is A($\pi/2$)-stable. A method is A($0$)-stable if the negative real axis itself lies in the
stability region.
Some numerical methods with A($\alpha$)-stability property which are appropriate for
stiff ODE problems have been developed. For example we can mention the {\bf implicit RK} methods and the {\bf backward differentiation formulas (BDF)} which will be introduced in sections \ref{sect:implicit_rk} and \ref{sect:bdf}. 
To see how to solve stiff ODEs using MATLAB, refer to the section \ref{sect:matlab} below.
\vsp

\subsection{L-stability}
As we pointed out, both implicit Euler and trapezoidal methods are A-stable but
there is a major difference between the stability regions of these methods.
The trapezoidal method is stable only in the left half-plane, whereas implicit Euler's method is
stable not only in the left-half plane but also over much of the right half-plane. See Figures \ref{fig:eulerregabs} and \ref{fig:regionabstrap}.
Let us solve the stiff ODE \eqref{eq:stiff1} on interval $[0,10]$ for $\lambda=-10^4$ with both methods. First, we assume that $y(0)=1$ which corresponds to the smooth exact solution $y(t)=\cos t$. Let $h=0.2$. Both methods provide satisfactory results
with norm-infinity error $9.998\times 10^{-6}$ for the Euler's method and $3.346\times 10^{-7}$ for the trapezoidal method.
Remember that the explicit methods failed to numerically solve this problem even at much smaller step sizes.
It is expectable that the trapezoidal method is more accurate because its convergence order is $2$ compared with the convergence order of the implicit Euler's method which is only $1$.  However, this is not the whole story. Let us change the initial value to
$y(0)=1.5$ which corresponds to exact solution
$$
y(t)= \frac{1}{2}e^{-10000t} + \cos t
$$
which includes a fast transient term. The initial value $y(0)=1.5$ rapidly (at a time scale of about $10^{-4}$) decreases toward the particular solution $\cos t$. The approximate solutions with $h=0.2$ are plotted in Figure \ref{fig:stiff5}.
\begin{center}
\includegraphics[scale=0.7]{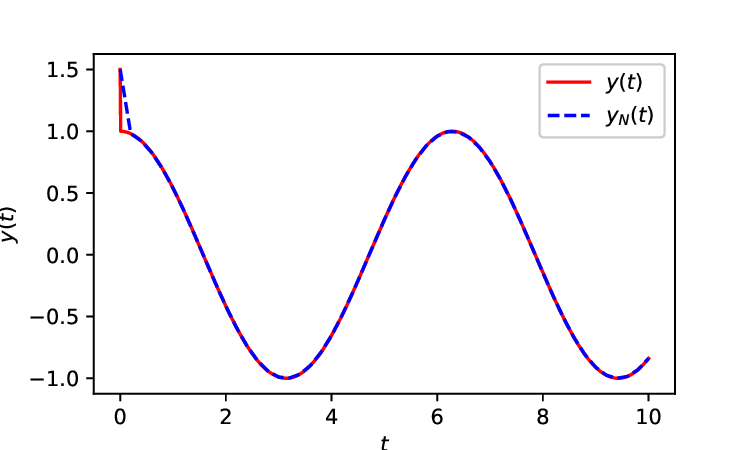}\includegraphics[scale=0.7]{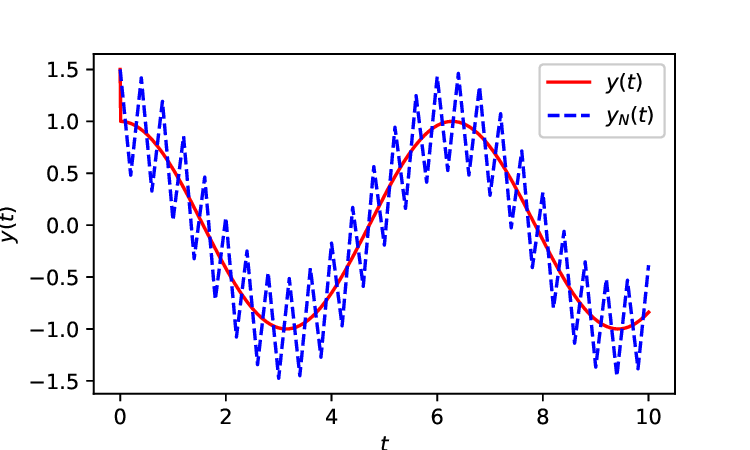}
\captionof{figure}{Numerical solution of stiff problem \eqref{eq:stiff1} with $\lambda=-10^4$ and $y(0)=1.5$ with step size $h=0.2$, the implicit Euler's method (left) and the trapezoidal method (right).}\label{fig:stiff5}
\end{center}
Both methods are still absolutely stable, but the
result of trapezoidal method shows unsatisfactory oscillations.
For absolute stability we test on equation $y'=\lambda y$ and obtain $y_{k+1}=R(z)y_k$ such that
\begin{equation*}
  R(z) = \frac{1}{1-z}, \quad |R(z)|\to 0 \;\mbox{as}\; z\to\infty,
\end{equation*}
for the implicit Euler's method and
\begin{equation*}
  R(z) = \frac{1+\frac{1}{2}z}{1-\frac{1}{2}z}, \quad |R(z)|\to 1 \;\mbox{as}\; z\to\infty,
\end{equation*}
for the trapezoidal method.
For problems with rapid transients, we aim for a method that can effectively damp in a single time step, because we intend to use a steplength much larger than the true decay time of the transient. To illustrate this point, refer to Figure \ref{fig:stiff6}, which provides a close-up comparison of the exact and numerical solutions obtained with each method over the time interval $[0,2]$.
\begin{center}
\includegraphics[scale=0.8]{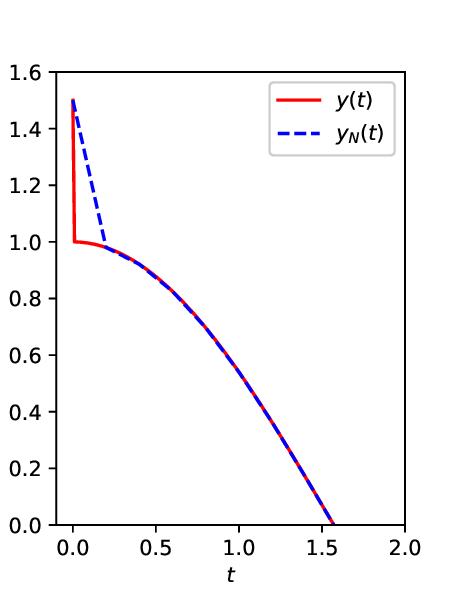}\includegraphics[scale=0.8]{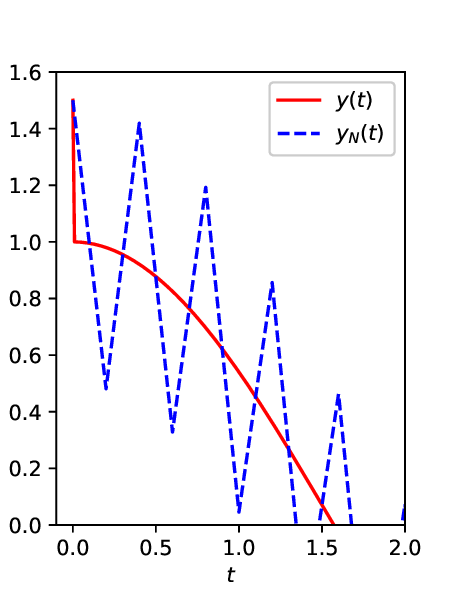}
\captionof{figure}{Closeup of solutions: the implicit Euler's method (left) and the trapezoidal method (right).}\label{fig:stiff6}
\end{center}
At the first time step, the implicit Euler's method damps very effectively toward the steady state solution $\cos t$ and continues
to produce very accurate results thereafter. In fact, this method, if is applied on stiff ODE \eqref{eq:stiff1}, yields
$$
y_{k+1} = \frac{1}{1+\lambda h} y_k-\frac{\lambda h}{1+\lambda h}\cos t_{k+1}-\frac{h}{1-\lambda h}\sin t_{k+1},
$$
which means
$$
y_{k+1}\approx \cos{t_{k+1}}\vsp
$$
because $\lambda h = -10^4\times 0.2=-2000$ and
\begin{equation*}
\frac{1}{1+\lambda h} = \frac{1}{2001}\doteq 0.0005\approx 0.
\end{equation*}
This implies that by the second time step we approximately fall on the steady state solution $\cos t$.
The trapezoidal method is also stable and the
results stay bounded, however, we have
\begin{equation*}
\frac{1+\frac{1}{2}\lambda h}{1-\frac{1}{2}\lambda h} = -\frac{999}{1001}\doteq -0.9980\approx -1.
\end{equation*}
By applying on stiff ODE \eqref{eq:stiff1}, this method gives
$$
y_{k+1} = \frac{1+\frac{1}{2}\lambda h}{1-\frac{1}{2}\lambda h}y_k - \frac{\frac{1}{2}\lambda h}{1-\frac{1}{2}\lambda h}[\cos t_{k+1}+\cos t_k]-\frac{\frac{1}{2} h}{1-\frac{1}{2}\lambda h}[\sin t_{k+1}+\sin t_k],
$$
or,
$$
y_{k+1}\approx -y_k + [\cos t_{k+1}+\cos t_k]. \vsp
$$
At the first time step, we have $y_1\approx  -1.5 + [\cos 0+\cos 0.2]\approx 0.48$ which falls below the steady state solution. Moving to the second time step, $y_2$ is $0.48 + [\cos 0.2+\cos 0.4] \approx 2.38$, overshooting the steady-state solution.
This pattern persists in subsequent time steps, resulting in the zigzag-shaped solution observed on the left-hand side of Figure \ref{fig:stiff5}.

 The results of the trapezoidal method can be improved if a small enough steplength $h$ is used. In Figure \ref{fig:stiff7}, numerical solutions  with $h=10^{-2}$ and $h=10^{-4}$ are shown.
 \begin{center}
\includegraphics[scale=0.7]{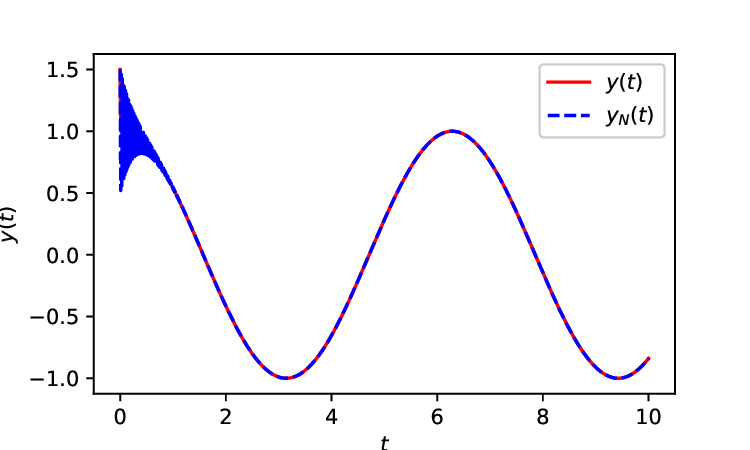}\includegraphics[scale=0.7]{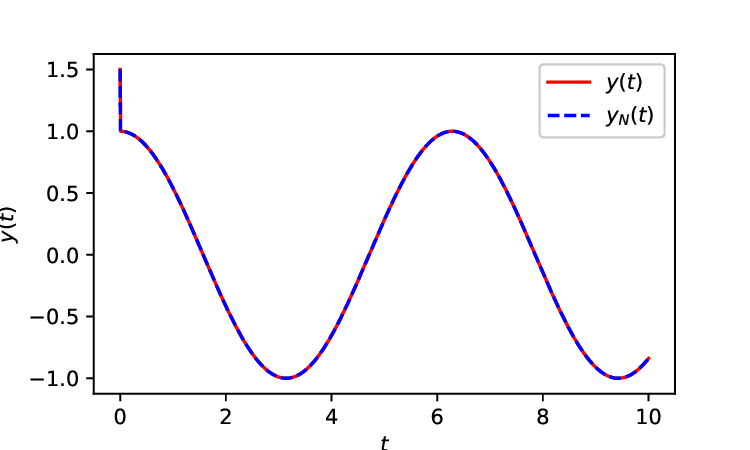}
\captionof{figure}{Numerical solution of stiff problem \eqref{eq:stiff1} with $\lambda=-10^4$ and $y(0)=1.5$ using the trapezoidal method, with $h=10^{-2}$ (left) and $h=10^{-4}$ (right).}\label{fig:stiff7}
\end{center}
With step length  $h=10^{-2}$ we still have oscillations near the initial time, while with step length  $h=10^{-4}$ a perfect numerical solution is observed. We note that, $h=10^{-4}$ is identical with the time scale of the transient term. This means that the trapezoidal method works perfect provided that the step size is chosen equal or smaller than the time scale of the transient terms in the solution.

However, our primary aim was to develop efficient numerical methods that be able to integrate the stiff problems with a rather large step length. What which makes the implicit Euler's method different from the trapezoidal method is the following property that the implicit Euler's method possesses while the trapezoidal method does not.
\begin{definition}
A one-step method is called {\bf L-stable} if it is A-stable and
$\displaystyle \lim_{z\to\infty}|R(z)| = 0$.
\end{definition}
The implicit Euler method is L-stable.
We will introduce some higher order L-stable methods in the forthcoming sections.

\section{Adaptive time stepping}

All solvers presented up to here use a single stepsize $h$ in all iterations. 
It is desirable to have an algorithm that can adjust the stepsize $h_k$ in each step $k$ to ensure that the local
truncation errors at all steps remain below a certain tolerance:
\begin{equation}\label{adapt:localerr}
|\tau_k| \leqslant \varepsilon,\; \mbox{for all }\; k.
\end{equation}
In certain parts of the domain, where the derivatives of $y$ have small amplitudes, larger stepsizes may suffice, resulting in a reduced computational expense due to a lower number of function evaluations. Since $y$ and its derivatives are unknown the truncation error can not be calculated analytically. 
To do so, we need an {\bf estimate} for $\tau_k$ at each step, and a way to select a new stepsize that will ensure that the estimated error is acceptably small.

\begin{remark}
It is more reasonable to keep the {\bf global error} under control rather than the {\bf local errors} \eqref{adapt:localerr}. Because bounding the $\tau_k$'s individually does not lead to a natural bound on the global error, since
it ignores the propagation of the error in each step. But local control is simple and is good enough in
practice if we take local tolerances smaller than what it seems necessary.
See the following.   
\end{remark}

Now we address the question of how to estimate the local errors and compute the new time step. 
\subsection{Using two methods of different order}
A good strategy is to use two methods of different orders. Informally, we use the more accurate one as an estimate for the exact solution, to estimate the error for the less accurate method.
Assume that for a given stepsize $h$ we have two methods:
\begin{itemize}
\item Method A: with local order $p$ 
and truncation error $\tau_{k+1}(h)$ to compute $y_{k+1}$,
\item Method B: with local order $p+1$ to compute a more accurate value $\tilde y_{k+1}$. 
\end{itemize}
Then the estimated local truncation error is 
$$
e_{k+1}: = \tilde y_{k+1} - y_{k+1}.
$$
If we suppose that the value at time step $k$ is exact, i.e. $y_k = y(t_k)$,
then by definition of local truncation errors we have 
\begin{equation}\label{adapt:locloc}
\begin{split}
y(t_{k+1}) &= y_{k+1} + \tau_{k+1}(h)\\
y(t_{k+1}) &= \tilde y_{k+1} + \mathcal O(h^{p+1}).
\end{split} 
\end{equation}
Since the local error of Methpd A is $\mathcal O(h^{p})$, we can expand the truncation error $\tau_{k+1}(h)$ as
$$
\tau_{k+1}(h) = C h^{p} + \mathcal O(h^{p+1}),
$$
for a constant $C$.
By subtracting two equations in \eqref{adapt:locloc} and taking absolute value we obtain
$$
|e_{k+1}|= |\tilde y_{k+1}-y_{k+1}| = |C|h^p + \mathcal O(h^{p+1}).
$$
By dropping the term $\mathcal O(h^{p+1})$, we approximately have
\begin{equation}\label{adapt:estimtau}
|\tau_{k+1}(h)|\approx |C|h^p \approx |e_{k+1}|.
\end{equation}
Consequently, if the error is small, i.e., 
\begin{equation}\label{adapt:errestimate}
|e_{k+1}|\leqslant \varepsilon
\end{equation}
then the step is {\em accepted} and the more accurate solution $\tilde y_{k+1}$
is assigned as an approximation for $y(t_{k+1})$. The algorithm then goes to the next time step $t_{k+1}+h$. Otherwise the tentative value $\tilde y_{k+1}$
is {\em rejected} and the step must be redone by a smaller steplength $h_{new}$. To estimate $h_{new}$ we use this fact that it must satisfy
$$
|\tau_{k+1}(h_{new})|\approx |C|h_{new}^p \leqslant \varepsilon.
$$
Taking the ratio of this and the estimate \eqref{adapt:estimtau} we obtain
$$
\left(\frac{h_{new}}{h}\right)^{p} \lesssim \frac{\varepsilon}{|e_{k+1}|}
$$
which gives an estimate for $h_{new}$ as 
$$
h_{new} \lesssim h \left( \frac{\varepsilon}{|e_{k+1}|} \right)^{1/p}.
$$
In practice, because this is just an estimate, one puts an extra coefficient in to be safe, typically something like
\begin{equation}\label{adapt:hnewformula}
h_{new} := 0.8 h \left( \frac{\varepsilon}{|e_{k+1}|} \right)^{1/p}.\vsp
\end{equation}
This formula decreases the stepsize if the error is large. 
This process should be repeated by replacing $h$ by $h_{new}$ until \eqref{adapt:errestimate} is satisfied and the step is accepted. 
Sometimes, other controls are added; like not decreasing or increasing $h$ by too much per time step. 
\vsp

\subsection{Embedded RK methods}
A disadvantage of the above adaptive scheme is that each step involving error estimation incurs a cost that is roughly twice that of a fixed step method because two methods are run. Nonetheless, by carefully selecting the methods, we can generate some computational overlap and thereby reduce the workload. 

An approach is to use an {\bf embedded pair} of RK formulas where most of the 
$f$ values are the same for both methods A and B.
For example, suppose we wish to create a pair for estimating the error in Euler's method (RK1) as Method A.
We also need a method with order $2$ (local order $3$) as Method B. Recall the RK2 method \eqref{rk2-2} which can be written as 
\begin{align*}
f_1& = f(t_k,y_k)\\
f_2& = f\Big(t_k+\frac{h}{2}, y_k + \frac{h}{2}f_1\Big)\\
\tilde y_{k+1} &= y_k + hf_2
\end{align*}
Then the value $y_{k+1}$ is computed via the Euler's method by
$$
y_{k+1} = y_k + hf_1,
$$ 
which requires essentially no extra function evaluation; the value of $f_1$ is used for both. Then following the rule \eqref{adapt:hnewformula} the adapted stepsize is computed as 
$$
h_{new} = 0.8 h \left( \frac{\varepsilon}{|e_{k+1}|} \right)^{1/2}
$$ 
where $e_{k+1}=\tilde y_{k+1}-y_{k+1} = h(f_2-f_1)$.
That value $y_{k+1}$ is not actually needed to compute because 
the error estimation $e_{k+1}$ is obtained in terms of $f_k$ values, and
the more accurate value $\tilde y_{k+1}$ is selected as an approximation for $y(t_{k+1})$. 
The initial stepsize $h = \varepsilon^{1/3}$ can be used because the local error in the first step is supposed to be of order $h^3$. 
The MATLAB code for such scheme is given here.
\begin{lstlisting}[style = matlab]
function [T,Y] = ODE12(f,y0,tspan,tol)
% Adaptive embedded method with RK1 and a RK2 formulas 
% Inputs: 
%  f: right hand side function f(t,y)
%  y0: initial condition 
%  tspan: [t0, tfinal] 
%  tol: predefined error 
% Output:
%  T: vector of time steps
%  Y: solution 
t = tspan(1); Y = y0; T = t;
h = tol^(1/3);
while t <= tspan(2)
   f1 = f(t,y0);
   f2 = f(t+h/2,y0+h/2*f1);
   e = norm(h*(f2-f1),inf);   
   if e < tol
      yt = y0 + h*f2;
      y0 = yt; t = t + h; 
      Y = [Y yt]; T = [T t];
      h = tol^(1/3);
   else
      h = 0.8*h*(tol/e)^(1/2);
   end
end
\end{lstlisting}

\begin{example}
As an example, consider the IVP 
\begin{align*}
&y' = \frac{1}{y^2+0.01}, \quad 0\leqslant t\leqslant 3, \\
& y(0)= 0
\end{align*}
 We solve this equation using the above MATLAB code with tolerance $\varepsilon = 10^{-3}$. 

\begin{lstlisting}[style = matlab1]
[t,y] = ODE12(@(t,y) 1/(y^2+0.01),0,[0 3],10e-3); 
plot(t,y(1,:),'-ob','MarkerFaceColor','b')
\end{lstlisting}

\begin{center}
\includegraphics[scale=.7]{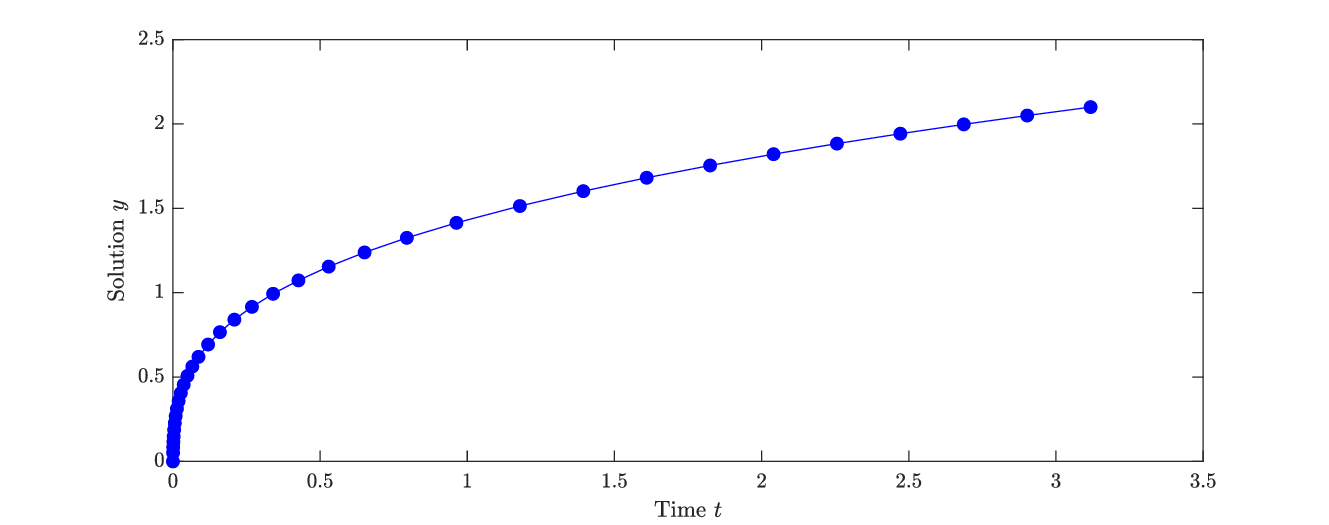}
\captionof{figure}{Numerical solution using the adaptive scheme}\label{fig:adapt2}
\end{center}


Plot of the solution is given in Figure \ref{fig:adapt2}. The scheme adapts much smaller stepsizes at the beginning of the time interval (where the solution takes more action) but remarkably larger stepsizes at remaining parts. This leads to saving memory and time specially when the underlying problem is a large system of equations with a long time interval.   
\end{example}
\vsp 

Embedded methods of higher order can be constructed by the right choice of coefficients.
One popular embedded pair is the {\bf Runge-Kutta-Fehlberg} method, which uses a
fourth-order and fifth-order formula that share most of the $f_k$ values. A formula of this form with stepsize selected by \eqref{adapt:hnewformula}
is the strategy employed, for instance, by MATLAB's \texttt{ode45}.

\section{Implicit Runge-Kutta methods}\label{sect:implicit_rk}
While high convergence order and ease of implementation are advantages of explicit RK methods and make them popular for solving various types of ODEs, their bounded stability regions render them impractical for a variety of important and challenging problems, such as stiff differential equations. Therefore, it is natural to develop, among other techniques, a class of implicit RK methods.

An $s$-stage {\bf implicit Runge-Kutta} method has the form
\begin{align*}
  z_\ell & =y_k + h\sum_{j=1}^{s}a_{\ell j} f(t_k+c_jh,z_j), \quad \ell=1,2,\ldots, s, \\
  y_{k+1} & = y_k + h\sum_{j=1}^{s}b_{j} f(t_k+c_jh,z_j).
\end{align*}
The Butcher's tableau for this formula has the form
\begin{equation*}
  \begin{array}{r|cccc}
    c_1 & a_{1,1} & a_{1,2} &\cdots  & a_{1,s}   \\
    c_2 & a_{2,1} & a_{2,2} &\cdots  & a_{2,s}   \\
    \vdots &\vdots & \vdots &  &  \vdots  \\
    c_s & a_{s,1} & a_{s,2} &\cdots  & a_{s,s}  \\ \hline
     & b_1 & b_2 & \cdots   & b_s
  \end{array}
\end{equation*}
Here, the diagonal and upper diagonal parts of the tableau may have nonzero values.
To advance form time level $t_k$ to $t_{k+1}$ using an $s$-stage implicit RK method we should solve a system of $s$ nonlinear equations
for $s$ unknowns $z_1,z_2,\ldots,z_s$. For a system of equations  $\y'=f(t,\y)$ with $m$ differential equations we must solve a system of $sm$ equations in $sm$ unknowns at each time step.

There typically are some ways to derive coefficients $\{b_j,c_\ell,a_{\ell, j}\}$ for a given accuracy, provided the number of stages is sufficiently
large. The dominant approach converts the IVP in to integral equation
$$
y(t) = y(t_k) + \int_{t_k}^{t} f(\tau,y(\tau))d\tau, \quad t\in[t_k,t_{k+1}],\vsp
$$
uses a polynomial interpolation of order $s$ for $f(\tau,y(\tau))$ on a predefined set of interpolation points $\tau_1,\ldots,\tau_s\in[t_k,t_{k+1}]$, and collocate the resulting equation at $t=\tau_1,\ldots,\tau_s$
to predict $z_1,z_2,\ldots,z_s$ and $y_{k+1}$. This approach is called {\em collocation}.
We omit derivation details but we present the Butcher's tableaus for some formulas instead.
An extensive theory has been developed by Butcher for analyzing these
methods.

The tableau for a two-stage formula with good convergence and stability properties is

\begin{equation*}
  \begin{array}{r|cc}
    (3-\sqrt 3)/6 & 1/4 & (3-2\sqrt 3)/12   \\
    (3+\sqrt 3)/6 & (3+2\sqrt 3)/12 & 1/4  \\ \hline
      & 1/2 & 1/2
  \end{array}
\end{equation*}

This method is also called {\em two-stage Gauss method} because the transferred Gauss-Legendre points $\{-\tfrac{\sqrt 3}{3},\tfrac{\sqrt 3}{3}\}$ are used for the polynomial interpolation. It can be shown that the local truncation error for this method has size $\mathcal O(h^5)$, and the global error behaves like $\mathcal O(h^4)$.

\begin{workout}\label{wo3-guassstab}
Show that the absolute stability region of the two-stage Gauss method is
$$
S= \left\{z\in \C: \left|\frac{1+\tfrac{1}{2}z+\tfrac{1}{12}z^2}{1-\tfrac{1}{2}z+\tfrac{1}{12}z^2}\right| \leqslant 1   \right\}.
$$
Show that the negative part of the real line falls in $S$. More generally, show that the left-half plane falls in $S$.
\end{workout}

The nice stability feature of the two-stage Gauss method, as outlined in Workout \ref{wo3-guassstab}, comes at the cost of solving a nonlinear system of algebraic equations for each time step.
In general, a fully implicit RK method where each $z_\ell$ value depends on all the $z_\ell$ values, can be costly to implement for a system of ODEs. This is because a nonlinear system of $sm$ equations in $sm$ unknowns must be solved at each time step.
However, a subclass of implicit methods that are simpler to implement are the {\bf diagonally implicit Runge–Kutta} methods (DIRK methods), in which 
$a_{\ell, j}=0$ for $j>\ell$, i.e.,
$z_\ell$ depends on $z_j$ for $j=1,\ldots,\ell$.
An 
$s$-stage DIRK method, when applied to a system of 
$m$ nonlinear ODEs, requires solving a sequence of $s$ nonlinear systems, each of size $m$, rather than a coupled set of $sm$ equations.

As an example of a $3$-stage DIRK method, we can mention a method with following Butcher tableau:

\begin{equation*}
  \begin{array}{r|ccc}
     0&0  &  &  \\
     1/2& 1/4& 1/4  &\\
    1 & 1/3 & 1/3  & 1/3\\ \hline
     &1/3 &1/3 &1/3
  \end{array}
\end{equation*}
\vsp 

This method is of second order accuracy, and also is known as  the TR-BDF2 method.

%% file: lec3_part3.tex
\vsp 

\section{Multistep methods}\label{sect:adams}
All methods we considered so far are {\em single-step} or {\em one-step} methods, since at a typical step $y_{k+1}$ is determined solely from $y_k$. In this section we study the {\em multistep methods} where the solution $y_{k+1}$ depends on more than one previous values, i.e., $y_k,y_{k-1},\ldots$. Three families of such methods are widely used; {\bf Adams-Bashforth (AB)} and {\bf Adams-Moulton (AM)} methods and {\bf backward differentiation formulas (BDF)}.

Again, we reformulate the solution of IVP $y'(t)=f(t,y)$ at $t=t_{k+1}$ as
\begin{equation}\label{multi-eqintegral}
  y(t_{k+1}) = y(t_k) + \int_{t_k}^{t_{k+1}}f(t,y(t)) dt.
\end{equation}
A numerical scheme for computing $y(t_{k+1})$ can be obtained if the integral on the right-hand side of \eqref{multi-eqintegral} is approximated by a numerical quadrature. A numerical quadratures for computing a definite integral of the form
$$
\int_{t_{k}}^{t_{k+1}}g(t)dt \vsp
$$
can be developed by replacing $g$ with an interpolation polynomial $p$ of a certain degree.
In a $q$-step AB method we assume that $p$
interpolates $g$ at points $\{t_{k-q},\ldots, t_{k-1},t_k\}$ while in a $q$-step AM method at points $\{t_{k-q+1},\ldots,t_k,t_{k+1}\}$.
In AM methods the solution $y_{k+1}$ depends on already computed values $y_k,y_{k-1},\ldots, y_{k-q}$ thus AB methods are explicit.
In AM methods, on the other hand, $y_{k+1}$ depends on $q-1$ previous values and $y_{k+1}$ itself, thus AM methods are implicit.

\subsection{Adams-Bashforth methods}
Let $q=1$. The linear interpolant of $g$ at points $\{t_{k-1},t_{k}\}$ then is
\begin{equation*}
  p_1(t) = \frac{1}{h}[(t_k-t)g(t_{k-1})+(t-t_{k-1})g(t_k)],
\end{equation*}
with error function
\begin{equation*}
 e(t) = g(t)- p_1(t) = \frac{1}{2}(t_k-t)(t-t_{k-1})g''(\zeta_k),
 \vsp
\end{equation*}
for some $t_{k-1}\leqslant \zeta_k\leqslant t_{k}$. Integrating over $[t_k,t_{k+1}]$ yields
$$
\int_{t_{k}}^{t_{k+1}}g(t)dt = \int_{t_{k}}^{t_{k+1}}p_1(t) + \int_{t_{k}}^{t_{k+1}}e(t) = \frac{h}{2}[3g(t_k)-g(t_{k-1})] + \frac{5h^3}{12}g''(\xi_k)
$$
for some $t_{k-1}\leqslant \xi_k\leqslant t_{k+1}$. Applying this to \eqref{multi-eqintegral}, gives
\begin{equation}\label{multi-ab1ex}
  y(t_{k+1})=y(t_k) + \frac{h}{2}\big[3f(t_k,y(t_k))-f(t_{k-1},y(t_{k-1}))\big] + \frac{5h^3}{12}y'''(\xi_k).
\end{equation}
Dropping the error term and replacing the exact values $y(t_k)$ by the approximate values $y_k$, we obtain the $2$-step AB formula
\begin{shaded}
\vspace*{-0.1cm}
\begin{equation}\label{multi-ab2}
  y_{k+1}=y_k + \frac{h}{2}\big[3f(t_k,y_k)-f(t_{k-1},y_{k-1})\big], \quad k=1,2,\ldots.
\end{equation}
\vspace*{-0.3cm}
\end{shaded} 
At the initial level $k=1$, computing $y_2$ requires both $y_0$ and $y_1$, yet we only have $y_0$ from the initial value. The value of $y_1$ must be computed using another method. If $y_1$ is obtained in such a way that $|y_1-y(t_1)|\leqslant c_1h^2$ then it can be proved that the global error of the method \eqref{multi-ab2} satisfies
$$
|e_N|\leqslant Ch^2
\vsp
$$
provided that $h$ is sufficiently small, $f(t,y)$ is Lipschitz continuous and $y$ is $3$ times continuously differentiable. 
See section \ref{sect:error_anal_multi}. 
To calculate $y_1$ we have at least two possibilities from previous sections. The explicit Euler's method gives
$$
y_1 = y_0 + hf(t_0,y_0)
$$
with error $|y(t_1)-y_1|\leqslant c_1h^2$,
and the RK2 method
$$
y_1 = y_0 + \frac{h}{2}[f(t_0,y_0) + f(t_{0}+h,y_0+hf(t_0,y_0))]
$$
with error $|y(t_1)-y_1|\leqslant c_1h^3$, which is more than adequate.

The $3$-step AB method is obtained by interpolating $g$ at points $\{t_{k-2},t_{k-1},t_k\}$ and then integrating over $[t_k,t_{k-1}]$ as is required in \eqref{multi-eqintegral}. The interpolant will be
$$
p_2(t)= \ell_0(t)g(t_{k})+\ell_1(t)g(t_{k-1})+\ell_2(t)g(t_{k-2})
$$
with Lagrange functions
\begin{align*}
  \ell_0(t) & = \frac{1}{2h^2}(t-t_{k-1})(t-t_{k-2}), \\
  \ell_1(t) & = \frac{1}{h^2}(t-t_{k})(t-t_{k-2}), \\
  \ell_2(t) & = \frac{1}{2h^2}(t-t_{k})(t-t_{k-1}),
\end{align*}
and interpolation error
$$
e(t)=g(t)-p_2(t) = \frac{1}{6}(t-t_{k})(t-t_{k-1})(t-t_{k-2})g'''(\zeta_k)
$$
for some $t_{k-2}\leqslant \zeta_k\leqslant t_{k}$. Exact integration of the polynomial $p_2$ and error term $e$ reveals that
$$
\int_{t_k}^{t_{k+1}}g(t)dt = \frac{h}{12}\left[23g(t_k)-16g(t_{k-1})+5g(t_{k-1}) \right]+\frac{3h^4}{8}g'''(\xi_k)
\vsp
$$
for some $t_{k-2}\leqslant \xi_k\leqslant t_{k+1}$. Using this quadrature for \eqref{multi-eqintegral} and dropping the error term, we obtain the $3$-step AB method
\begin{shaded}
\vspace*{-0.1cm}
\begin{equation}\label{AB3}
  y_{k+1}=y_k + \frac{h}{12}\left[23y'_k -16y'_{k-1}+5y'_{k-2} \right], \quad k=2,3,4,\ldots,
\end{equation}
\vspace*{-0.2cm}
\end{shaded} 
where $y'_k = f(t_k,y_k)$. In this formula $y_1$ and $y_2$ values must be obtained separately by other methods with errors at most of order $h^3$ (such as RK2 method) to keep the global error of \eqref{AB3} of order $h^3$.

Higher order AB methods can be derived similarly. In Table \ref{tb:ABmethods} the AB methods of order $1$ through $4$ are listed. Local truncations errors are given in the last column. See section \ref{sect:error_anal_multi} for an error analysis.

\begin{table}[!th]
  \centering
    \caption{Adams-Bashfoth methods. Here by $y'_k$ we mean $f(t_k,y_k)$.}\label{tb:ABmethods}
  \begin{tabular}{lclc}
    \hline
    $q$ & Order & $\qquad\qquad \qquad$ Method & $\tau_{k+1}$ \\ \hline \\
    $0$ & $1$ & $y_{k+1}= y_k+hy'_k$& $\tfrac{1}{2}h^2y''(\xi_k)$ \\\\
    $1$ & $2$ & $y_{k+1}= y_k+\tfrac{1}{2}h[3y'_k-y'_{k-1}]$& $\tfrac{5}{12}h^3y'''(\xi_k)$ \\\\
    $2$ & $3$ & $y_{k+1}= y_k+\tfrac{1}{12}h[23y'_k-16y'_{k-1}+5y'_{k-2}]$&$ \tfrac{3}{8}h^4y^{(4)}(\xi_k)$ \\\\
    $3$ & $4$ & $y_{k+1}= y_k+\tfrac{1}{24}h[55y'_k-59y'_{k-1}+37y'_{k-2}-9y'_{k-3}]$& $\tfrac{251}{720}h^5y^{(5)}(\xi_k)$ \\\\
    \hline
  \end{tabular}
\end{table}

Compared to Runge-Kutta methods with the same order of accuracy, multistep methods require fewer evaluations of $f$ at each time step. For instance, in the explicit RK4 method \eqref{rk4}, we need $4$ function evaluations per time step, whereas in the explicit $4$-step AB method, only $1$ function evaluation is needed in each time step, provided that previous function values of $f$ are reused.

\subsection{Adams-Moulton methods}
The implicit Euler's method can be obtained by setting $q=0$, i.e., by using  constant interpolation $p_0=t_{k+1}$ in the AM process. Besides,
with $q=1$ the AM method is identical with trapezoidal method \eqref{trap-method} because the resulting quadrature in $[t_k,t_{k+1}]$ is just the well-known trapezoidal rule which is obtained by linear interpolation on points $\{t_k,t_{k+1}\}$.
For $q=2$ the interpolant $p_2$ should set up on points
$\{t_{k-1},t_{k},t_{k+1}\}$ and integration should apply on interval $[t_k,t_{k+1}]$. The resulting method is the $2$-step AM method
listed in the second row of Table \ref{tb:AMmethods}. Other AM methods up to order $4$ are given in this table with the corresponding truncation errors in the last column. As we observe, the $h$ powers in local truncation errors are the same as their counterparts in the table of AM methods. However, the constants behind the $h$ are remarkably smaller in the AM methods. The error analysis is given in section \ref{sect:error_anal_multi} below. 

\begin{table}[!th]
  \centering
    \caption{Adams-Moulton methods. Here by $y'_k$ we mean $f(t_k,y_k)$.}\label{tb:AMmethods}
  \begin{tabular}{lclc}
    \hline
    $q$ & Order & $\qquad\qquad \qquad$ Method & $\tau_{k+1}$ \\ \hline \\
    $0$ & $1$ & $y_{k+1}= y_k+hy'_{k+1}$& $-\tfrac{1}{2}h^2y''(\xi_k)$ \\\\
    $1$ & $2$ & $y_{k+1}= y_k+\tfrac{1}{2}h[3y'_{k+1}-y'_{k}]$& $-\tfrac{1}{12}h^3y'''(\xi_k)$ \\\\
    $2$ & $3$ & $y_{k+1}= y_k+\tfrac{1}{12}h[5y'_{k+1}+8y'_{k}-y'_{k-1}]$&$ -\tfrac{1}{24}h^4y^{(4)}(\xi_k)$ \\\\
    $3$ & $4$ & $y_{k+1}= y_k+\tfrac{1}{24}h[9y'_{k+1}+19y'_{k}-5y'_{k-1}+y'_{k-2}]$& $-\tfrac{19}{720}h^5y^{(5)}(\xi_k)$ \\\\
    \hline
  \end{tabular}
\end{table}

Discussion on using other methods with consistent accuracies for calculating few initial values to bootstrap the AB methods is applicable here for the AM methods as well.

Since AM methods are implicit, an initial guess $y_{k+1}^{(0)}$ and an iteration on $y_{k+1}$ are needed in each time step.
A choice for initial guess can be a solution of an AB method of the same order. For example, to implement the AM method of order $3$ we may write
\begin{align*}
  y_{k+1}^{(0)} & =y_k+\tfrac{1}{12}h[23y'_k-16y'_{k-1}+5y'_{k-2}], \\
  y_{k+1}^{(1)} & =y_k+\tfrac{1}{12}h[5f(t_{k+1},y_{k+1}^{(0)})+8y'_{k}-y'_{k-1}].
\end{align*}
If $h$ is sufficiently small, we can accept $y_{k+1}^{(1)}$ as the solution $y_{k+1}$. Otherwise, we can proceed with more iterations at the expense of additional evaluations of $f(t,y)$.
\vsp

\subsection{Backward differentiation formulas (BDF)}\label{sect:bdf}
Another class of efficient numerical methods with excellent stability properties is the backward differentiation formulas (BDF). As the name implies, they are backward (implicit) formulas.
For constructing a BDF of order $q$, we assume that $p$ is a polynomial of degree $q$ that interpolates $y(t)$ at points
$\{ t_{k+1}, t_k,\ldots, t_{k-q+1}\}$ for $k\geqslant q-1$:
\begin{equation}\label{bdf:poly}
  p(t) = \sum_{j=-1}^{q-1} \ell_j(t) y(t_{k-j})
\end{equation}
where $\ell_j(t), j=-1,\ldots,q-1$ are Lagrange polynomials on points $\{ t_{k+1}, t_k,\ldots, t_{k-q+1}\}$. Then we use
\begin{equation}\label{bdf:deriv}
  p'(t_{k+1})\approx y'(t_{k+1}) = f(t_{k+1},y(t_{k+1})).
\end{equation}
Since the interpolation points are equidistance with distance $h$, we can use the change of variable
$$
\theta = \frac{t-t_{k+1}}{h}
\vsp
$$
to simplify the notation. This change of variable maps the interpolation points to integer set $\{0,-1,-2,\ldots,-(k-q)\}$.
In particular, $t=t_{k+1}$ is mapped to $\theta=0$.
If Lagrange functions on this scaled points are denoted by $\tilde{\ell}_j(\theta)$, we can simply show that
$$
\ell_j(t) = \tilde{\ell}_j(\theta) \quad \mbox{and}\quad \ell'_j(t)=\frac{1}{h}\tilde \ell_j(\theta).
$$
On the other hand, the Lagrange interpolation error is
\begin{align*}
e(t) =& \frac{(t-{t_{k+1}})(t-t_k)\cdots(t-t_{k-q+1})}{(q+1)!}y^{(q+1)}(\zeta_k(t))\\
 = &\frac{h^{q+1}\theta(\theta-1)\cdots(\theta-k+q)}{(q+1)!}y^{(q+1)}(\zeta_k(t))
\end{align*}
for some $t_{k-q+1}\leqslant \zeta_k(t)\leqslant t_{k+1}$. The error of differentiation in \eqref{bdf:deriv} at $t=t_{k+1}$ (or $\theta=0$) then is
\begin{equation}\label{bdf:errorderiv}
  e'(t_{k+1}) = \frac{1}{h}\frac{h^{q+1}(-1)^{k-q}(k-q)!}{(q+1)!}y^{(q+1)}(\xi_k) =:\tilde \tau_{k+1},
\end{equation}
where $\xi_k=\zeta_k(t_{k+1})$.
Combining \eqref{bdf:poly} and \eqref{bdf:deriv} and adding the error term give
\begin{equation*}
  \frac{1}{h}\tilde \ell'_{-1}(0) y(t_{k+1}) + \frac{1}{h}\sum_{j=0}^{q-1} \tilde \ell'_j(0) y(t_{k-j}) = f(t_{k+1},y(t_{k+1}))+\tilde \tau_{k+1}.
\end{equation*}
By setting
\begin{equation}\label{bdf:alphabeta}
\beta = \frac{1}{\tilde \ell'_{-1}(0)}, \quad \alpha_j = -\frac{\tilde \ell'_j(0)}{\tilde \ell'_{-1}(0)}, \quad j=0,\ldots,q-1,
\end{equation}
we obtain
\begin{equation*}
  y(t_{k+1}) = \sum_{j=0}^{q-1} \alpha_j y(t_{k-j}) +h\beta f(t_{k+1},y(t_{k+1}))+h\tilde \tau_{k+1}.
\end{equation*}
By replacing the exact values $y(t_k)$ by approximate values $y_k$ when dropping the truncation error
\begin{equation*}
  \tau_{k+1} = h\tilde{\tau}_{k+1} = \frac{(-1)^{k-q}(k-q)!}{(q+1)!}h^{q+1} y^{(q+1)}(\xi_k),
\end{equation*}
the $q$-step BDF method is obtained as
\begin{shaded}
\vspace*{-0.1cm}
\begin{equation}\label{bdf-method}
  y_{k+1} = \sum_{j=0}^{q-1} \alpha_j y_{k-j} +h\beta f(t_{k+1},y_{k+1}), \quad k=q-1,q,q+1,\ldots.
\end{equation}
\vspace*{-0.2cm}
\end{shaded} 
The coefficients $\beta$ and $\alpha_j$ can be simply computed (for example using a symbolic toolbox such as Maple) from \eqref{bdf:alphabeta}. In Table \ref{tb:BDFmethods0} the $q$-step BDF methods for $q=1,2,\ldots,5$ are given.



\begin{table}[!th]
  \centering
    \caption{BDF methods. Here by $y'_{k+1}$ we mean $f(t_{k+1},y_{k+1})$.}\label{tb:BDFmethods0}
  \begin{tabular}{llc}
    \hline
    $q$  & $\qquad\qquad \qquad\qquad\qquad$ Method                                        & $\tau_{k+1}$ \\ \hline \\
    $1$  & $y_{k+1}= y_k + hy'_{k+1}$                                & $-\tfrac{1}{2}h^2y''(\xi_k)$ \\\\
    $2$  & $y_{k+1}= \tfrac{4}{3}y_k-\tfrac{1}{3}y_{k-1} + \tfrac{2}{3}hy'_{k+1}$ & $-\tfrac{2}{9}h^3y'''(\xi_k)$ \\\\
    $3$  & $y_{k+1}= \tfrac{18}{11}y_k-\tfrac{9}{11}y_{k-1} +\tfrac{2}{11}y_{k-2}+ \tfrac{6}{11}hy'_{k+1}$ &$ -\tfrac{3}{22}h^4y^{(4)}(\xi_k)$ \\\\
    $4$  & $y_{k+1}= \tfrac{48}{25}y_k-\tfrac{35}{25}y_{k-1} +\tfrac{16}{25}y_{k-2}-\tfrac{3}{25}y_{k-3}+ \tfrac{12}{25}hy'_{k+1}$ & $-\tfrac{12}{625}h^5y^{(5)}(\xi_k)$ \\\\
    $5$  & $y_{k+1}= \tfrac{300}{137}y_k-\tfrac{300}{137}y_{k-1} +\tfrac{200}{137}y_{k-2}-\tfrac{75}{137}y_{k-3}+\tfrac{12}{137}y_{k-4}+ \tfrac{60}{137}hy'_{k+1}$ & $-\tfrac{10}{137}h^6y^{(6)}(\xi_k)$ \\\\
    \hline
  \end{tabular}
\end{table}

The case $q=1$ is simply the implicit Euler's method \eqref{imeuler:method}. All discussions concerning the initial values $y_1,\ldots,y_{q-1}$ and iteration for nonlinearity that we outlined for the Adams-Moulton methods are applicable here for BDF methods.

As the local truncation error for a $q$-step BDF method behaves as $\mathcal{O}(h^{q+1})$, it is predictable that the global error of such method behaves as $\mathcal{O}(h^q)$. A general error analysis is given in section \ref{sect:error_anal_multi}.
\vsp 
\begin{workout} \label{wo:BDF3forSys}
Consider solving the linear system of equations $y'(t)=Ay(t) + f(t)$ with initial condition $y(0)=y_0$, given constant matrix $A\in \R^{n\times n}$ and know vector function $f$. Write down the BDF3 scheme for solving this system. 
\end{workout}
\vsp 
\begin{labexercise}
Let's revisit Example \ref{ex:heat} (MOL solution of the diffusion equation). Impose the initial condition
$$
u_0(x) = \sin(\pi x), \quad 0\leqslant x\leqslant 1.\vsp 
$$
With a spatial step size of $\Delta x = 0.001$, apply the BDF3 method on the resulting system of ordinary differential equations (ODEs) and compute the PDE solution with steplengths $h=0.1,0.05,0.025,0.0125,$ and $0.00625$ up to the final time $t = 0.5$. Plot the errors, compute the convergence orders, and compare them with the theoretical order $q=3$.
To initiate BDF3 iterations, in addition to the initial condition $\y_0$, we require $\y_1$ and $\y_2$, which should be computed using a one-step method with a truncation error of at least $\mathcal{O}(h^3)$. For instance, the trapezoidal or RK2 methods can be employed for this purpose.
The exact solution of this equation with the given initial condition is $u(x,t) = e^{-\pi^2 t}\sin x$, which we can use to compute the errors and orders.

Explain the advantages of BDF3 over explicit Euler, implicit Euler, and RK methods. Recall that when the spatial step size $\Delta x$ decreases, the size of the ODE system increases and the system becomes stiff. 
\end{labexercise}
\vsp 

\subsection{Error analysis of multistep methods}\label{sect:error_anal_multi}

In general, a multistep method, including AB, AM, and BDF methods, can be formulated as
\begin{equation}\label{multi-general}
y_{k+1} = \sum_{j=0}^{q} a_jy_{k-j} + h \sum_{j=-1}^{q} b_j f(t_{k-j},y_{k-j}), \quad k\geqslant q.
\vsp
\end{equation}
In both Adams methods (AB and AM methods) $a_0=1$ and $a_j=0$ for $j=1,\ldots,q$. In AB methods $b_{-1}=0$ and in AM methods $b_q=0$.
In BDF methods $b_{-1}=\beta$, $b_j=0$ for $j=0,\ldots,q$, and $a_q=0$.

We present the analysis for the general form \eqref{multi-general} with only a restrictive condition on coefficients $a_j$. The analysis will be valid for all three classes of methods mentioned above. It also works for some one-step schemes such as implicit Euler and trapezoidal methods as they are special cases of AM methods. The truncation error for formula \eqref{multi-general} is defined as
\begin{equation}\label{multi-truncationerr}
\tau_{k+1} = y(t_{k+1})-\sum_{j=0}^{q} a_jy(t_{k-j}) + h \sum_{j=-1}^{q} b_j y'(t_{k-j}), \quad k\geqslant q.
\end{equation}
For function $\tau(h)$ defined by
$$
\tau(h) = \frac{1}{h}\max_{q\leqslant k\leqslant N} |\tau_k|,
\vsp
$$
if we have $\tau(h)\to 0$ as $h\to 0$ then we say the method \eqref{multi-general} is {\em consistent}. If
$$
\tau(h) = \mathcal{O}(h^m)
$$
for some $m\geqslant 1$, we say the order of consistency is $m$. The following theorem gives necessary and sufficient conditions on coefficients in
\eqref{multi-general} to achieve the consistency.
\begin{theorem}[\cite{Atkinson-et-al:2009}]\label{thm:multi_consistency}
  Let $m\geqslant 1$ be a given integer. For $\tau(h)\to 0$ for all continuously differentiable function $y$, that is, for method \eqref{multi-general} to be consistent, it is necessary and sufficient that
  \begin{align}
    \sum_{j=0}^q a_j = & 1, \label{sum_a_j1} \\
    -\sum_{j=0}^q ja_j + \sum_{j=-1}^{q}b_j= & 1 \label{sum_b_j1}.
  \end{align}
  Furthermore, to have the consistency order $m$, i.e., $\tau(h)=\mathcal O(h^m)$ for all functions $y$ that are $m+1$ continuously differentiable, it is necessary and sufficient that \eqref{sum_a_j1} and \eqref{sum_b_j1} hold and that
  \begin{equation}\label{sum_aj_bj}
   \sum_{j=0}^q (-j)^ia_j + i\sum_{j=-1}^{q}(-j)^{i-1}b_j= 1, \quad i=2,3,\ldots,m.
  \end{equation}
\end{theorem}
\proof For the proof just write the degree $m$ Taylor expansion of $y(t)$ around $t_k$ and manipulate the truncation error \eqref{multi-truncationerr}. It is left as an exercise.
$\qed$
\vsp 

\begin{workout}
Write the proof of Theorem \ref{thm:multi_consistency}
\end{workout}
\vsp 

The following theorem proves the convergence of the multistep method \eqref{multi-general}. The theorem does not cover all the multistep methods but includes all methods we considered so far such as AB, AM and BDF schemes. 

\begin{theorem}[\cite{Atkinson-et-al:2009}]
Consider solving the initial value problem $y'(t) =  f(t,y(t))$ for $t\geqslant t_0$ with initial condition $y(t_0)=y_0$
using the multistep method \eqref{multi-general}. Assume that $f(t,y)$ is continuous and satisfies the Lipschitz condition with respect to variable $y$ with Lipschitz constant $L>0$. Assume that the initial errors satisfy 
\begin{equation*}
  \eta(h):= \max_{0\leqslant k\leqslant q} |y(t_k)-y_k|\to 0, \; \mbox{as} \; h\to0.
\end{equation*}
Moreover, assume that $y$ is continuously differentiable and the method is consistent, that is, $\tau(h)\to 0$ as $h\to0$. Finally, assume that all coefficients $a_j$ are nonnegative:
\begin{equation}\label{multi:aj_pos}
  a_j\geqslant 0,\quad j=0,1,\ldots, q.
\end{equation}
  Then the method \eqref{multi-general} is convergent and
  \begin{equation}\label{multi:conv1}
    \max_{0\leqslant k\leqslant N} |y(t_k)-y_k|\leqslant c_1 \eta(h) + c_2\tau(h)
  \end{equation}
  for some suitable constants $c_1$  and $c_2$. If the solution $y(t)$ is $m+1$ continuously differentiable, the method \eqref{multi-general} is of consistency order $m$, and the initial errors satisfy $\eta(h)=\mathcal O(h^m)$ then the convergence order is $m$, i.e., the error is of size $\mathcal{O}(h^m)$. 
\end{theorem}

\proof Let $e_k=y(t_k)-y_k$. By subtracting \eqref{multi-general} from \eqref{multi-truncationerr}, we obtain
\begin{equation*}
e_{k+1}=\sum_{j=0}^{q} a_je_{k-j} + h \sum_{j=-1}^{q} b_j \big[f(t_{k-j},y(t_{k-j}))-f(t_{k-j},y_{k-j})\big] + \tau_{k+1}.
\end{equation*}
Taking absolute value, using the Lipschitz condition, and using assumption \eqref{multi:aj_pos} yield
$$
|e_{k+1}|\leqslant \sum_{j=0}^{q} a_j|e_{k-j}| + hL \sum_{j=-1}^{q} |b_j|\, |e_{k-j}| + |\tau_{k+1}|.
$$ 
By introducing the notation
$$
E_k = \max_{0\leqslant i\leqslant k}|e_i|, \quad k=0,1,\ldots,N
\vsp 
$$
we can write 
$$
|e_{k+1}|\leqslant E_k\sum_{j=0}^{q} a_j + hL E_{k+1} \sum_{j=-1}^{q} |b_j| + h\tau(h).
\vsp 
$$
From \eqref{sum_a_j1} we know that the sum of $a_j$s is $1$, thus we have
$$
|e_{k+1}|\leqslant E_k + hc E_{k+1} + h\tau(h), \quad c= L \sum_{j=-1}^{q} |b_j|. 
$$
Since the right hand side is trivially a bound for $E_k$ and $E_{k+1}=\max\{|e_{k+1}|,E_k\}$, we simply have 
$$
E_{k+1}\leqslant E_k + hcE_{k+1}+ h\tau(h).
\vsp 
$$
For $hc\leqslant \frac{1}{2}$, which is true as $h\to 0$, we obtain 
$$
E_{k+1}\leqslant \frac{E_k}{1-hc} +\frac{h}{1-hc}\tau(h)\leqslant (1+2hc)E_k + 2h \tau(h). 
$$
If we proceed as in the analysis of the Euler's method in section \ref{sect:error_anal_euler} we finally have
\begin{equation*}
  E_{N}\leqslant e^{2c(b-t_0)}\eta(h) + \left[ \frac{e^{2c(b-t_0)}-1}{c} \right]\tau(h),
\end{equation*}
which complete the proof.
$\qed$

To obtain the rate of convergence of $\mathcal O(h^m)$ it is necessary that $\tau_{k+1}=\mathcal O(h^{m+1})$ in each step $k$ (i.e., $\tau(h)=\mathcal O(h^m)$ which needs relation \eqref{sum_aj_bj} to hold), and the initial values $y_0,y_1,\ldots,y_q$ need to be computed with an accuracy of only $\mathcal O(h^m)$, i.e., $\eta(h)=\mathcal O(h^m)$.
\vsp 

\subsection{Stability regions of multistep methods}
Again consider the general multistep method \eqref{multi-general}. 
If we apply it on test equation $y'(t)=\lambda y(t)$ we obtain
\begin{equation*}
  y_{k+1} = \sum_{j=0}^{q} a_jy_{k-j} + h\lambda \sum_{j=-1}^{q} b_j y_{k-j}.
  \vsp 
\end{equation*}
Letting $z=\lambda h$ and rearranging the above equation give
\begin{equation}\label{multi:abs-diff}
  (1-zb_{-1})y_{k+1} - \sum_{j=0}^{q} (a_j+zb_j)y_{k-j}=0.
\end{equation}
This is a homogenous linear {\em difference equation} of order $q+1$. For absolute stability, we need to find the general solution $y_k$ and impose some condition on $z$ in order to have a bounded $|y_k|$ as $k\to \infty$. The general theory for solving linear difference equation is given in the Appendix A.
The general solution can be found by looking for solutions of the special form
$$
y_k = r^k, \quad k\geqslant 0.
$$
Substituting this special solutions in to \eqref{multi:abs-diff} and cancelling the factor $r^{k-q}$ yield
\begin{equation}\label{multi:abs-diff}
  (1-zb_{-1})r^{q+1} - \sum_{j=0}^{q} (a_j+zb_j)r^{q-j}=0.
\end{equation}
This equation is called the {\em characteristic equation}. For example, the characteristic equations for $3$-step AB, AM, and BDF methods
are
\begin{align*}
  & r^{3} - (1+\tfrac{23}{12}z)r^2 +\tfrac{16}{12}z r - \tfrac{5}{12}z = 0,\\
  & (1-\tfrac{9}{24}z)r^3-(1+\tfrac{19}{24}z)r + \tfrac{5}{24}zr-\tfrac{1}{24}z = 0, \\
  & (1-\tfrac{6}{11}z)r^3-\tfrac{18}{11}r^2 + \tfrac{9}{11}r-\tfrac{2}{11} = 0,
\end{align*}
respectively. 

Assume that the characteristic polynomial has roots $r_0,r_1,\ldots,r_q$. As the roots depend on $z$, we denote them by
$r_0(z),r_1(z),\ldots,r_q(z)$. The boundary of the stability region is where all roots have magnitude $1$ or less, and at least one root has magnitude $1$. Roots with magnitude $1$ can be represented as $r=e^{i\theta}$ for $0\leqslant \theta<2\pi$, where $i$ is the imaginary number. To obtain the boundary of the stability region we can find all $z$ where \eqref{multi:abs-diff} holds true with $r=e^{i\theta}$.
One can separate the terms containing $z$ and write \eqref{multi:abs-diff} in the form
\begin{equation*}
  \rho(r)-z\sigma(r) = 0
\end{equation*}
where
$$
\rho(r) = r^{q+1}-\sum_{j=0}^{q}a_j r^{q-j},\quad \sigma(r) =b_{-1}r^{q+1}+ \sum_{j=0}^{q}b_jr^{q-j}.
\vsp 
$$
Now, for all $\theta\in[0,2\pi)$ we compute the complex values 
$$
z = \frac{\rho(e^{i\theta})}{\sigma(e^{i\theta})} =: u+iv
$$
and plot $v$ against $u$. The plot includes the boundary of the stability region. With a little care we can identify the true stability region. 

\begin{figure}[ht!]
  \centering
  \includegraphics[scale = .05]{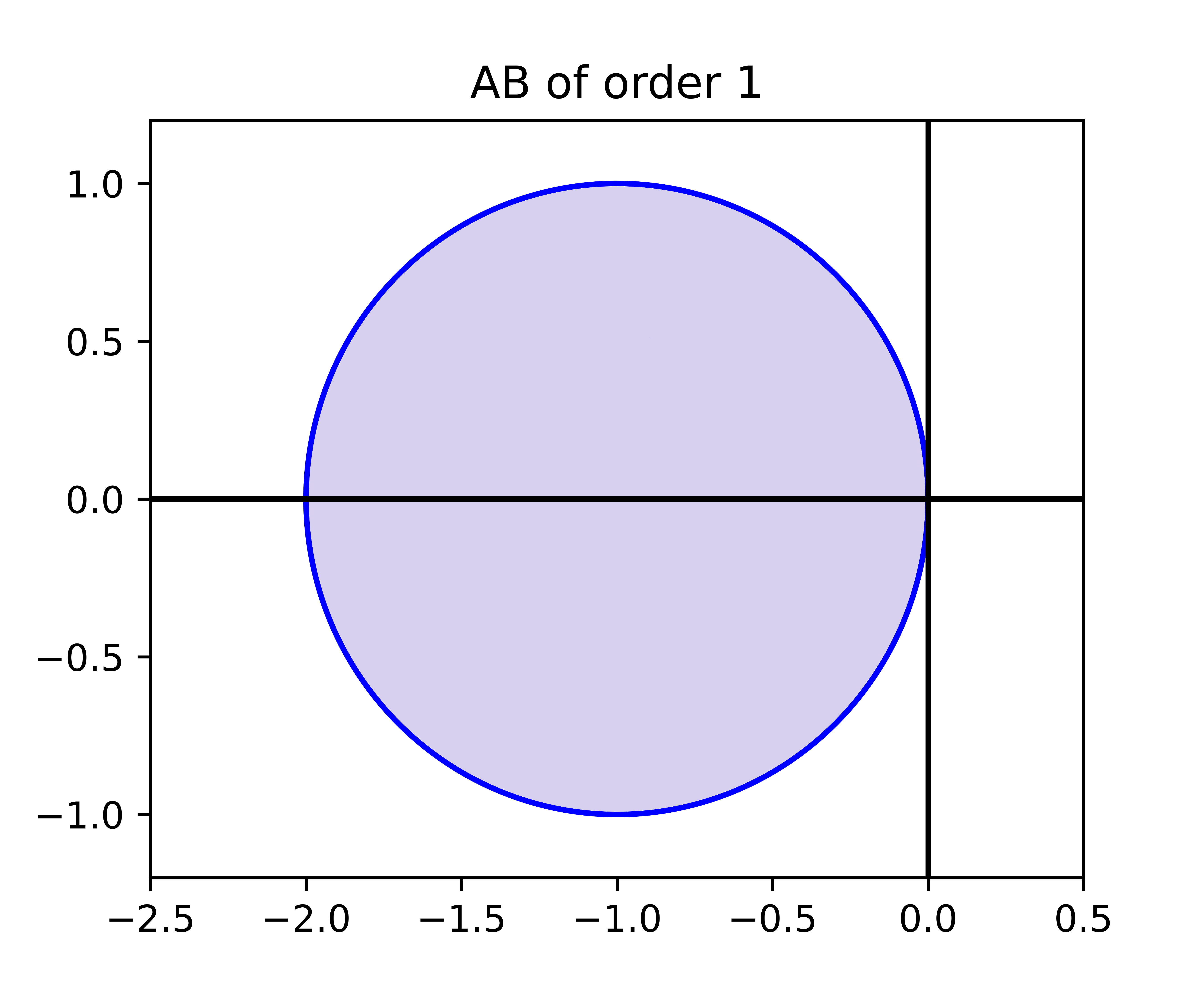}\includegraphics[scale = .05]{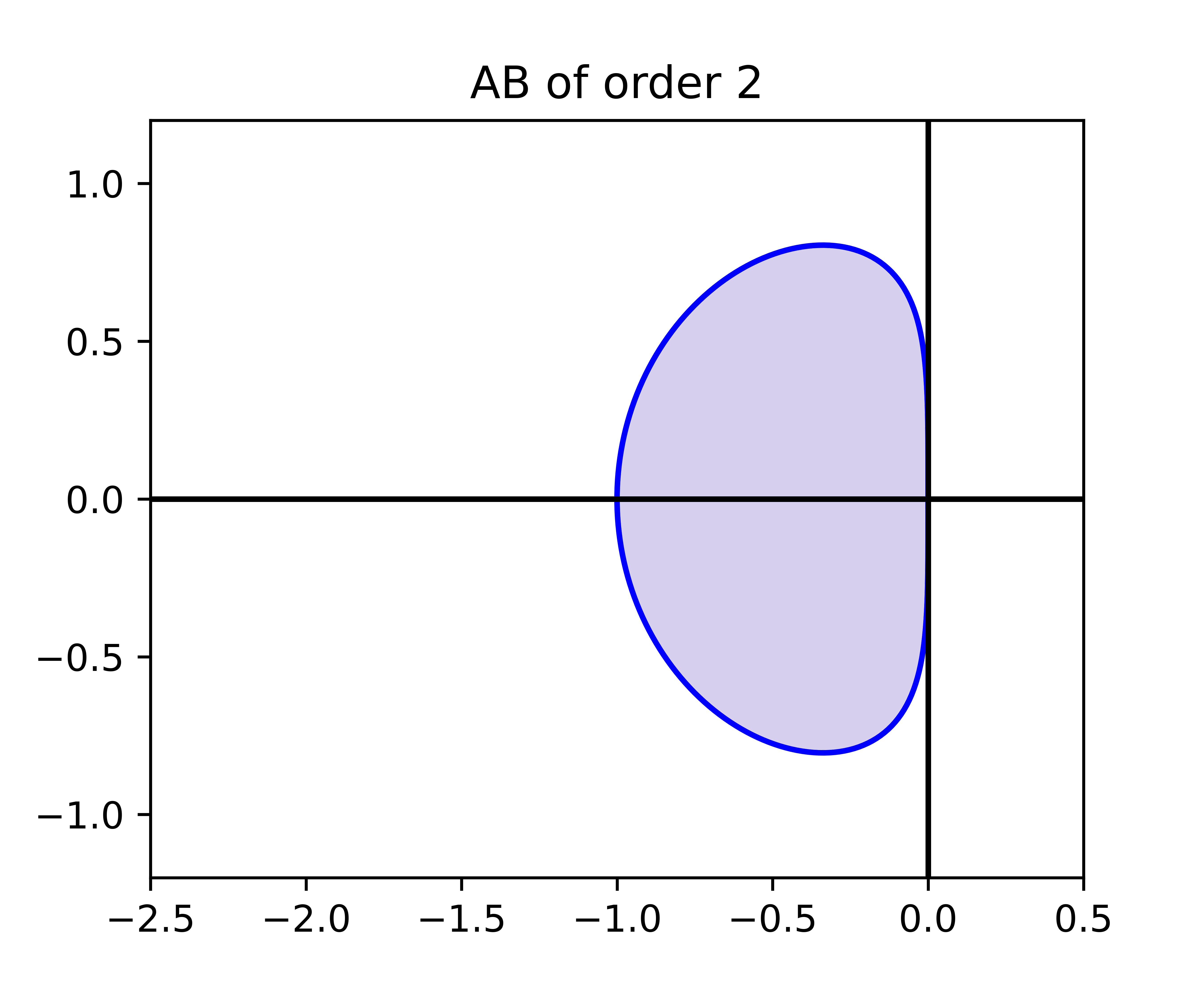}
  \includegraphics[scale = .05]{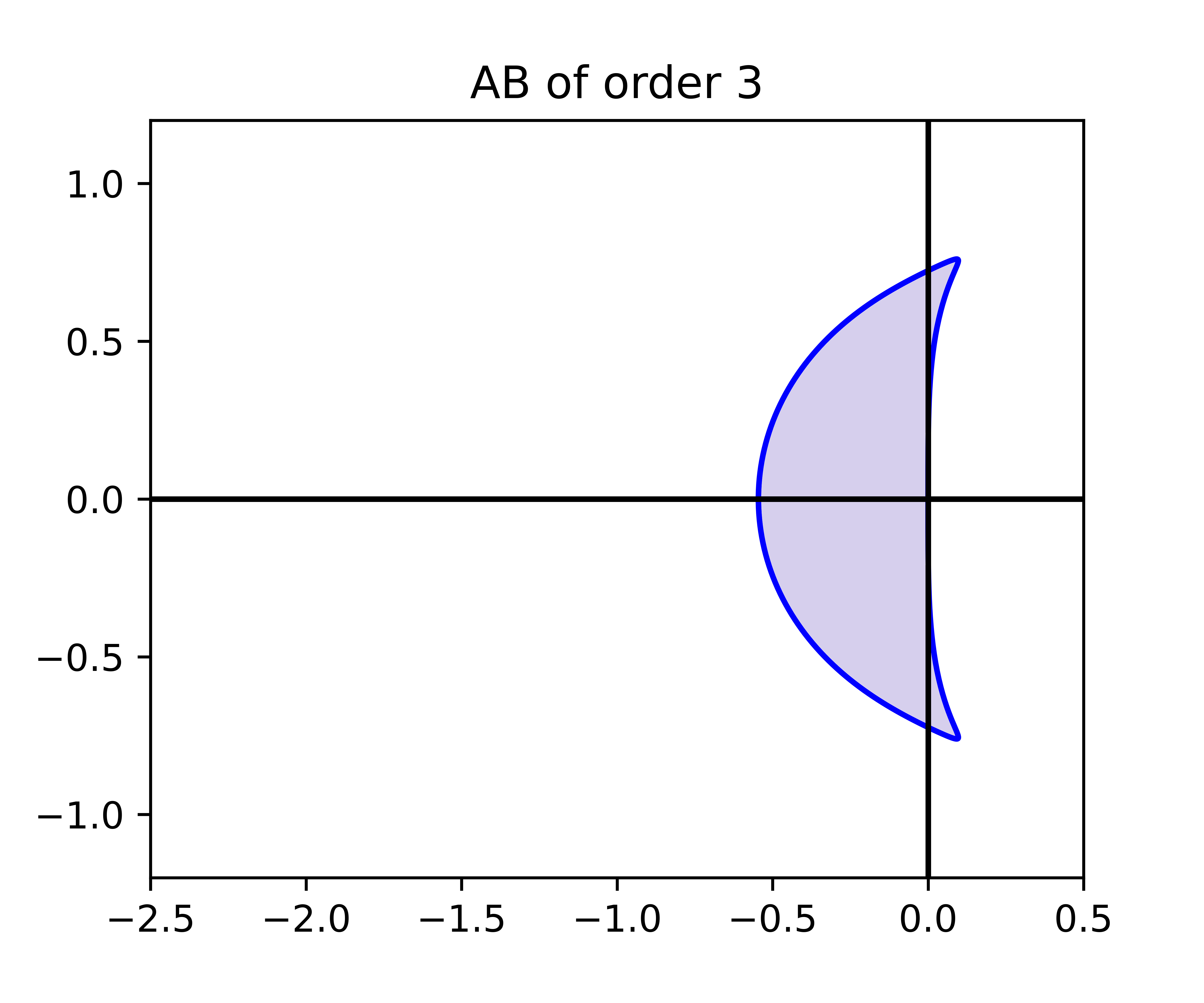}\includegraphics[scale = .05]{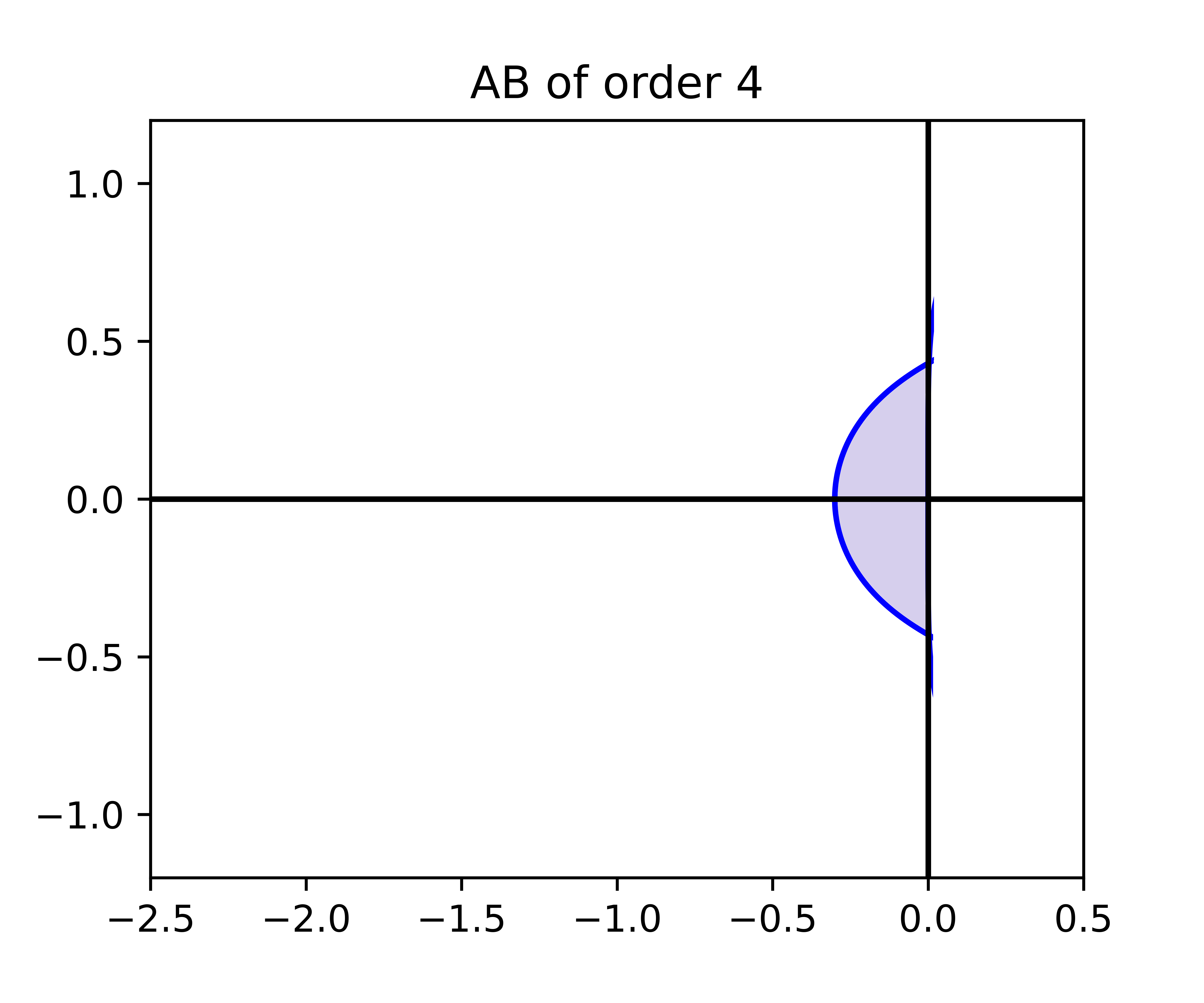}
  \caption{Stability regions of Adams-Bashforth methods of orders $1$ to $4$. }\label{fig:AB-stab-reg}
\end{figure}

\begin{figure}
  \centering
  \includegraphics[scale = .05]{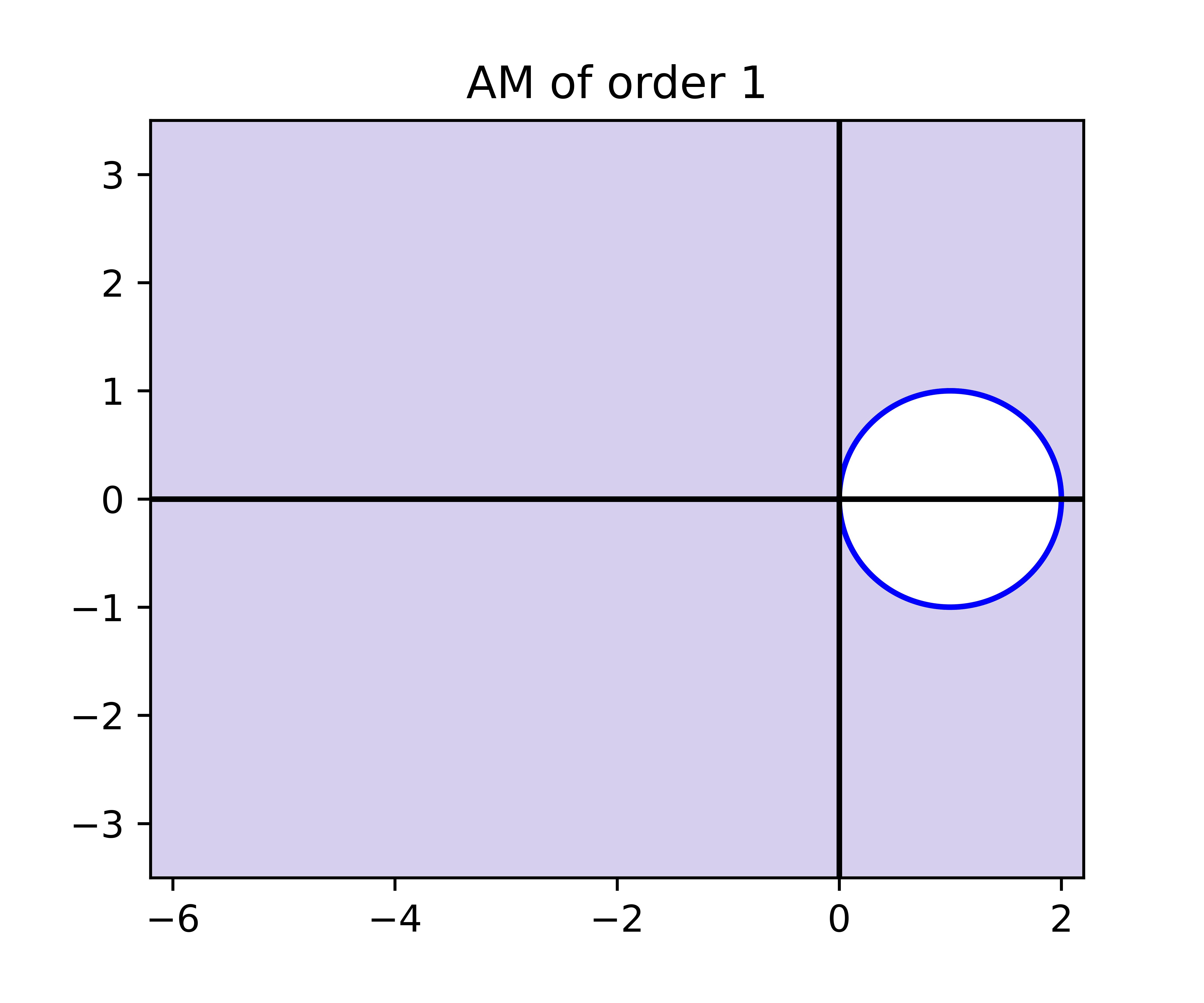}\includegraphics[scale = .05]{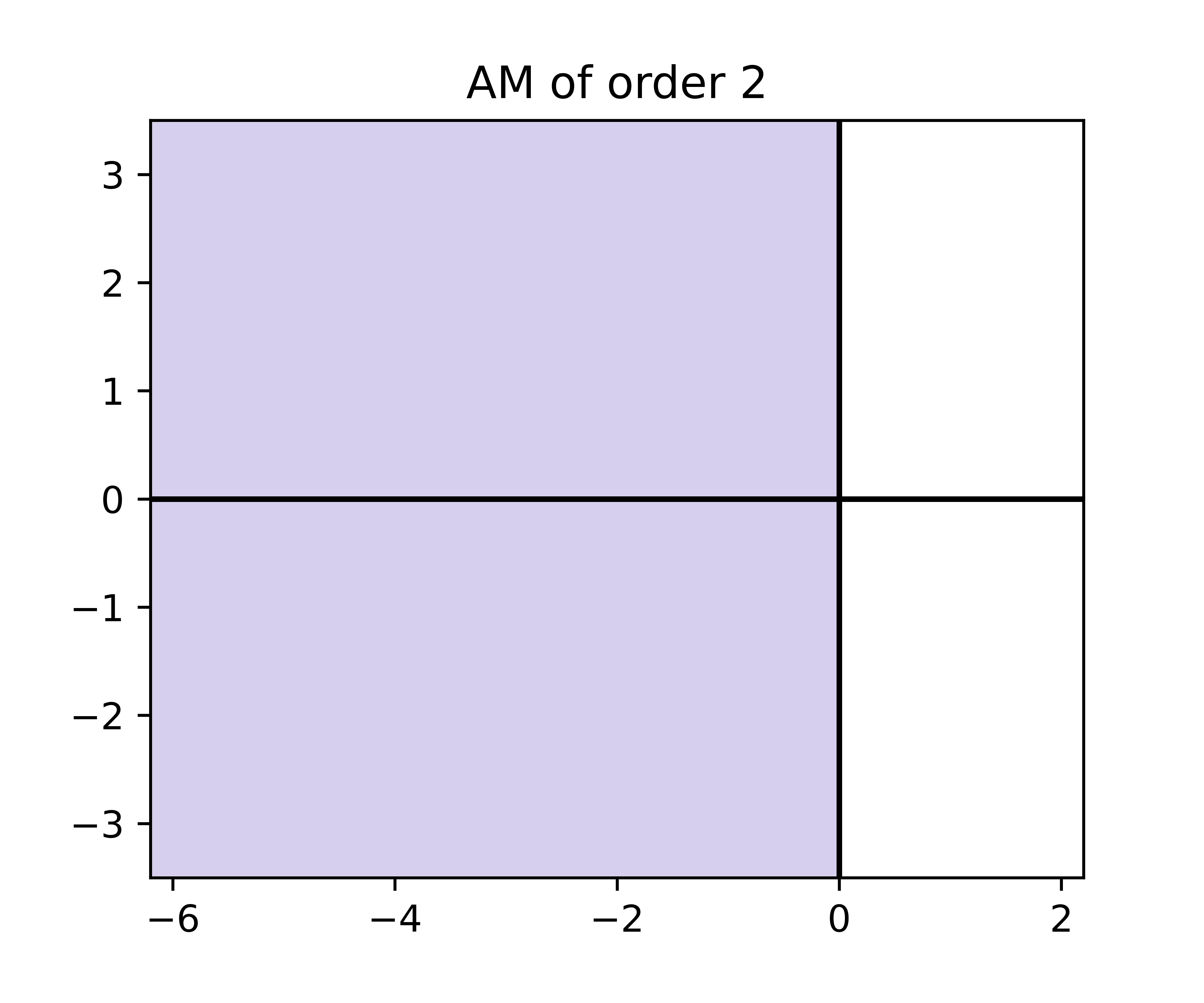}
  \includegraphics[scale = .05]{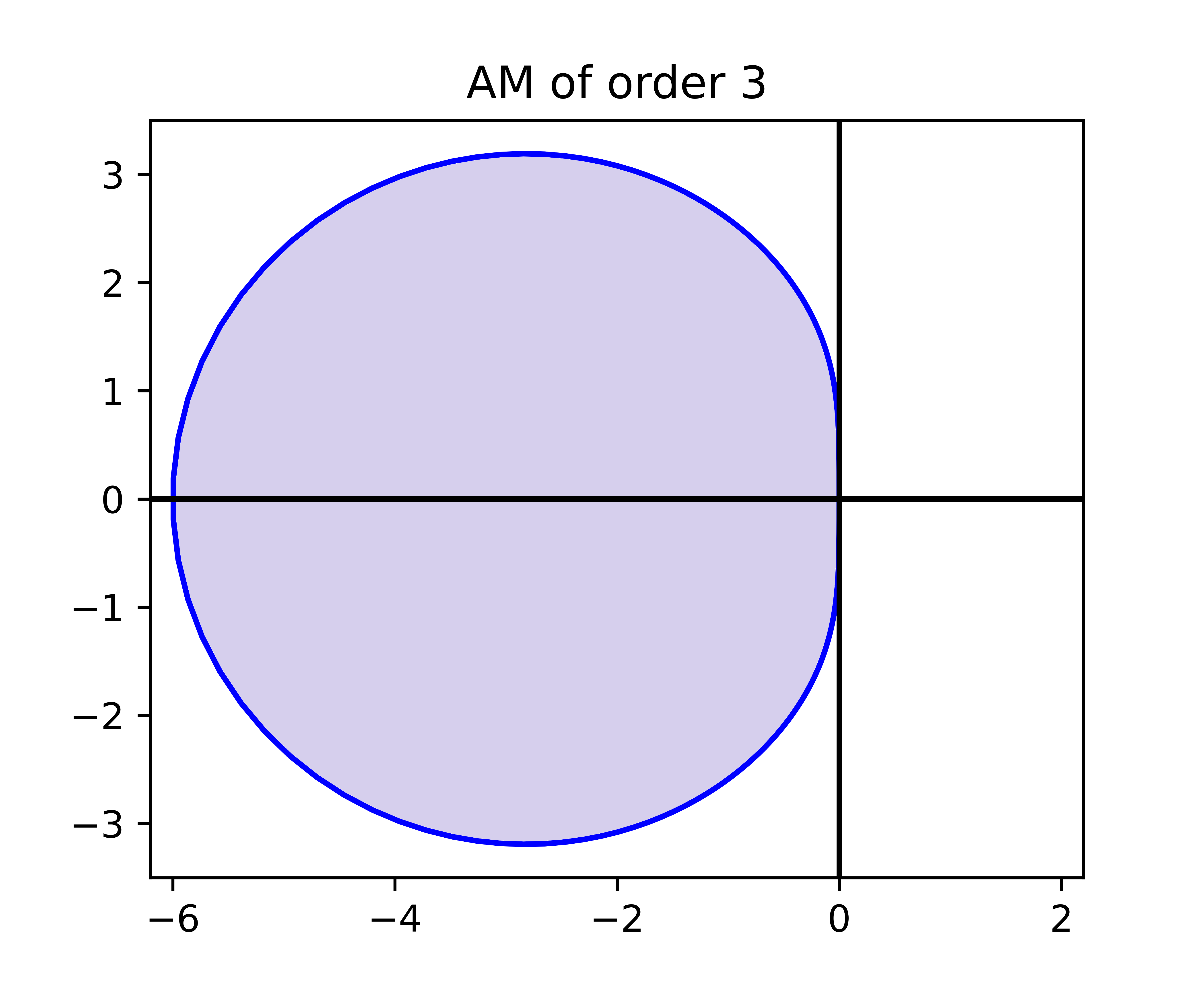}\includegraphics[scale = .05]{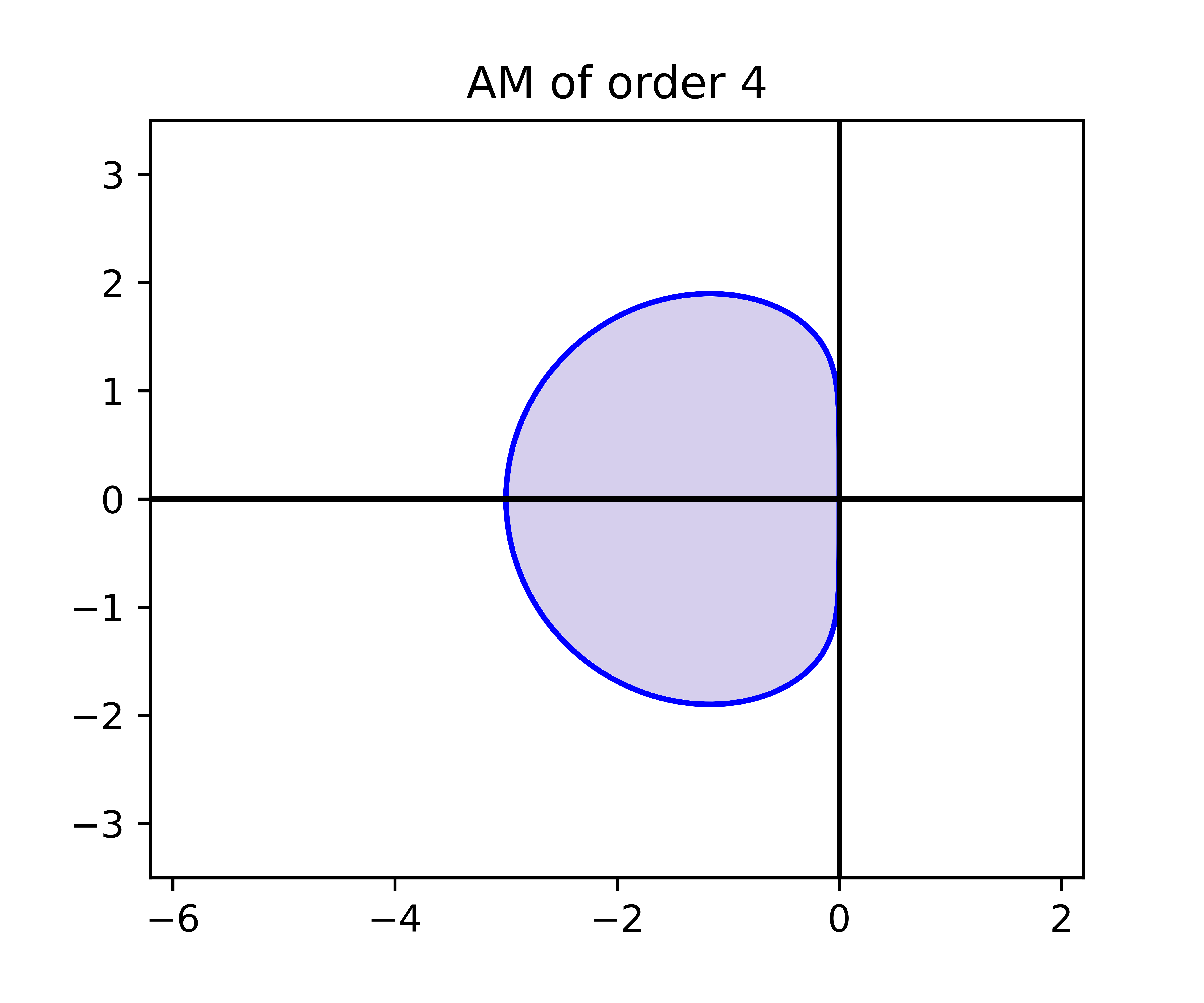}
  \caption{Stability regions of Adams-Moulton methods of orders $1$ to $4$. Note the different scale on the axes comapred to Figure \ref{fig:AB-stab-reg}.}\label{fig:AM-stab-reg}
\end{figure}

Stability regions of the Adams-Bashforth and Adams-Moulton methods are plotted in Figures \ref{fig:AB-stab-reg} and \ref{fig:AM-stab-reg}.
It can be seen that for AB and AM methods of equal order, the AM method has a larger absolute stability region. Unless the AM methods of order $1$ and $2$ which are the implicit Euler's and the trapezoidal methods, stability regions of other AB and AM methods do not contain the left-half plane and even do not encompass the whole negative real line. Consequently, these methods are not adequate for solving stiff differential equations.

\begin{figure}[ht!]
  \centering
  \includegraphics[scale = .05]{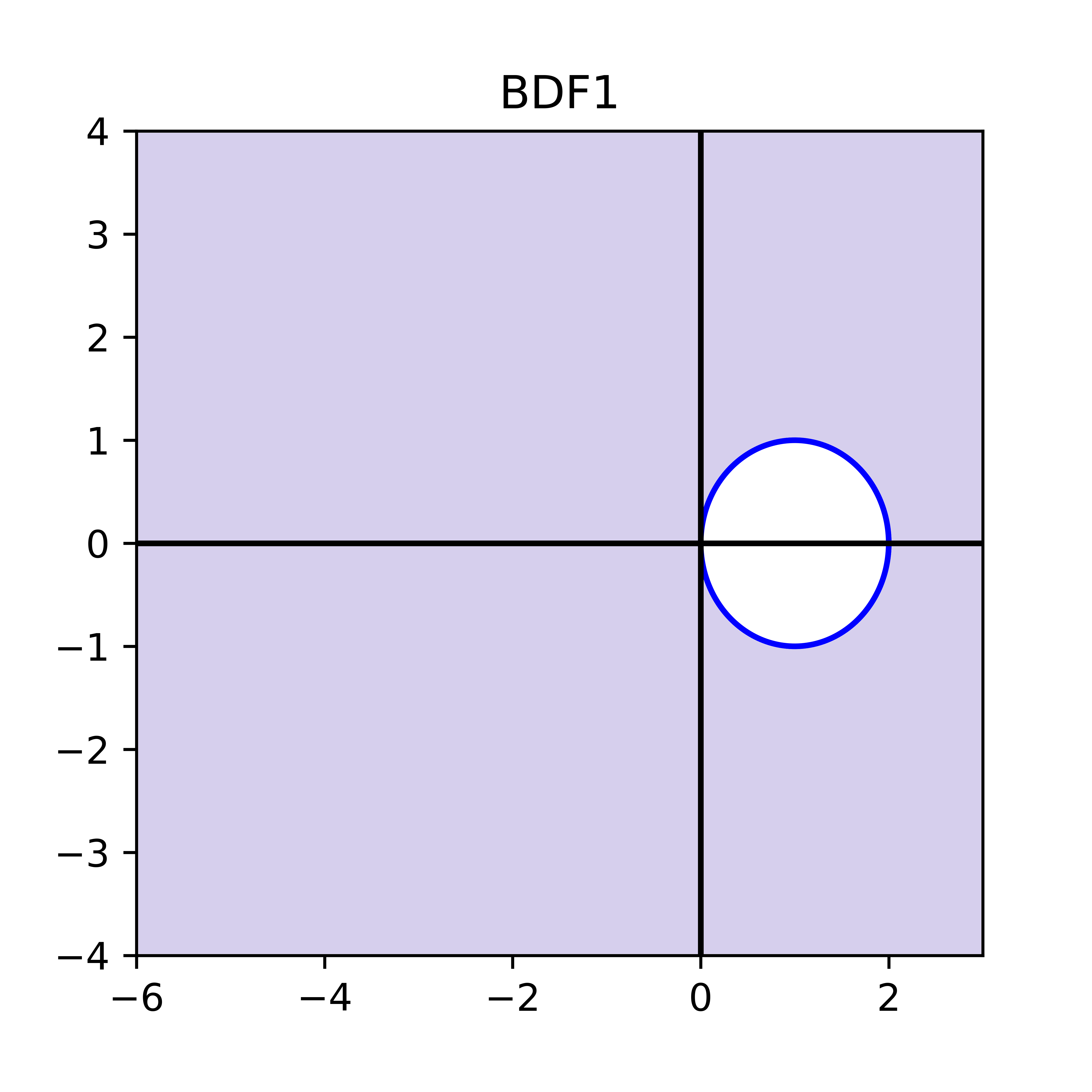}\includegraphics[scale = .05]{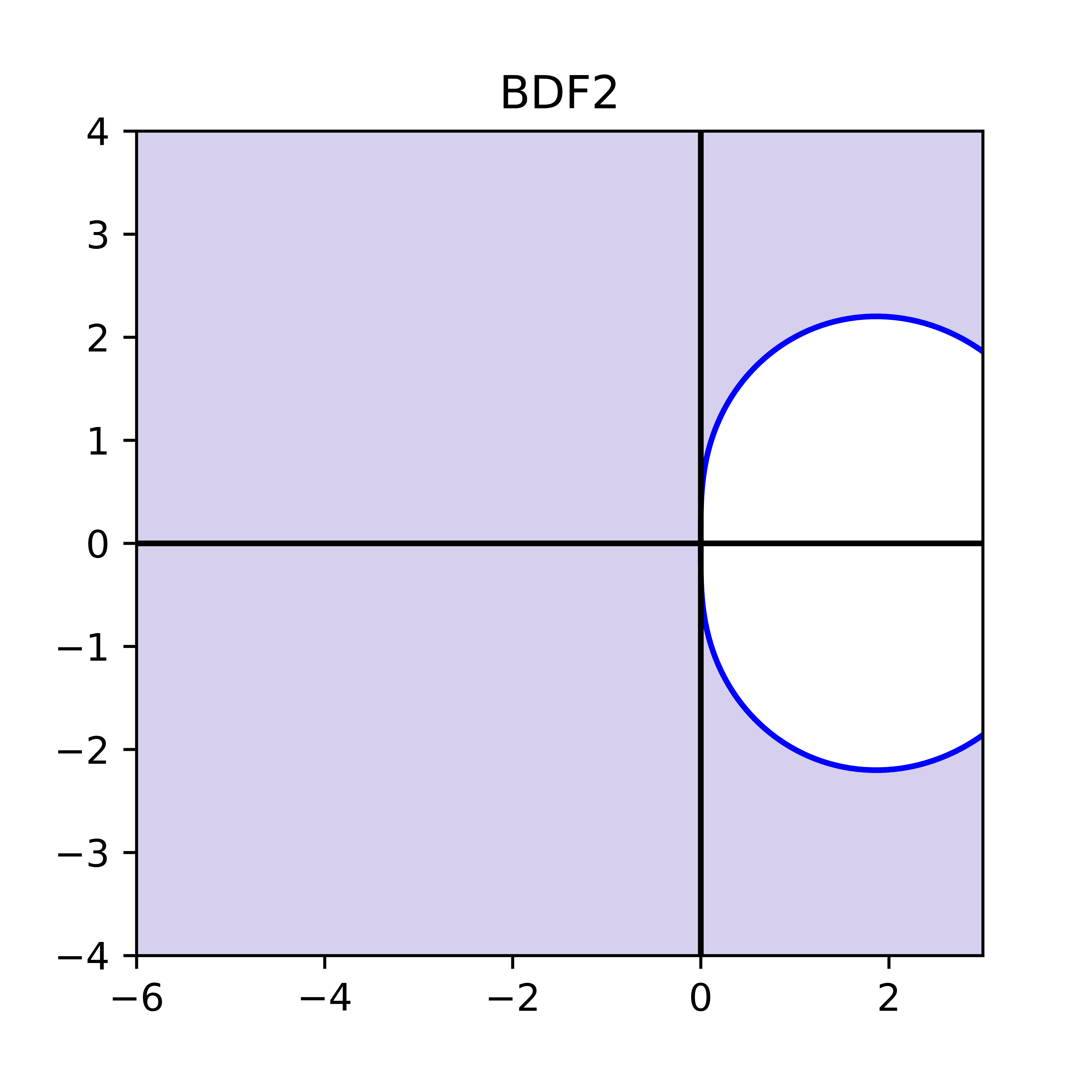}\includegraphics[scale = .05]{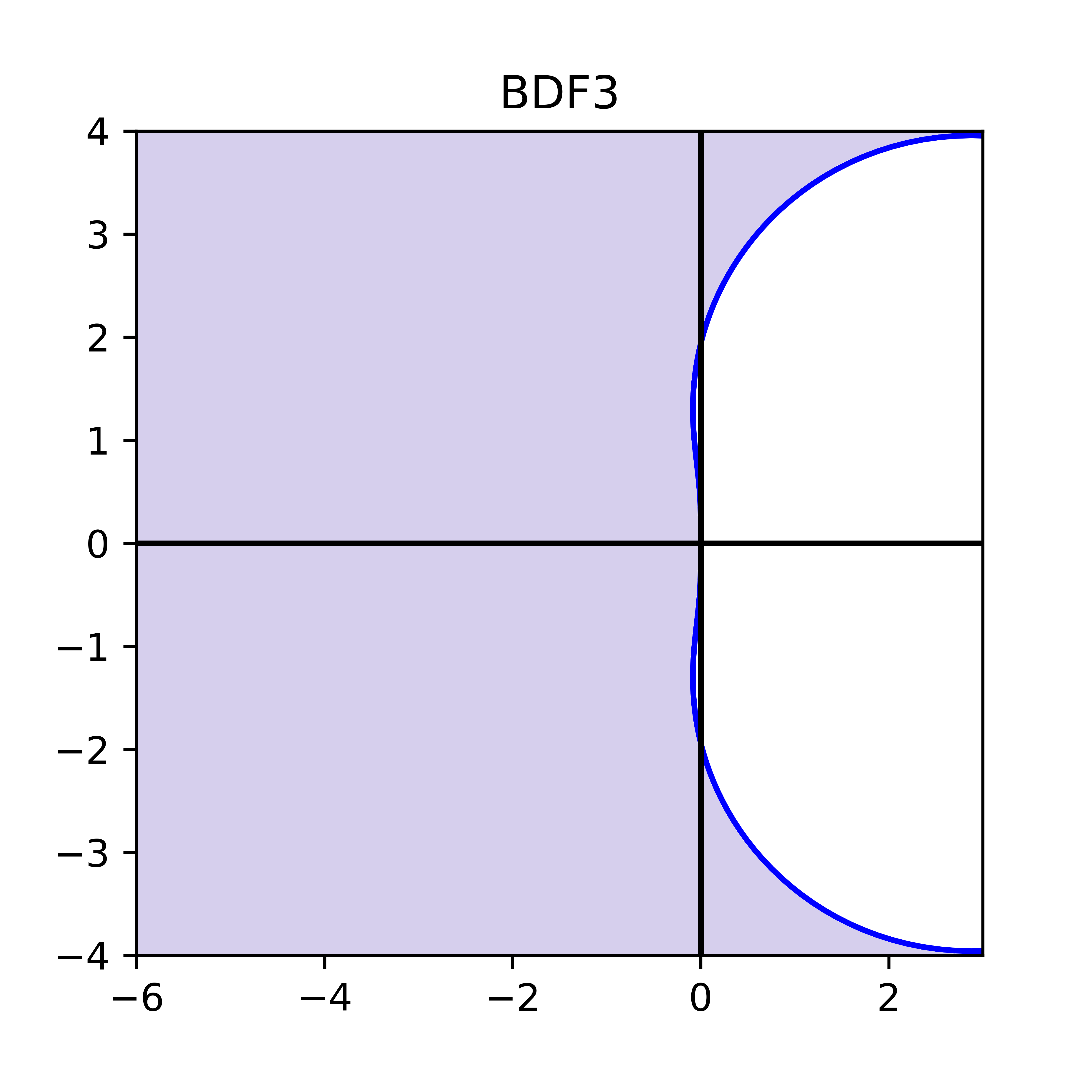}
  \includegraphics[scale = .05]{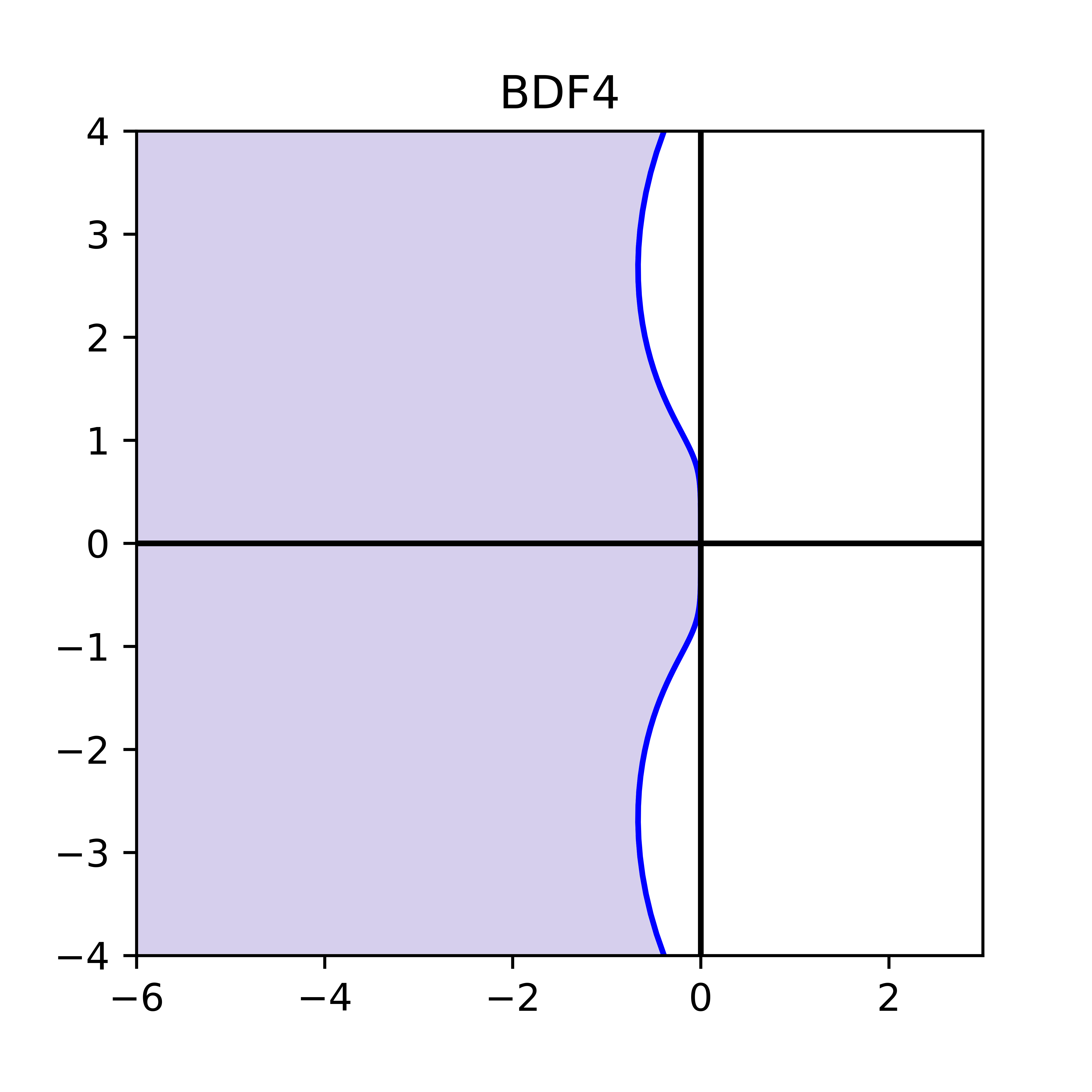}\includegraphics[scale = .05]{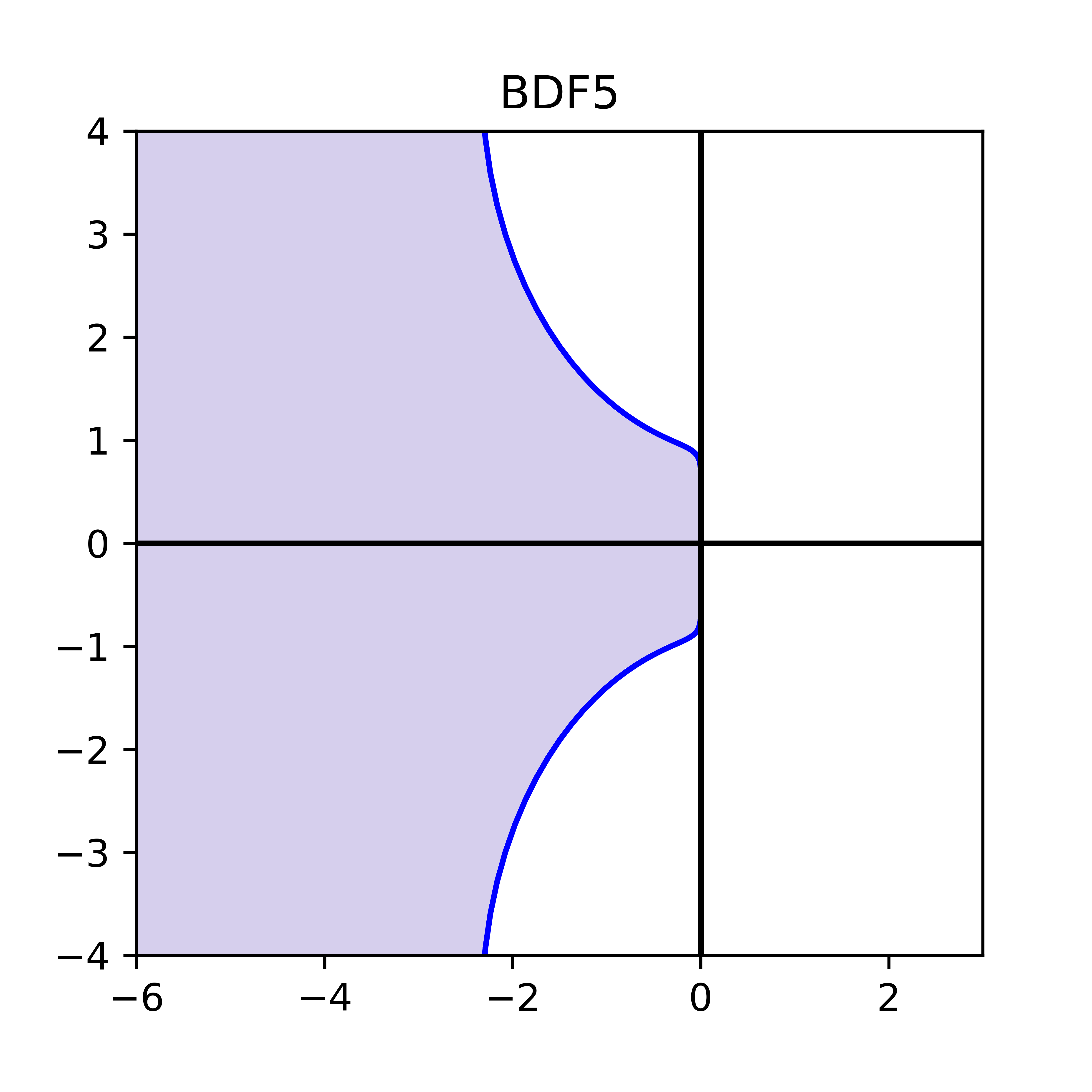}\includegraphics[scale = .05]{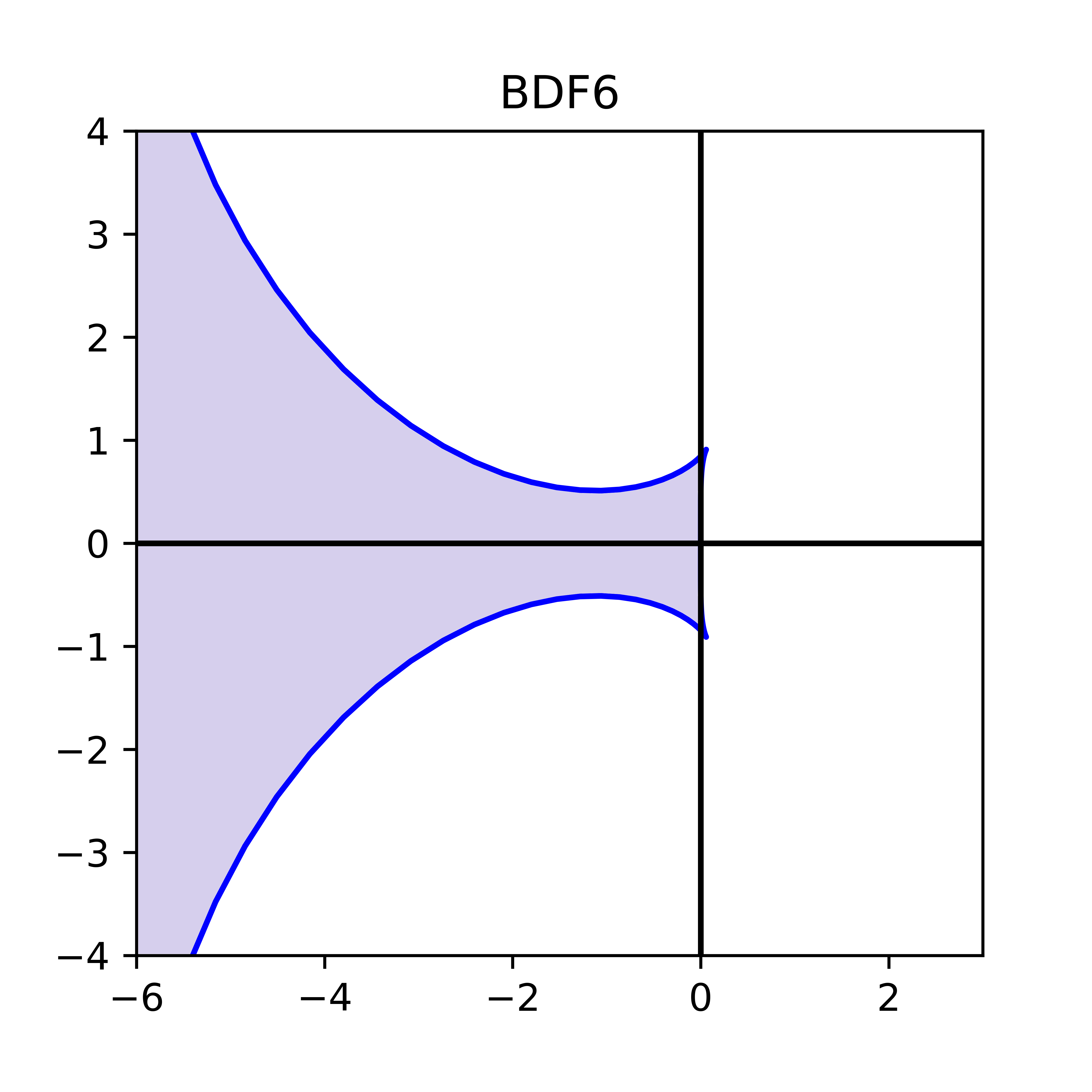}
  \caption{Stability regions of BDF methods of orders $1$ to $6$. BDF1 and BDF2 are $A$-stable and BDF3-BDF6 are $A(\alpha)$-stable. }
  \label{fig:BDF-stab-reg}
\end{figure}

The stability regions of BDF schemes are plotted in Figure \ref{fig:BDF-stab-reg}. Notably, BDF1 (implicit Euler's method) and BDF2 are $A$-stable. As the order of the method increases, the region of absolute stability diminishes, although it still encompasses the entire negative real line and maintains $A(\alpha)$-stability. The angle $\alpha$ decreases with increasing order from 1 to 6. Nonetheless, BDF1 to BDF6 are suitable schemes for solving stiff differential equations. However, schemes of orders $q\geqslant 7$ lose this property and are not adequate for solving stiff problems.
\vsp 

\subsection{One-step versus multistep methods}
Finally, we compare one-step and multistep methods for numerical solution of ODEs. 
Some advantages of one-step methods over multistep methods are listed below.
\begin{itemize}
\item One-step methods are self-starting: the available initial guess $y_0$ is enough to start and compute next values $y_1,y_2,\ldots$.
\item Adaptivity is easier with one-step methods. The step length  $h$ can be changed at any point in one-step method, based on an error estimate. In a multistep methods, however, more care is required since the previous values are assumed to be equally spaced in the standard form of these methods given in Table \ref{tb:ABmethods} and \ref{tb:AMmethods}.
\item If the solution $y(t)$ is not smooth at some isolated point $t^*$ then with a one-step method it is often possible to get
full accuracy simply by ensuring that $t^*$ is a grid point. This is impossible with multistep methods.
\end{itemize}
On the other hand, one-step methods have some disadvantages. The disadvantage of
Taylor series methods is that they require differentiating the given equation and are cumbersome
and often expensive to implement. RK methods only use evaluations of
the function $f$, but a higher order RK method requires evaluating $f$ several times
each time step. While, in multistep methods only one new
$f$ evaluation is required in each time step.

%% file: lec3_part4.tex
\section{MATLAB's ODE suite}\label{sect:matlab}
A suite of ODE solvers was introduced with version 5 of MATLAB and continued to current versions. The suite now contains 
seven solvers.
Here we introduce some of them.
\begin{itemize}
\item \texttt{ode23}

This is a low-order adaptive embedded solver for {\em nonstiff} ODEs. 
The `23' in the function name indicates that two simultaneous single-step formulas, one of second order and one of third order, are involved.

\item \texttt{ode45}

This is a high order embedded Runge-Kutta-Fehlberg solver which uses a
fourth-order and fifth-order formulas. For differential equations with smooth solutions, \texttt{ode45} is often more accurate than \texttt{ode23} but still works for {\em nonstiff} ODEs.  
The MATLAB documentation recommends \texttt{ode45} as the first choice and Simulink blocks set \texttt{ode45} as the default solver.

\item \texttt{ode23t}

This solver is an implementation of the trapezoidal rule using a free interpolant. Use this solver if the problem is only {\em moderately stiff} and you need a solution without numerical damping. This solver can solve differential algebraic equations (DAEs) as well.  

\item \texttt{ode15s}

This solver is suitable for solving {\em stiff} ODEs and DAEs. 
It is a variable order solver based on the numerical differentiation formulas. Optionally, it uses the backward differentiation formulas (BDFs, also known as Gear's method) which are a class of multistep solvers. (we did not learn multistep methods  in this course). Try \texttt{ode15s} when \texttt{ode45} fails, or is very inefficient. 

\item See also \texttt{ode113}, \texttt{ode78}, \texttt{ode89}, etc. in the MATLAB's help document.
\end{itemize}
You can explore other ODE and DAE solvers by referring to MATLAB's help documentation or other sources available.
The syntax for calling these solvers is one of the 
\begin{lstlisting}[style = matlab]
[t,y] = solver(odefun,tspan,y0)
[t,y] = solver(odefun,tspan,y0,options)
\end{lstlisting}
where \texttt{solver} is one of \texttt{ode45}, \texttt{ode23} and else. 
\texttt{odefun}
is a function that evaluates $f(t,y)$, the right-hand side of the differential equations. 
\texttt{y0} is the initial condition. 
\texttt{tspan}
is a two-vector specifying the interval of integration, $[t_0,t_{final}]$. 
To obtain solutions at specific times (all increasing or all decreasing), use 
$\texttt{tspan} = [t_0,t_1,\ldots,t_{final}]$.
For example the command
\begin{lstlisting}[style = matlab]
[t,y] = ode45(odefun,[0:0.01:0.5],[0 1]);
\end{lstlisting}
returns the function values at the specified vector $\texttt{[0:0.01:0.5]}$. But if the values at specified points are not required you can simply set
$\texttt{tspan = [0 0.5]}$. 

Optional integration arguments are created using the \texttt{odeset} function.
For example we can customize the error tolerances using 
\begin{lstlisting}[style = matlab]
options = odeset('RelTol',1e-6,'AbsTol',1e-10);
[t,y] = ode45(@myfunc,[0 0.5],[0 1],options);
\end{lstlisting}
This guarantees that the error at each step is less than
\texttt{RelTol} times the value at that step, and less than \texttt{AbsTol}.
More precisely 
$$
|e_k|\leqslant \max\{\texttt{RelTol}\times |y_k|, \texttt{AbsTol}\}. 
$$ 
Decreasing error tolerance can considerably slow the solver.
\vsp 
\begin{example}
To solve the scaler value equation
 $$y'(t) = t^3/y(t), \quad 0\leqslant t\leqslant 10, \quad y(0)=1$$ 
 using \texttt{ode23} we can write (without additional options):
\begin{lstlisting}[style = matlab1]
[t,y] = ode23(@(t,y) t^3/y,[0 10],1);
plot(t,y,'-ob')
\end{lstlisting}
The second example is a known nonstiff system of equations describing the motion of a rigid body without external forces:
\begin{align*}
y_1'(t) &= y_2(t)y_3(t)\\
y_2'(t) &= -y_1(t)y_3(t) \\
y_3'(t) & = -0.51 y_1(t)y_2(t)
\end{align*}
with initial conditions $y_1(0)=0$, $y_2(0)=1$ and $y_3(0)=1$. The time interval is $[0,12]$. With optional relative and absolute errors, we can compute the solutions by writing the following script.
\begin{lstlisting}[style = matlab1]
options = odeset('RelTol',1e-4,'AbsTol',[1e-4 1e-4 1e-5]);
[t,y] = ode45(@rigid,[0 12],[0 1 1],options);
plot(t,y(:,1),'-',t,y(:,2),'-.',t,y(:,3),'.');
\end{lstlisting}
The \texttt{rigid} function can be written in a separated script as a new function:
\begin{lstlisting}[style = matlab1]
function yprime = rigid(t,y)
yprime = zeros(3,1);   
yprime(1) = y(2)*y(3);
yprime(2) = -y(1)*y(3);
yprime(3) = -0.51*y(1)*y(2);
end
\end{lstlisting}
\end{example}
\vsp 
All solvers solve systems of equations in the form $\y'= \f(t,\y)$ or problems that involve a mass matrix, $M(t,\y)\y' = \f(t,\y)$. \texttt{ode15s} and \texttt{ode23t} can solve problems with a mass matrix that is singular, i.e., DAEs. The mass matrix can be imported via \texttt{odeset}. 
See doc \texttt{odeset} in the MATLAB's help for a list of options you can customize. 
\vsp
\begin{labexercise}
Repeat solving examples in Lab Exercise \ref{py_ex_euler} using \texttt{ode23} and \texttt{ode45}. Try to test different options. In each case plot the exact and numerical solutions in the same figure. 
\end{labexercise}
\vsp 

\begin{labexercise}
The Van der Pol oscillator is a self-maintained electrical circuit comprised of an inductor ($L$), an initially charged capacitor with capacitance ($C$), and a nonlinear resistor ($R$), all connected in series, as shown in the figure below.

\begin{center}
\includegraphics[scale=.85]{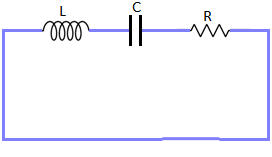}
\end{center}

By using the operational amplifier, the characteristic intensity-tension of the nonlinear resistance ($R$) is given as
\begin{equation}\label{Intensity}
U_R = -R_0 i_0 \left[\frac{i}{i_0}-\frac{1}{3}\left(\frac{i}{i_0}\right)^3 \right] \vsp
\end{equation}
where $i$ is the current and $i_0$ and $R_0$ are the current and the resistance of the normalization respectively.  By applying the Kirchhoff's law to the above figure we have 
\begin{equation}\label{law}
U_L + U_R + U_C = 0 
\end{equation}
where $U_L$ and $U_C$ are the tension to the limits of the inductor and capacitor, respectively, and are defined as
\begin{equation}\label{UL_UC}
U_L = L \frac{\mathrm{d}i}{\mathrm{d}t}, \quad U_C = \frac{1}{C} \int i\, \mathrm{d} t. \vsp 
\end{equation}
\begin{enumerate}
\item  Substitute \eqref{Intensity} and \eqref{UL_UC} into \eqref{law} and obtain an integral-differential equation. Then differentiate it and obtain the following second order ODE:
\begin{equation}\label{second_order1}
L\frac{\mathrm{d}^2i}{\mathrm{d}\tau^2} - R_0\left[1-\frac{i^2}{i_0^2} \right]
\frac{\mathrm{d} i}{\mathrm{d} \tau } + \frac{i}{C} = 0.
\end{equation}
Then use the change of variables $u = i/i_0$ and $t = \omega \tau$ where $\omega = 1/\sqrt{LC}$ (electric pulsation), to convert \eqref{second_order1} to the following equation:
\begin{equation*}
\frac{\mathrm{d}^2u}{\mathrm{d}t^2} - R_0\sqrt{\frac{C}{L}}\left(1-u^2 \right) 
\frac{\mathrm{d} u}{\mathrm{d} t} + u = 0.
\end{equation*}
\item
By setting $\mu = R_0\sqrt{\frac{C}{L}}$ and adding initial conditions obtain the well-known 
 {\em Van der Pol} equation
 \begin{equation}\label{vdp}
\begin{split}
&u''-\mu u'(1-u^2)+u = 0, \quad 0\leqslant t\leqslant b\\
&u(0) = u_0, \quad u'(0)=u'_0.
\end{split}
\end{equation}
and convert it 
to a system of first-order ODEs. 
\item 
 Apply RK2 and RK4 methods (your own codes) for solving the system obtained from equation \eqref{vdp} with different values $\mu = 1,10,100,1000$, and with initial conditions $u_0=2$ and $u'_0=0$. 
 For $\mu=1$ let $b=20$, for $\mu=10$ let $b=100$, for $\mu=100$ let $b = 500$ and for $\mu = 1000$ let $b = 5000$. For large values of $\mu$ use a supper small stepsize $h$
to get accurate results.
In each case plot the numerical solution $u$ and $u'$ in terms of $t$ in interval $[0,b]$ and report your observations.  Also report the executing times (sec.) in a table. 
 
\item  Now use ODE solvers \texttt{ode45}, \texttt{ode23t} and \texttt{ode15s} to solve this ODE again with different values of $\mu$ and $b$ given in item (1). In each case produce the plot of $u$ and $u'$ in terms of $t$, and compute the CPU time required (sec.) for executing the codes. Compare with the results of item (2).

\end{enumerate}

\end{labexercise}
\vsp

%% file: appendix.tex
\section{Appendix}
\subsection*{A: Difference equations}
Consider the following {\em homogeneous difference equation} of order $n$,
\begin{equation}\label{diff-eq_n}
  c_{k}y_{k} + c_{k-1}y_{k-1} + \cdots + c_{k-p} y_{k-p} = 0, \quad k\geqslant p
\end{equation}
with given initial values $y_0,y_1,\ldots,y_{p-1}$. The general solution of this equation is obtained by looking for solutions of the special form
$$
y_k =  r^k, \quad k\geqslant 0.
$$
If we can obtain $p$ linearly independent solutions, then an arbitrary linear combinations of this solutions give the general solution of \eqref{diff-eq_n}. Substituting $y_k=r^k$ into \eqref{diff-eq_n} and cancelling $r^{k-p}$, we obtain
\begin{equation}\label{characteristic_eq}
c_{p}r^{p} + c_{p-1}r^{p-1} + \cdots + c_1r+ c_0 = 0 \vsp 
\end{equation}
which is called {\em characteristic equation}, and the left side is {\em characteristic polynomial}. If \eqref{characteristic_eq} possesses $p$ distinct solutions (roots) $r_1,r_2,\ldots,r_p$ then the general solution of \eqref{diff-eq_n} is
\begin{equation}\label{differ:generalsol}
  y_k = \sum_{j=1}^{p} \beta_j r_j^k, \quad k\geqslant 0.
\end{equation}
The coefficients $\beta_0, \beta_1,\ldots, \beta_{p-1}$ are obtained by imposing the known initial values
$y_0,y_1,\ldots,y_{p-1}$ into the general solution \eqref{differ:generalsol}.
It is clear that if $|r_j|\leqslant 1$ then the solution $y_k$ does not grow as $k\to\infty$. 
However, if $r_j$ is a root of multiplicity $\nu$, i.e.,
$$
r_j = r_{j+1} = \cdots = r_{j+\nu-1},
$$
then the $\nu$ linearly independent solutions corresponding to these roots are
$$
r_j^k ,\, kr_j^k, \ldots, k^{\nu-1}r_j^k
$$
and in formula \eqref{differ:generalsol} the part $\beta_jr_j^k + \beta_{j+1}r_{j+1}^k+\cdots+\beta_{j+\nu-1}r_{j+\nu-1}^k$
should be replaced by
\begin{equation*}
  [\beta_j + k\beta_{j+1}+\cdots + k^{\nu-1}\beta_{j+\nu-1}]r_j^k.
\end{equation*}
In this case the solution $y_k$ remains stable as $k\to\infty$ provided that $|r_j|\leqslant 1$ for simple roots and $|r_j|<1$ for roots with multiplicity. 
 \vsp 
 
\begin{example}
The general solution of difference equation
$$
y_{k}+5y_{k-1}+6r_{k-2}=0
$$
is obtained in terms of roots $r_1=2$ and $r_2=3$  of characteristic polynomial
$r^2+3r+2$, 
$$
y_k = \beta_1 2^k  + \beta_2 3^k.
$$
If $y_0=0$ and $y_1=2$ are given then by solving $y_0=\beta_1+\beta_2$ and $y_1=2\beta_1+3\beta_2$ we simply obtain
$\beta_1=-2$ and $\beta_2=2$. The solution then is
$$
y_k = -2\times 2^k +2 \times 3^k,\quad k=0,1,2,3,\ldots.
$$
As second example consider the following difference equation
$$
y_{k}-5y_{k-1}+6y_{k-2}+4y_{k-3}-8y_{k-4}=0.
$$
The characteristic equation is $r^4 - 5r^3 + 6r^2 + 4r - 8=0$ with roots $r_1=-1$ and $r_2=r_3=r_4=2$. The general solution then is
$$
y_k = \beta_1(-1)^k + [\beta_2+k\beta_3+k^2\beta_4]2^k.
$$
With given initial values $y_0=-1$, $y_1=y_2=-7$ and $y_3=7$  we must solve 
\begin{align*}
  -1= & \beta_1+\beta_2 \\
  -7= & -\beta_1+2[\beta_2+\beta_3+\beta_4] \\
  7 =& \beta_1+4[\beta_2+2\beta_3+4\beta_4] \\
  7 = & -\beta_1+8[\beta_2+3\beta_3+9\beta_4]
\end{align*}
to obtain coefficients $\beta_j$. Solving this system gives $\beta_1=1,\beta_2=-2,\beta_3=-2$ and $\beta_4=1$.

\end{example}
\ \vsp